\documentclass{article}
\usepackage[english]{babel}
\usepackage{todonotes}
\usepackage{listings}
\usepackage{color} %red, green, blue, yellow, cyan, magenta, black, white
\usepackage{xcolor}
\definecolor{mygreen}{RGB}{28,172,0} % color values Red, Green, Blue
\definecolor{mylilas}{RGB}{170,55,241}
\usepackage{graphics} %new
\usepackage{graphicx}
\usepackage{caption}
\usepackage{subcaption}
\usepackage{amsmath}
\usepackage{amssymb}
\usepackage{mathtools}
\usepackage{overpic}
\usepackage{amsthm}
\usepackage{changepage}
\usepackage{hyperref}
\usepackage{comment}

\usepackage{algorithm} % algorithm package
\usepackage[noend]{algpseudocode}
\usepackage{float}      % suggested for floating Algorithms

\usepackage{enumerate}%

\usepackage{booktabs,array}

\usepackage{xcolor}

\usepackage{subcaption}
\usepackage{pgfplots}
\pgfplotsset{grid style={dotted,gray}}
\pgfplotsset{compat=newest}
\pgfplotsset{plot coordinates/math parser=false}
\newlength\figureheight
\newlength\figurewidth

\usepackage{stfloats} % for positioning of figure* on the same page
\usepackage{array, multirow}

\usepackage[pass]{geometry}
\usepackage[outdir=./]{epstopdf}
\usepackage{url}
\usepackage{tabstackengine}
\setstacktabbedgap{1.5ex}%               sets gap between columns
\setstackgap{L}{1.2\normalbaselineskip}% sets baselineskip of rows
\theoremstyle{definition}

\newtheorem{remark}{Remark}

\definecolor{KTHblue}{RGB}{25,105,188}
\definecolor{KTHlblue}{RGB}{22,159,219}
\definecolor{KTHyellow}{RGB}{251,186,0}
\definecolor{KTHred}{RGB}{176,9,48}
\definecolor{KTHlred}{RGB}{231,51,57}
\definecolor{KTHgreen}{RGB}{98,146,46}
\definecolor{KTHlgreen}{RGB}{175,202,11}
\definecolor{KTHpink}{RGB}{219,81,151}
\definecolor{usered}{RGB}{236, 91, 40}

\addto{\captionsenglish}{}
\renewcommand{\vec}[1]{\boldsymbol{#1}}

\addto{\captionsenglish}{}

\newcommand{\x}{\vec x}
\newcommand{\y}{\vec y}
\newcommand{\z}{\vec z}

\newcommand{\uu}{\vec u}
\newcommand{\bb}{\vec\beta}
\newcommand{\mm}{{\vec \mu}}
\newcommand{\be}{\begin{equation}}
\newcommand{\ee}{\end{equation}}
\newcommand{\ben}{\begin{enumerate}}
\newcommand{\een}{\end{enumerate}}

\newcommand{\bigO}{\mathcal O}
\newcommand{\mbf}[1]{{\mathbf #1}}
\newcommand{\pO}{{\partial\Omega}}    %  boundary
\newcommand{\R}{\mathbb{R}}

\newcommand{\partone}{i} %generic particle no, one
\newcommand{\parttwo}{k}  %generic particle no, two
\newcommand{\csi}{j} %coarse source index
\newcommand{\cci}{q} %coarse collocation index
\newgeometry{left=2.5cm,right =2.5cm,top= 2.5cm,bottom = 2.5cm}

\providecommand{\keywords}[1]{\textbf{Key words: } #1}

\usepackage[english]{nomencl}
\makenomenclature
\newcommand{\thickhline}{%
	\noalign {\ifnum 0=`}\fi \hrule height 1pt
	\futurelet \reserved@a \@xhline
}

\title{Preconditioning for near-contacts in large 2D Stokes flows: a locally compressed method of fundamental solutions}
\author{Anna Broms$^{*,1)}$, Anna-Karin Tornberg$^{2)}$, and Alex H. Barnett$^{3)}$
  \\\\ $^{1)}$ Department of Mathematics, Imperial College London, UK\\
    $^{2)}$ Department of Mathematics, KTH Royal Institute of Technology,
	Stockholm, Sweden\\	
  $^{3)}$ Center for Computational Mathematics, Flatiron Institute,
	New York, USA\\\\
	$^*$e-mail: a.broms@imperial.ac.uk}
\date{\today}

\begin{document}
	\lstset{language=Matlab,%
		%basicstyle=\color{red},
		%breaklines=true,%
		morekeywords={matlab2tikz},
		keywordstyle=\color{black},%
		morekeywords=[2]{1}, keywordstyle=[2]{\color{black}},
		identifierstyle=\color{black},%
		stringstyle=\color{mylilas},
		commentstyle=\color{mygreen},%
		showstringspaces=false,%without this there will be a symbol in the places where there is a space
		numbers=left,%
		numberstyle={\tiny \color{black}},% size of the numbers
		numbersep=9pt, % this defines how far the numbers are from the text
		emph=[1]{for,end,break},emphstyle=[1]\color{blue}, %some words to emphasise
		%emph=[2]{word1,word2}, 6 emphstyle=[2]{style},    
	}
	\maketitle

\begin{abstract}
  We tackle two key difficulties in
the simulation of the viscous hydrodynamics of a large dense collection of rigid particles:
%  simulating the relative motion of many small rigid particles moving close together in a viscous fluid:
%(i) the severe ill-conditioning of the discretized linear system as particle gaps shrink, and
%(ii) the issue of capturing lubrication-driven fine scales without resorting to a globally fine discretization.
(i) the poor convergence rate of an iterative solution of the discretized linear system as
particle gaps shrink, and
(ii) the large number of unknowns needed to accurately discretize the resulting lubrication-driven flows.
Our focus is the 2D Stokes resistance and mobility boundary value problems for nearly-touching disks.
To address both challenges, we introduce a general two-body preconditioning strategy,
% generalizing an idea of Cheng-Greengard.
and implement it with the method of fundamental solutions. %(MFS).
For each close particle pair, the hard-to-resolve interaction is represented in a basis
precomputed by solving a local boundary value problem on a fine grid.
In an iterative solve, the resulting flow field corrects that obtained from a coarse representation of all particles.
The local fine-grid correction can furthermore be compressed so that all particles except the pair itself are affected by an equivalent set of coarse sources.
Numerical experiments demonstrate rapid GMRES convergence in challenging multi-particle settings, with iteration counts remaining low even in densely packed suspensions. For example, the mobility problem is solved for a random close packing with area fraction $\varphi = 0.65$, $P = 10000$ monodisperse disks, and minimum separation $10^{-3}$, in 47 GMRES iterations, achieving five digits of accuracy with 72 vector unknowns per body.
% N fine, N coarse? 

%%%%%%%%%%%%%%%%%%%%%%%%%%% Clean up for submission
%   We tackle two key difficulties in
% the simulation of the viscous hydrodynamics of a large dense collection of rigid particles:
% (i) the poor convergence rate of an iterative solution of the discretized linear system as
% particle gaps shrink, and
% (ii) the large number of unknowns needed to accurately discretize the resulting lubrication-driven flows.
% Our focus is the 2D Stokes resistance and mobility boundary value problems for nearly-touching disks.
% To address both challenges, we introduce a general two-body preconditioning strategy,
% and implement it with the method of fundamental solutions. 
% For each close particle pair, the hard-to-resolve interaction is represented in a basis
% precomputed by solving a local boundary value problem on a fine grid.
% In an iterative solve, the resulting flow field corrects that obtained from a coarse representation of all particles.
% The local fine-grid correction can even be compressed so that all particles except the pair itself are affected by an equivalent set of coarse sources.
% Numerical experiments demonstrate rapid GMRES convergence in challenging multi-particle settings, with iteration counts remaining low even in densely packed suspensions. For example, the mobility problem is solved for a random close packing with area fraction $\varphi = 0.65$, $P = 10000$ monodisperse disks, and minimum separation $10^{-3}$, in just 47 GMRES iterations, achieving five digits of accuracy with 72 vector unknowns per body.

\end{abstract}
%\tableofcontents

\section{Introduction}\label{sec:intro}
    The Stokes equations---a set of linear elliptic PDEs---describe fluid flow in regimes where viscous forces dominate over inertia. This regime is typical for suspensions of small particles, from nanometer to micrometer in size, moving through a viscous fluid. At such scales, Stokes flow governs the motion of rigid bodies \cite{Broms2025,Bagge2021,crowder2025}, drops \cite{Wrobel2018,Sorgentone2018,Palsson2020}, vesicles \cite{Quaife2014, Barnett2015}, flexible fibers and filaments \cite{Schoeller2021,Maxian2022}, and even swimming microorganisms \cite{Swan2011,Usabiaga2021}.

%In this work, 
We focus on the exterior Stokes boundary value problem (BVP) in two dimensions, where the fluid domain is the unbounded region outside a collection of~$P$ rigid circular particles. Denoting the particles by $\Omega^{(i)} \subset \mathbb{R}^2$, with boundaries $\partial\Omega^{(i)}$, $i=1,\dots,P$, the union of all particles is $\Omega = \bigcup_{i = 1}^P \Omega^{(i)}$, and the fluid domain is $\mathbb{R}^2 \setminus \overline{\Omega}$. %, bounded by $\partial\Omega$. 
Many questions of physical interest can be addressed in this framework—for instance, \emph{how do the suspended particles affect the effective viscosity of the suspension?} Or, more generally, \emph{how do hydrodynamic interactions between nearby particles influence their collective motion?}

In 2D Stokes problems, of both resistance and mobility flavours,
with a fluid of constant dynamic viscosity \( \mu \), we have
\begin{equation}\label{stokeseq}
	\begin{aligned}
		-\mu\Delta \vec u + \nabla p &= \mathbf{0}, &&\text{ in }\mathbb R^2\setminus\overline{\Omega},\\
		\nabla\cdot \vec u &= 0, &&\text{ in }\mathbb R^2\setminus\overline{\Omega},\\
		\vec u &= \vec g^{(i)}, &&\text{ on }\partial\Omega^{(i)},\\
		\vec u(\vec x) &= \frac{\vec\Sigma}{4\pi \mu} \log \frac{1}{r} + \bigO(1), &&\text{ as } r \coloneqq \|\vec x\| \to \infty,
	\end{aligned}
\end{equation}
where \( \vec u \) is the velocity field and \( p \) the pressure.
The logarithmically unbounded term has been scaled so that
$\vec \Sigma $ is the total force on the fluid. % *** Anna: this is probably the standard notation, but $\Sigma$ is used also later on for the diagonal matrix of singular values.
There are no-slip conditions on the rigid particle boundaries. For particle~\( i \), the boundary velocity is given by
\begin{equation}\label{rigidbc}
	\vec g^{(i)}(\vec x) = \vec v^{(i)} + \omega^{(i)} (\vec x - \vec c^{(i)})^{\perp}, \quad \vec x \in \partial\Omega^{(i)},
\end{equation}
where \( \vec v^{(i)} \in \mathbb{R}^2 \) is a translational velocity, \( \omega^{(i)} \in \mathbb{R} \) an angular velocity, and \( \vec c^{(i)} \) the center of particle~\( i \). The perpendicular map \( (x_1, x_2)^{\perp} \coloneqq (-x_2, x_1) \) rotates a vector \( 90^\circ \) counterclockwise and represents the 2D analog of the 3D cross product between \(\vec e_3 \) and a vector in the plane.
For each particle, the hydrodynamic force and torque exerted by the fluid are given by integrals of the traction $\boldsymbol{\sigma} \vec n$ over the boundary:
\begin{equation}\label{eq:forces}
\begin{aligned}
	\vec f^{(i)} &= \int_{\partial \Omega^{(i)}} \boldsymbol{\sigma} \vec n \, ds,\\ %\label{force}\\
	t^{(i)} &= \int_{\partial \Omega^{(i)}} (\vec x - \vec c^{(i)})^\perp \cdot \boldsymbol{\sigma} \vec n \, ds, %\label{torque}
\end{aligned}
\end{equation}
where \( \boldsymbol{\sigma} := -p\mathbf{I} + \mu(\nabla \vec u + (\nabla \vec u)^T) \) is the Cauchy stress tensor, and \( \vec n \) is the unit normal pointing outward from each particle.
%and into the surrounding fluid.

% We will solve both the \emph{resistance problem} and its inverse, termed the \emph{mobility problem}. For resistance, the translational and angular velocities \( (\vec v^{(i)}, \omega^{(i)}) \) of each particle are prescribed, and the resulting hydrodynamic forces and torques \( (\vec f^{(i)}, t^{(i)}) \) are to be computed. In the mobility problem, the net force and torque on each particle are given, and one seeks the velocities that produce them. There is an important difference in these two settings that will impact how to solve them: In the resistance problem, \eqref{stokeseq} can be solved with the boundary conditions in \eqref{rigidbc} determined from the given velocities, and forces and torques are determined via \eqref{eq:forces} in a post-processing step. For mobility, to the contrary, given forces and torques enter as \emph{constraints} to a problem with unknown velocities in the boundary condition \eqref{rigidbc}.

We solve both the \emph{resistance problem}
and its inverse, the \emph{mobility problem}
(see Remark~\ref{r:uinfty} below for certain details).
In the resistance problem, the translational and angular velocities \( (\vec v^{(i)}, \omega^{(i)}) \) of each rigid particle are prescribed, which fixes
the boundary velocity Dirichlet data \eqref{rigidbc}.
Once the exterior Dirichlet BVP \eqref{stokeseq} is solved,
the desired hydrodynamic forces and torques \( (\vec f^{(i)}, t^{(i)}) \) are extracted
via \eqref{eq:forces} in a post-processing step.
Applications of the resistance problem include porous media and microfluidic devices.
In the mobility case the roles are reversed: the net forces and torques are specified, and one seeks the resulting rigid particle velocities and angular velocities.
Thus the given forces and torques enter as \emph{constraints} that must be enforced alongside \eqref{stokeseq} and \eqref{rigidbc} (see Sec.~\ref{sec:mob} for the numerical approach).
Applications of the mobility problem include sedimentation, rheology, motile swimmers and active fluids.

\vspace{2ex}
\textbf{Computational challenges.}  
When rigid particles undergo relative motion at close separations, their hydrodynamic interactions become especially challenging to determine, both in the resistance and mobility settings \cite[p.~175]{Lefebvre-Lepot2015,Kim1991}. Two distinct difficulties appear:
\begin{itemize}
\item \textbf{Lubrication-driven fine scales.}  
In narrow inter-particle gaps, the fluid velocity develops steep gradients. Accurate capture of these gradients requires very high spatial resolution:  
volume-based methods, such as finite elements, must deploy extremely fine meshes~\cite{Lefebvre-Lepot2015}.  
Potential based methods, such as  boundary integral equations (BIEs) or the \emph{method of fundamental solutions} (MFS), avoid the need for volume discretization, but not the demand for local refinement. Relative motion drives sharp peaks in the surface force density.  
For two disks, this peak scales like $\mathcal{O}(\sqrt{\delta})$ with gap width $\delta$~\cite{Barnett2015, Sangani1994}, reflecting the rapid amplification of lubrication forces between the particles.
If the density peaks are under-resolved,  the error in computed hydrodynamic quantities---forces and torques in a resistance problem, or velocities in a mobility problem---can be substantial. %, since these outputs are commonly obtained as integrals over the layer density. 
 % Local refinement in close-to-touching areas, however, increases the computational cost significantly, due to the long-range nature of the elliptic set of PDEs.
Moreover, refinement in close-to-touching regions significantly increases computational cost because elliptic kernels are globally coupled: locally introduced fine-scale degrees of freedom must interact with the entire suspension.
\item \textbf{Ill-conditioning.}  
Regardless of discretization, the resulting linear system inherits the singular nature of the physics.  
As $\delta\to0$, or as the number of particles grows, 
% *** Anna: this is not necessarily true for mobility. AHB: should we change?
the conditioning deteriorates
\cite{Barnett2015, Quaife2018}.
When using an iterative solver such as GMRES, the number of iterations needed to reach fixed accuracy rises without limit as $\delta\to 0$ \cite{malhotra2023}; for instance
in the 3D resistance setting \cite{Broms2025}
this was empirically found to be $\mathcal{O}(\delta^{-1/2})$.
The problem is therefore not merely one of accuracy, but also of efficiency.
% *** AHB I added some; this is a but vague - be more specific? Anna: Other people's work?
\end{itemize}
\textbf{Related work.} 
% Despite recent progress in developing accurate quadrature schemes for boundary integral methods \cite{Helsing2008,Bremer2010}, resolving the layer density, behaving much like the surface force density, remains difficult. Adaptive refinement, such as dyadic clustering of quadrature nodes near the points of closest approach between particles, is needed. 
% An approach to circumvent the resulting large number of degrees of freedom  is recursively compressed inverse preconditioning (RCIP), developed by Helsing \cite{HelsingRCIP} and extended to Stokes flow in \cite{Bystricky2019,chunkie}. The technique stabilizes the discretization of boundary integral equations near geometric singularities or near-singularities, such as those arising at corners, or in our context for close-to-touching bodies. In a preconditioning step, recursive refinement is done in a region of the boundary centered around the exact location of the singularity or near-singularity, followed by compression. In a dynamic setting, however, the mobility problem has to be solved repeatedly, at each time-step with an updated geometry. For general such time-dependent settings, the RCIP  technique is challenging and costly to apply, as also the singularity locations are updated with time. Another caveat is that RCIP is currently available only in 2D and not readily extendable to 3D. The compression shares similarities with other schemes to handle singularities arising at corners (developed for the Laplace equation), where an efficient basis reduces the number of unknowns \cite{Bremer2012b,Hoskins2019,Hoskins2020}.
%
We briefly overview numerical methods for Stokes flows with near-contacts,
focusing on PDE-based approaches. % (see also \cite{Pozrikidis1992,Yan2020,Broms2025}).
We first note the long tradition of using far-field approximations of the
hydrodynamic interactions between bodies (usually spheres in 3D),
such as in Stokesian dynamics \cite{Brady1988},
with lubrication corrections added in a pair-wise manner \cite{Sangani1994,Lefebvre-Lepot2015}.
% *** explain why pair is not enough?
Regularized effective particle interactions are also popular, as in the rigid multiblob \cite{Balboa2017,Broms2022}, force coupling \cite{Su2024} or regularized Stokeslets \cite{cortez2001} methods.
While useful, neither method is {\em convergent}, in the sense that the
error in solving the BVP \eqref{stokeseq} may be reduced as 
close to zero as desired by adjusting numerical parameters.

PDE-based approaches---which in contrast do converge to the true Stokes solution---either
discretize the fluid volume (e.g.~finite element methods with a conforming mesh),
exploit potential theory and linearity to discretize only the boundaries
(as in BIE and MFS), or combine both (as in immersed boundary methods and cut finite element methods).
Despite recent advances in accurate quadrature for BIEs \cite{Helsing2008,Bremer2010,Barnett2015,AfKlinteberg2016},
resolving the layer density---which closely mirrors the physical force density---remains challenging.
In near-contact situations, the density develops sharp peaks that demand refinement, for example via dyadic (exponential) clustering of quadrature nodes near points of closest approach.
One way to curb the resulting growth in degrees of freedom is recursively compressed inverse preconditioning (RCIP), introduced for corners 
% *** Anna: Helsing's tutorial also includes close to touching so I cite both Helsing-Ojala and the tutorial here
by Helsing \& Ojala \cite{HelsingOjala}, demonstrated for close-to-touching bodies in \cite{HelsingRCIP} and extended to Stokes flow in \cite{Bystricky2019,chunkie}. 
Starting from a dyadic local refinement in a neighborhood of the singularity, RCIP
uses a sequence of small dense direct linear solves to compress the refined system to a much smaller effective one.
For static geometries, this is highly effective, but in dynamic problems
where the geometry changes at every time step this can be expensive.
%and with it the location of the near-singularities; recomputing the compression at each step can become prohibitively expensive.
A further major limitation is that RCIP has not been generalized to 3D.

The MFS (also known as the method of auxiliary sources or the charge simulation method~\cite{Alves2004,Fairweather2005,Barnett2007,Alves2009,Liu2016,Karageorghis2019,Antunes2022}) moves the
potential sources off the boundary,
removing the need for singular quadratures that complicates BIE.
%and into the nonphysical domain.
An exterior solution is represented as a linear combination of Stokeslets placed inside each particle. For instance, for the unit disk, a standard choice is to place ``proxy'' sources uniformly on an interior curve of radius $R_p<1$, enforcing boundary conditions in a least-squares sense at collocation points. The resulting rectangular matrix
%target–from–source matrix
becomes exponentially ill-conditioned in the high-accuracy regime,
so that a backward-stable solve---typically via %a truncated
dense singular value decomposition (SVD)---is essential.

%%sssssssssssssss
% The method of fundamental solutions (MFS), also known as the method of auxiliary sources or the charge simulation method \cite{Alves2004,Fairweather2005,Barnett2007,Alves2009,Liu2016,Karageorghis2019,Antunes2022}, shares much of the structure of boundary integral methods. The singular integrals associated with the latter are however entirely avoided by representing the solution to an elliptic PDE as a linear combination of fundamental solutions centered at source points strictly inside each particle.  
% For unit-radius circles, so-called \emph{proxy} sources are typically placed uniformly on an interior curve of radius $R_p < 1$. Boundary conditions are enforced in the least-squares sense at a slightly larger set of target collocation points on the particle boundary. The resulting target-from-source matrix becomes exponentially ill-conditioned as $R_p \to 0$. To control this, the least-squares problem must be solved in a backward-stable fashion, typically via a truncated singular value decomposition (SVD).  

% our prior work:
In near-contact situations, simply increasing the number of proxy sources to resolve lubrication forces becomes prohibitive. For Stokes spheres in 3D \cite{Broms2025},
%\cite{Broms2024c} <- this paper seems not relevant to your point. I edited down too. Anna: yes, was wrong ref.
% AHB: don't need this here: inspired by Cheng \& Greengard’s image representation for close conducting discs \cite{Cheng1998}
we enhanced the basic MFS setup for every pair of particles by adding various types of sources clustered toward \emph{image accumulation points}.
%whose locations depend only on the particle-particle separation.
We also stabilized the multi-particle solve via right-preconditioning using the particle self-interaction matrix blocks (``one-body'' or rectangular block-Jacobi preconditioning) \cite{Liu2016,Stein2022}.
This enabled accurate resolution of lubrication effects with a modest number of degrees of freedom for all gaps larger than $10^{-3}R$, covering what we argued to be all physically relevant separations for Stokes flow. For the mobility problem we introduced a ``recompleted'' formulation in \cite{Broms2024c} in which force/torque constraints are automatically satisfied, allowing an unconstrained least-squares solve. Combined with one-body preconditioning and fast multipole acceleration, this yielded a robust linear-complexity scheme for general smooth shapes, demonstrated on clusters of $10^4$ ellipsoids.

There are other PDE-based approaches that share with MFS the idea of
least-squares collocation. Examples include Crowdy et al.'s conformal mapping method for Laplace problems exterior to close cylinders~\cite{Crowdy2016} and the Stokes lightning method of Trefethen and collaborators~\cite{Brubeck2022,Xue2024}. These exploit complex analysis to express the solution in terms of analytic Goursat functions (for Stokes) or Laurent series (for Laplace). The expansions in \cite{Crowdy2016} even include terms centered at image accumulation points, similar to our MFS enhancement. While accurate, they are difficult to accelerate with fast summation techniques and, to our knowledge, have not been applied to close-to-touching rigid bodies in Stokes flow.
Such complex analytic methods are also intrinsically tied to 2D.

We know of no existing convergent method that combines the flexibility of the MFS with a preconditioning strategy able to resolve lubrication effects at physically relevant small gaps, while avoiding a global fine discretization---crucial for simulations with very large numbers of particles. Here, we present one.

% FFFFFFFFFFFFFFFFFFFFFFFFFFFFFFFFFFFFFFFFFFFFFFFFFFFFFFFFFFFFFFFFFFFFFFFFFFFFFFFf
\begin{figure}[htpb]
    \centering
    \hspace*{-4ex}
    \begin{subfigure}[t!]{0.34\textwidth}
        \includegraphics[trim={2cm 11.5cm 4.5cm 6cm},clip,width = 1.2\textwidth]{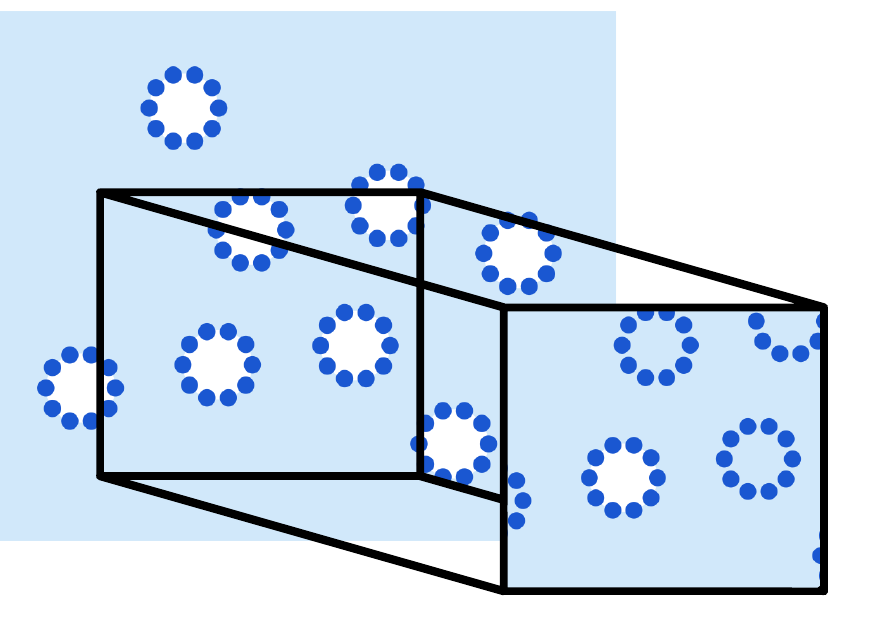}
        \caption{Dilute case}
      \label{fig:one-body}   \end{subfigure}~~~~~~~
    \begin{subfigure}[t!]{0.31\textwidth}
    \vspace*{5ex}
    \reflectbox{
        \includegraphics[trim={0cm 0cm 0cm 0cm},clip,angle=0,width = 0.95\textwidth]{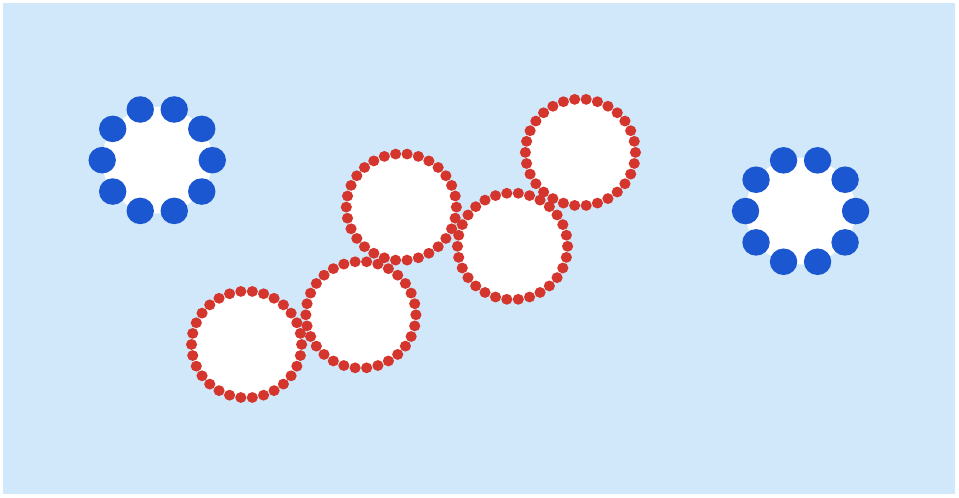}}
        \vspace*{8ex}
        \caption{Global fine grid}
    \end{subfigure}\hspace*{-5ex}
       \begin{subfigure}[t!]{0.32\textwidth}
       \hspace*{2ex}
        \includegraphics[trim={0cm 11.5cm 2cm 2cm},clip,angle=180,width = 1.05\textwidth]{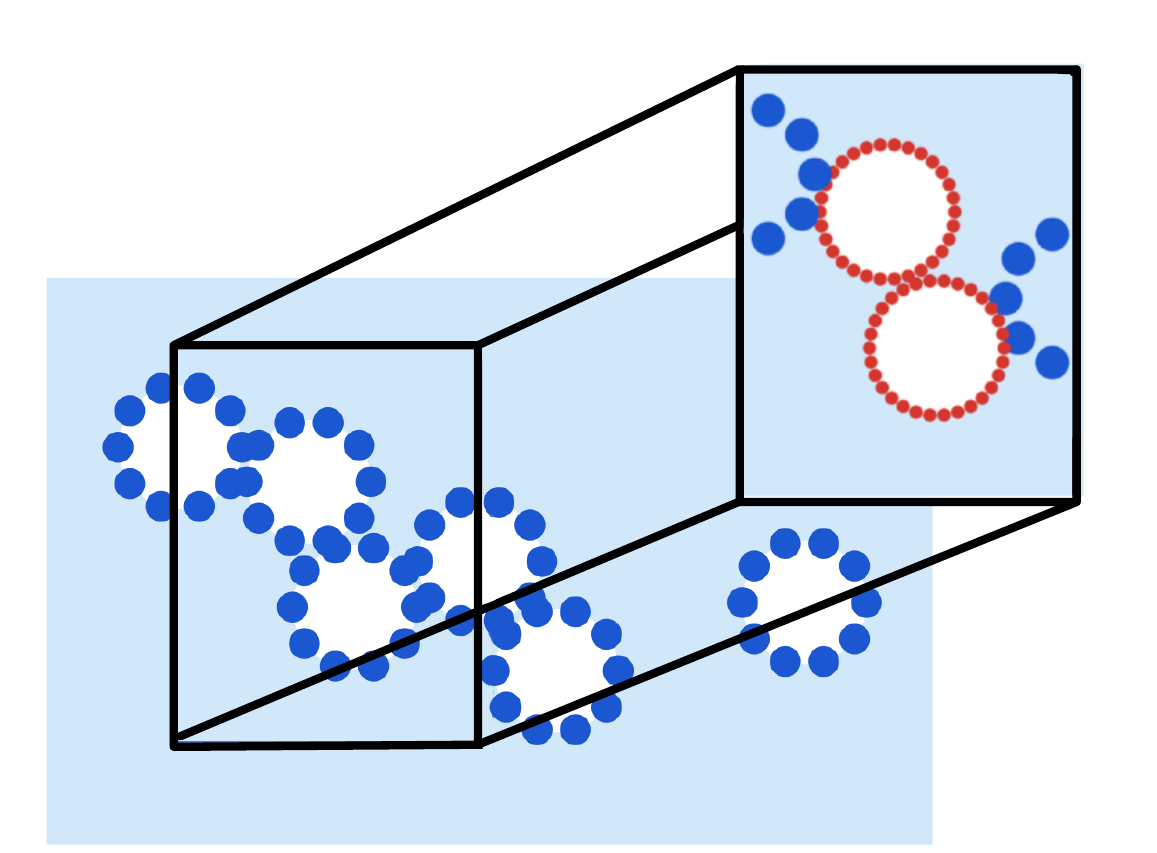}  \hspace*{4ex}
         \par\vspace*{2.2ex}
        \makebox[\textwidth][c]{%
    \hspace*{5ex}%
    \parbox{0.9\textwidth}{%
      \caption{Proposed two-body basis idea}
      \label{fig:two-body}
    }%
  }
        % \hspace*{4ex}\caption{Proposed two-body basis idea}
        % \label{fig:two-body}
    \end{subfigure}
       \caption{Sketch of the two-body preconditioner idea.
         Large blue dots indicate coarse surface discretization
         (MFS collocation nodes; source points are not shown).
         Small red dots indicate a more expensive fine discretization.
         (a) In dilute suspensions, a coarse discretization of each particle
         suffices to resolve all interactions, and a (one-body) basis for the flow field is obtained by using MFS to solve each body in isolation (inset).
         (b) In denser systems, any particle undergoing a near contact
         demands a fine discretization to capture local hydrodynamics.
         (c) The proposed scheme constructs a basis in which such fine discretizations are used only at the level of isolated pairs (inset), with each such interaction compressed to an equivalent coarse representation in the global solve.
         Fine grids are retained only locally for resolving near-field interactions and post-processing.
         In each panel and each inset, light blue indicates the fluid domain for the BVP solved.
       }
   \label{fig:idea}
\end{figure}

%%%%%%%%%%%%%%%%%%%%%%%%%%%%%%%%%%%%%%%%%%
\subsection{Summary of contributions}

The paper presents a {\em two-body preconditioner} that mitigates both the severe ill-conditioning,
and growth in the number of unknowns, that usually occur when seeking accurate Stokes solutions with close rigid particles in relative motion.
It is a hybrid of direct and iterative methods:
%In a similar flavor to RCIP for integral equations,
it directly solves a set of local finely-discretized BVPs---each involving
only two nearly-touching particles---in order to construct a (preconditioned) coarse discretized global system involving far fewer unknowns.
The global system is then solved iteratively using fast multipole (FMM) acceleration.
%The resulting hybrid scheme combines direct compression at the finely-resolved scales
%with fast accelerated iterative solution for the large scales.
%The spirit is similar to one-level multigrid schemes, and recent one-level hybrid fast direct solvers such as that of Lorca et al.~\cite{Gillman24}.
%It stabilizes the global solve while retaining a coarse global discretation, in spirit akin to a one-level multigrid scheme, or a one-level hybrid fast direct solver
The idea generalizes the
pairwise Laplace image-sum basis functions of Cheng \& Greengard \cite{Cheng1998, Cheng2000}
to ``two-body basis functions'' which may be precomputed with {\em any} convenient
%numerical
BVP solver.
We use MFS for this task in the present work, due to its excellent performance
for large scale Stokes flows \cite{Broms2025,Broms2024c}.
An overview of the scheme is given in Figure \ref{fig:idea}.

Our contribution has two main elements, each targeting one of the above-mentioned
%computational
challenges:

\begin{enumerate}
\item \emph{Stabilizing the ill-conditioned global system.} The preconditioner is built by solving directly for pairwise corrections to the one-body basis for each particle in the suspension as necessary, with each correction capturing the
  near-field interaction with a close neighbor.
  The result is a two-body basis representation that can approximate all possible
  near-contact interactions.
  The unknown coefficients of these basis functions are then solved for
  iteratively, using an FMM-accelerated global matrix-vector product.
  %we show (using another level of proxy point compression on ``peanut'' shaped
  The use of two-body bases greatly reduces GMRES iteration counts,
  even for gaps as small as $10^{-3}R$, in both resistance and mobility problems. 
\item \emph{Capturing lubrication-driven fine scales locally.} Each pairwise correction is obtained by solving a small but {\em high-resolution} BVP with only two particles (Figure~\ref{fig:two-body}, inset), fully resolving the steep surface-force peaks induced by lubrication. This is achieved using an image-enhanced MFS---a 2D version of that of \cite{Broms2025}---which also uses a hairpin curve to give high accuracy with only a single source type (Stokeslets).
  %adapted here to 2D and based solely on Stokeslet sources rather than a mixture of source types;
  We apply near-contact MFS image enhancement to the mobility problem for the first time. A subsequent %``peanut'' collocation nodes
  compression step replaces the fine pair representation with an equivalent coarse set of sources (Figure~\ref{fig:two-body}, background), so that the global discretization remains coarse.
  This addresses the local resolution requirements without increasing the global degrees of freedom (passed to the FMM) in the iterative solve.
\end{enumerate}

\begin{remark}[Connection to hybrid PDE solvers]\label{r:hybrid} % rrrrrrrrrrrrrrrrrrrrrrrrr
  % ***AHB add preconditioned?
  This ``hybrid'' of
  {\em direct} solution operators for local problems, compressing them to reduce
their effective number of unknowns, and then solving the remaining global problem
{\em iteratively} (block-diagonal preconditioning) recurs in diverse settings. These include: % in the numerical PDE literature:
\begin{itemize}
\item
  RCIP \cite{HelsingOjala} and other BIE corner-compressions schemes for elliptic PDE
  \cite{Bremer2012b,Hoskins2019,Hoskins2020} use essentially this idea, where the
  refined panel geometry at a single corner determines which fine unknowns are compressed.
  Other local preconditioners can be advantageous with BIE \cite{quaifepreco}. 
\item
  The two-body basis functions of Cheng--Greengard \cite{Cheng1998,Cheng2000} (mentioned above)
  fit into this framework, where the fine-scale ``solve'' is instead a direct image sum,
  and no acceleration of the global matrix-vector multiply was used.
\item
  The {\em scattering matrix} of a single body is a certain set of one-body basis functions
  (in our terminology). After such a direct one-body solution is built,
  it is commonly compressed using proxy and/or collocation points, and again
  used in a global
  iterative solve; this is sometimes called {\em fast multi-particle scattering}
  \cite{FMPS}. See \cite{Liu2016,Martinsson2026} for MFS versions for well-separated bodies.
\item
  {\em Fast direct solvers} (FDS) \cite{Martinsson2019} use this idea hierarchically,
  where the grouping of fine unknowns is often done via
  a quad- or oct-tree \cite{hackbusch,Martinsson2005}.
  Sometimes the coarse nodes are a subset of the fine; this is known as
  {\em skeletonization}. %and common for BIE but not appropriate for MFS - we don't use it
  A low-accuracy or incomplete FDS may also be used as a successful preconditioner,
  resulting in a hybrid scheme. 
  For porosity problems without relative body motion, Quaife et al.~\cite{Quaife2018} presented a hybrid solver of this type. A notable recent volumetric hybrid solver is due to Lorca et al.~\cite{Gillman2024}. 
  Our proposal is a pairwise extension of a hybrid 1-level FDS
  (iterative solution of \cite[Eq.~(13.14)]{Martinsson2019}).
\end{itemize}
    %or for treating corner singularities in Laplace problems \cite{Bremer2012b,Hoskins2019,Hoskins2020}, in which an adapted specialized basis captures the singular behavior and thereby reduces the number of unknowns and iteration count.  % <- Nice, I like this.
  In each of these settings, both the number of unknowns and the iteration count is reduced
  compared to a global solution using the fine-scale discretization.
  In contrast to almost all of the above,
  our proposal requires {\em rectangular least-squares} rather than square solves because of the nature of the MFS; this complicates the linear algebra but removes the need for a quadrature scheme for singular integrals. 
  %simplifies the quadrature schemes. % (with an eye towards 3D).
\end{remark}

\begin{figure}[h!]
  \centering
  \includegraphics[trim = {0cm 0cm 0cm 0cm},clip,width = \textwidth]{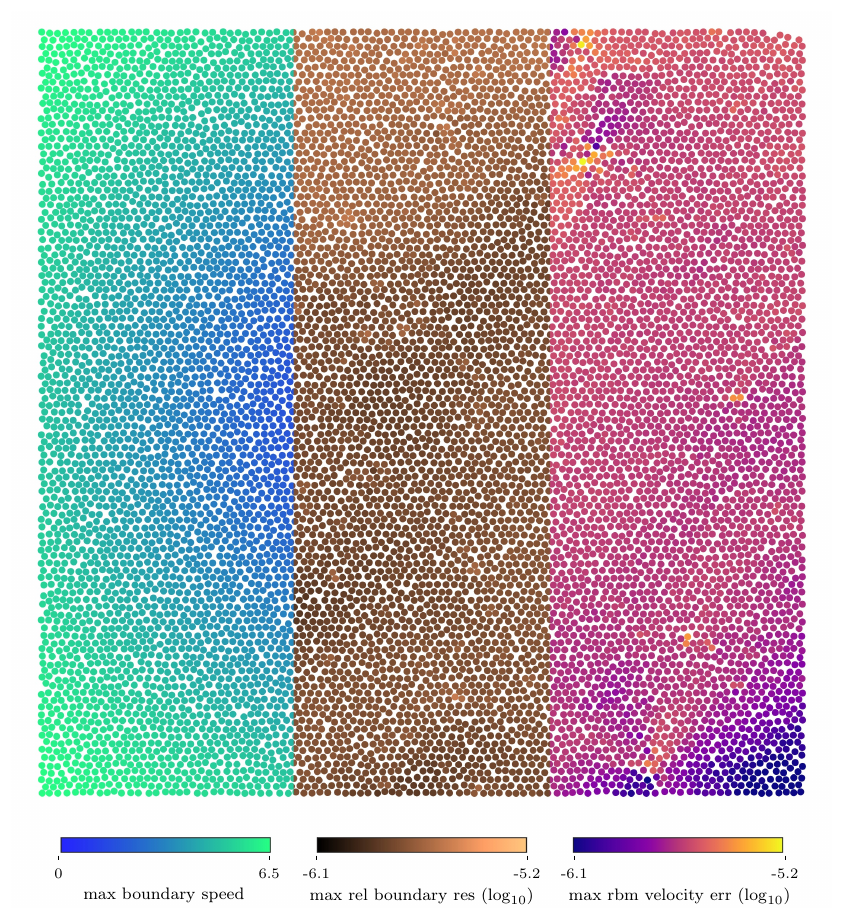}
%   \begin{subfigure}[t]{\textwidth}
%     \centering
%     \includegraphics[width = \textwidth]{figures/err_random_mc_P10000_phi065.pdf}
%     \caption{Max relative boundary residual}
%   \end{subfigure}
%   \hfill
%   \begin{subfigure}[t]{0.2\textwidth}
%     \centering
% \includegraphics[width = \textwidth]{figures/vel_random_mc_P10000_phi065.pdf}
%     \caption{Max boundary velocity}
%   \end{subfigure}
\caption{
Mobility solution for a random close packing of $10{,}000$ unit circles (packing density $\varphi = 0.65$, minimum separation $10^{-3}$), with each particle driven by a unit torque. The global iterative solve uses 72 vector unknowns per body and converges in 47 GMRES iterations and 36\,s on a single compute node. From left to right, a third of the particles is colored by velocity magnitude, followed by maximum relative boundary residual, and relative rigid body motion error, as defined in \eqref{speed},\eqref{residual} and \eqref{rel_vel_err}. The velocity field shows slower motion in the interior of the cluster. The boundary residual is measured with respect to the boundary data, and the velocity error is computed relative to a finer reference solution. }
   \label{large_ex}
\end{figure}

Together, these components yield a fast, accurate, and robust MFS framework for dense suspensions. Numerical experiments confirm stable performance in challenging multi-particle configurations.  As a motivational example, the mobility problem is solved for a system of 10,000 unit circles (packing fraction $\varphi = 0.65$) in Figure \ref{large_ex}. Despite minimum particle separations as small as $10^{-3}$, with $9920$ close pairs requiring local resolution, the interaction is resolved to a relative surface residual uniformly below $10^{-5}$. Convergence is achieved in 47 GMRES iterations, with a parallel solve time of 36\,s on a single compute node
(see Remark~\ref{r:cpu}).
Although demonstrated in 2D, the ideas are readily applicable in 3D.
% AHB add parallel timing? Anna: Done. Also reported in caption and returned to in the results section. Too much repetition?

\paragraph{Paper overview.}
Section~\ref{sec:general} sets the stage by introducing the two-body preconditioning idea in a general BVP framework. We then specialize to the resistance problem for geometries of circular bodies using the MFS in Section~\ref{sec:res}, where the construction proceeds in three steps: (1) identify one-body basis functions using the preconditioned MFS of \cite{Liu2016,Broms2024c,Broms2025}; (2) introduce two-body corrections to the one-body basis using fine sources; and (3) apply ``peanut compression'' by matching, on a peanut-shaped separation surface, the flow field from the fine sources of a particle pair to that of a set of coarse proxy sources. The resulting linear system involves only coarse unknowns. The same construction is then applied to the mobility problem in Section~\ref{sec:mob}. The choice of source points for the fine pairwise representation is discussed in Section~\ref{sec:fine}, while numerical experiments in Section~\ref{sec:num} validate the accuracy and efficiency of the method for both resistance and mobility problems. Finally, Section~\ref{sec:conc} concludes and outlines directions for future work.

% \alert{PLEASE IGNORE! THIS IS FROM OLD 2D TEXT IN THESIS}
% 	We present a technique where a low number of source points (and unknown MFS coefficients) per object in the domain are enough to solve the BVPs in \eqref{laplaceeq} and \eqref{stokeseq} to controllable accuracy. For closely interacting objects, such results are obtained by placing additional source points at approximate image locations. Image points and circle inversion will be the topic of Section \ref{inversion}. Section \ref{2D} describes the choice of source, collocation and approximate image points, and types of singularities needed at these locations to reach a target accuracy of $10^{-6}$, for circular particles and the Stokes BVP.  Techniques for solving the ill-conditioned and overdetermined linear system resulting from sampling the source and collocation points are presented in Section \ref{solving}. Then, numerical examples for Laplace problems are presented in Section \ref{Laplace}, with image points used also for non-circular particles, even if developed as a tool for circles. Example Stokes problems are solved in Section \ref{StokesNum}, where the parameter selection from Section \ref{2D} is tested for circular particles. We also demonstrate that dilute suspensions of non-circular particles can be accurately simulated. Numerical results are compared to a reference solution computed with a boundary integral method equipped with special quadrature in the exterior domain, or to the known value at the boundary. Our findings are summarised in Section \ref{conclusions}. 

\begin{remark}[Details on uniqueness and translational reference frame] % rrrrrrrrrrrrrrrrrrrrrrrrrrrrrrrrrrrrrrrrrrrrrrrrr
  \label{r:uinfty}
  Unlike in 3D where a constant flow is always a valid limit as $r\to\infty$,
  in 2D an exterior flow $\vec u$ as in \eqref{stokeseq}
  with nonzero total force
   $\vec \Sigma$ grows without bound, leading to the so-called {\em Stokes paradox}
  \cite{heywood74}
  \cite[\S6.14]{Pozrikidis2011} \cite[\S5.1]{Graham}.
  It is easy (e.g.\ by expanding the exterior Green's representation
  \cite[(2.3.20)]{Hsiao2008}) to show that the last line of \eqref{stokeseq}
is equivalent to
  \be
  \vec u(\vec x) = \frac{1}{4\pi \mu}
  \biggl( \log \frac{1}{r} + \frac{\vec x \vec x^T}{r^2} \biggr) \vec\Sigma
  + \vec w + o(1), \qquad r := \|\vec x\| \to \infty.
  \label{winfty}
  \ee
%  where $\vec w$ may 
  In the mobility problem $\vec\Sigma=\sum_i \vec f^{(i)}$ is specified,
  but also for uniqueness the constant $\vec w$ must be specified \cite[\S3]{Rachh2016}
  (physically the latter is needed to select a specific Galilean reference frame).
  In our mobility tests we always set $\vec\Sigma=\mbf{0}$ and $\vec w =\mbf{0}$, the
  latter being enforced numerically by our pure-Stokeslet MFS representation.
  Turning to resistance, 
  its standard BVP with Dirichlet data and $\vec\Sigma$ specified
  always has a unique solution \cite[\S2.3]{Hsiao2008}, from which $\vec w$ could be
  extracted via \eqref{winfty}.
  In this work our representation
  in fact solves a nonstandard
  resistance problem where Dirichlet data and $\vec w=\mbf{0}$ are specified
  (this avoids the complication of an additional constant
  in the flow representation \cite{Alves2004}).
  This nonstandard BVP is uniquely solvable except for a set of zero measure geometries
  (e.g., one disk of radius $\sqrt{e}$), but these can easily be bypassed by rescaling space.
  Moreover the standard resistance solution could be recovered by solving cases with
  $\vec u|_\pO \equiv (1,0)$ and $(1,0)$ then inverting the resulting $2\times 2$ linear map between
  $\vec w$ and $\vec \Sigma$. For simplicity---and since
  most applications involve mobility and/or 3D---we do not dwell on this further.
\end{remark} % rrrrrrrrrrrrrrrrrrrrrrrrrrrrrrrrrrrrrrrrrrrrrrrrrrrrrrrrrrrrrrrrrrrrrrrrrr

%%ssssssssssssssssssssssssssssssssssssssssssssssssssssssssssssssssssssssssssssssssss
\section{General framework for two-body basis construction}\label{sec:general}

We construct an efficient basis for the flow field due to $P$ particles in three steps:
\begin{enumerate}
\item Solve one-body BVPs for particle $i$ in isolation (see inset in Figure \ref{fig:one-body}) to obtain the one-body basis functions associated with particle $i$.
  % ***AHB: note the "basis" is really plural, even for 1 body (it has 2M functions). Anna: Later (e.g. caption Figs 3 and 5), we write that "each correction is computed by solving a BVP with..." Some inconsistency between "BVPs" and "a BVP".
\item
  For each close neighbor of particle $i$, add a correction to this basis involving
  BVPs for the pair comprising $i$ and the neighbor (see inset in Figure \ref{fig:two-body}).
\item Represent the flow due to all $P$ particles as a superposition of such pair-corrected basis functions.
\end{enumerate}
 
In this section we outline these steps separately for the resistance and mobility problems, in a deliberately high-level manner.
We assume only the existence of a local BVP solver (for one or two particles)
whose input is discretized surface data.
The reader is reassured that the abstraction will be made concrete
in Section~\ref{sec:res} and the sequel, where we specialize to the MFS as the local BVP solver
and to disks as the particles.
% Anna: ok, I like this!
%%%%%%
\subsection{The resistance (Dirichlet) problem}\label{sec:general_res}
%*** Anna: I also like the new structure here! 
%We approximate interpolation from its vector of samples
Consider the $i$th body, and 
let $\{\x^{(i)}_j\}_{j=1}^M$ be its set of \emph{coarse} boundary nodes, % on $\pO^{(i)}$,
sufficient for accurate discretization of $\pO^{(i)}$ in isolation.
We discretize boundary velocity data on these nodes to give the vector
$\mm^{(i)} = \{\vec \mu_1^{(i)},\dots,\vec \mu_M^{(i)}\}$.
Smooth boundary functions may then be accurately interpolated from these nodes.
Let $\vec\phi^{(i)}$ denote the {\em velocity solution operator} for the $i$th body in isolation,
meaning that it maps $\mm^{(i)}$ to the resulting flow field
$\uu(\x) = \vec\phi^{(i)}[\vec \mu^{(i)}](\x)$
which solves
the Stokes BVP \eqref{stokeseq} in the exterior of $\Omega^{(i)}$ alone, with Dirichlet
data $\mm^{(i)}$.
(Recall Remark~\ref{r:uinfty} for the condition at infinity.)
%$\vec\phi^{(i)}$ may be applied by any accurate PDE solver of choice, and of course has an analytic solution for the disk.
This means that for any vector $\mm^{(i)}$ the reproducing property holds:
$\vec\phi^{(i)}[\mm^{(i)}](\x^{(i)}_j) = \mm^{(i)}_j$, $j=1,\dots,M$.
Note that here and beyond we use ``representation'' notation
$\vec\phi^{(i)}[\vec \mu^{(i)}]$, similar to that used in potential theory.

By linearity, one may view $\vec\phi^{(i)}$ as a set of \emph{one-body basis functions}
$\{\vec\phi^{(i)}_j\}_{j=1}^{2M}$, where
$\vec\phi^{(i)}_j$ is the flow solution $\uu$ when the data vector $\mm^{(i)}$
is set to the $j$th unit vector in $\R^{2M}$.
For the full system of $P$ particles, the flow field can be expressed as a superposition of fields generated by each body,
\begin{equation}\label{evaluate}
  \uu(\x) = \sum_{i=1}^P \vec\phi^{(i)}[\vec \mu^{(i)}](\x),
  \qquad \x \in \R^2 \setminus \overline{\Omega}
\qquad \mbox{ (one-body representation).}
\end{equation}
The full resistance problem \eqref{stokeseq} may now
be solved by treating \eqref{evaluate} as an ansatz with $\{\vec \mu^{(i)}\}_{i=1}^P$
as \emph{unknowns}, and using collocation to impose that $\uu$ match
the given velocity data $\{\vec g^{(i)}\}_{i=1}^P$ at all coarse nodes.
The resulting linear system has a $2PM\times 2PM$
system matrix with identity blocks along the diagonal, because of the one-body solution property.
In practice, one never forms the matrix, but applies it to vectors using a fast algorithm
and block-diagonal corrections \cite{Liu2016,Stein2022,Broms2025}.

\begin{remark}%[One-body preconditioning]
In the context of boundary-based PDE solvers (BIE, MFS, etc) the above method
is known as one-body (block diagonal) preconditioning, and transforms the unknowns from densities (or source strengths) to surface data. In the MFS case we recap this
in Section \ref{1B_precond}. %Such a preconditioner yields a well-conditioned system only \emph{in the absence} of strong lubrication effects \cite{Broms2025}.
%, but GMRES iteration counts grow rapidly as particle gaps become small 
\end{remark}

However, when two particles become close (forming a \emph{near contact}), two problems arise:
i) the boundary functions become
nonsmooth, thus cannot be accurately discretized or interpolated with the coarse nodes, and
ii) the conditioning of the above linear system deteriorates.
One solution to problem i) is simply to use a finer set of one-body boundary nodes (larger $M$);
% near-contact adapted?
however, this would slow down each matrix-vector apply, while doing
nothing to address ii).
This motivates a two-body preconditioning method. This uses pairwise BVPs
discretized on fine nodes to modify the above coarse one-body
bases $\vec\phi^{(i)}$ to two-body bases denoted by $\vec\psi^{(i)}$, giving the
global representation for the flow field
\begin{equation}\label{rep2}
  \uu(\x) = \sum_{i=1}^P\vec \psi^{(i)}[\mm^{(i)}](\x),
  \qquad \x \in \R^2 \setminus \overline{\Omega}
  \qquad \mbox{ (two-body representation).}
\end{equation}
The coefficients \(\vec\mu^{(i)}\) are then still determined by collocation at the coarse boundary nodes, so that the linear system remains of size $2PM\times 2PM$.

For each body (without loss of generality we describe this for body $i=1$),
its set of $2M$ two-body basis functions are constructed as follows.
If body 1 has no near contacts, we simply set \(\vec\psi^{(1)} = \vec\phi^{(1)}\).
If it has a single near contact, with body 2 (say), we define
$\vec \psi^{(1)} =\vec \phi^{(1)} + \vec \eta^{(1,2)}$. Here
\(\vec\eta^{(1,2)}\) denotes a \emph{correction basis} that solves the exterior Stokes BVP
for the pair $\pO^{(1)} \cup \pO^{(2)}$ with Dirichlet velocity data
\begin{equation}\label{local_BVP}
\vec\eta^{(1,2)}[\vec\mu^{(1)}](\vec x) =
\begin{cases}
    \vec 0, & \vec x\in\partial\Omega^{(1)}\\
    -\vec\phi^{(1)}[\vec\mu^{(1)}](\vec x), & \vec x\in\partial\Omega^{(2)}
\end{cases}
\qquad\mbox{ (correction data),}
\end{equation}
which is solved numerically using a fine discretization of both particle boundaries.
In particular, the $j$th correction function, $\vec\eta^{(1,2)}_j$, is found
by solving the pair BVP using boundary data given by evaluating the $j$th one-body basis $\vec\phi^{(1)}_j$ with a negative sign on body 2,
%(this requires accurate evaluation close to the source curve)
and zero data on body 1.
The correction thus cancels the one-body Dirichlet data on body 2.
Thus if bodies 1 and 2 have no other close neighbors, performing the analogous
correction for body 2 completes a pair of two-body bases $\vec\psi^{(1)}, \vec\psi^{(2)}$ 
that \emph{completely solves the pairwise Dirichlet BVP}, obeying the reproducing property
\begin{equation}
    \vec\psi^{(1)}[\vec\mu^{(1)}](\vec x_j^{(k)}) = \begin{cases}
        \vec\mu^{(1)}_j, & k = 1,\\
        \vec 0, & k = 2,\\
    \end{cases}\qquad
    \vec\psi^{(2)}[\vec\mu^{(2)}](\vec x^{(k)}_j) = \begin{cases}
        \vec 0,& k = 1,\\
        \vec \mu^{(2)}_j, & k = 2,
    \end{cases}
    \qquad j = 1,\dots, M.
\end{equation}
In this case the system matrix would now have a $4M \times 4M$ identity block on the
diagonal (for the unknowns of both bodies 1 and 2).
This bypasses any accuracy loss and ill-conditioning associated with the
one-body representation of the near contact,
at the cost of solving $4M$ finely-discretized local BVPs.

%If the pair is far away from all other bodies, matching to the boundary data yields $\vec \mu^{(1)} = \vec g^{(1)}$ and $\vec \mu^{(2)} = \vec g^{(2)}$. In contrast to the one-body representation, this is well-conditioned for pairs of particles close to each other.  

\begin{figure}[h!]
  \centering

  % Left column: (a) above (b)
  \begin{minipage}{0.19\textwidth}
    \centering

    \begin{subfigure}{\linewidth}
      \centering
      \includegraphics[width=1.1\linewidth]{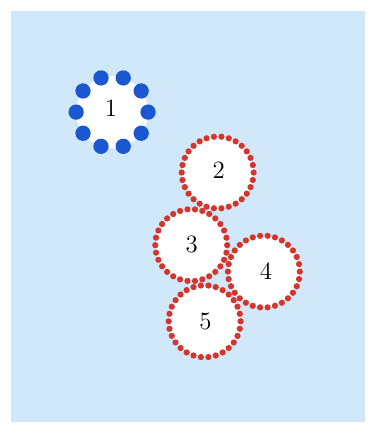}
      \caption{}
      \label{fig:a}
    \end{subfigure}

    \vspace*{5ex}

    \begin{subfigure}{\linewidth}
      \centering
      \includegraphics[width=\linewidth]{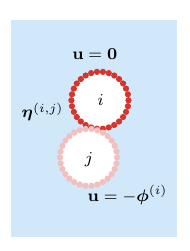}
      \caption{}
      \label{fig:b}
    \end{subfigure}
  \end{minipage}
  \hspace*{1ex}
  % Right column: (c)
  \begin{minipage}{0.78\textwidth}
    \centering

    \begin{subfigure}{\linewidth}
      \centering
      \includegraphics[width=\linewidth]{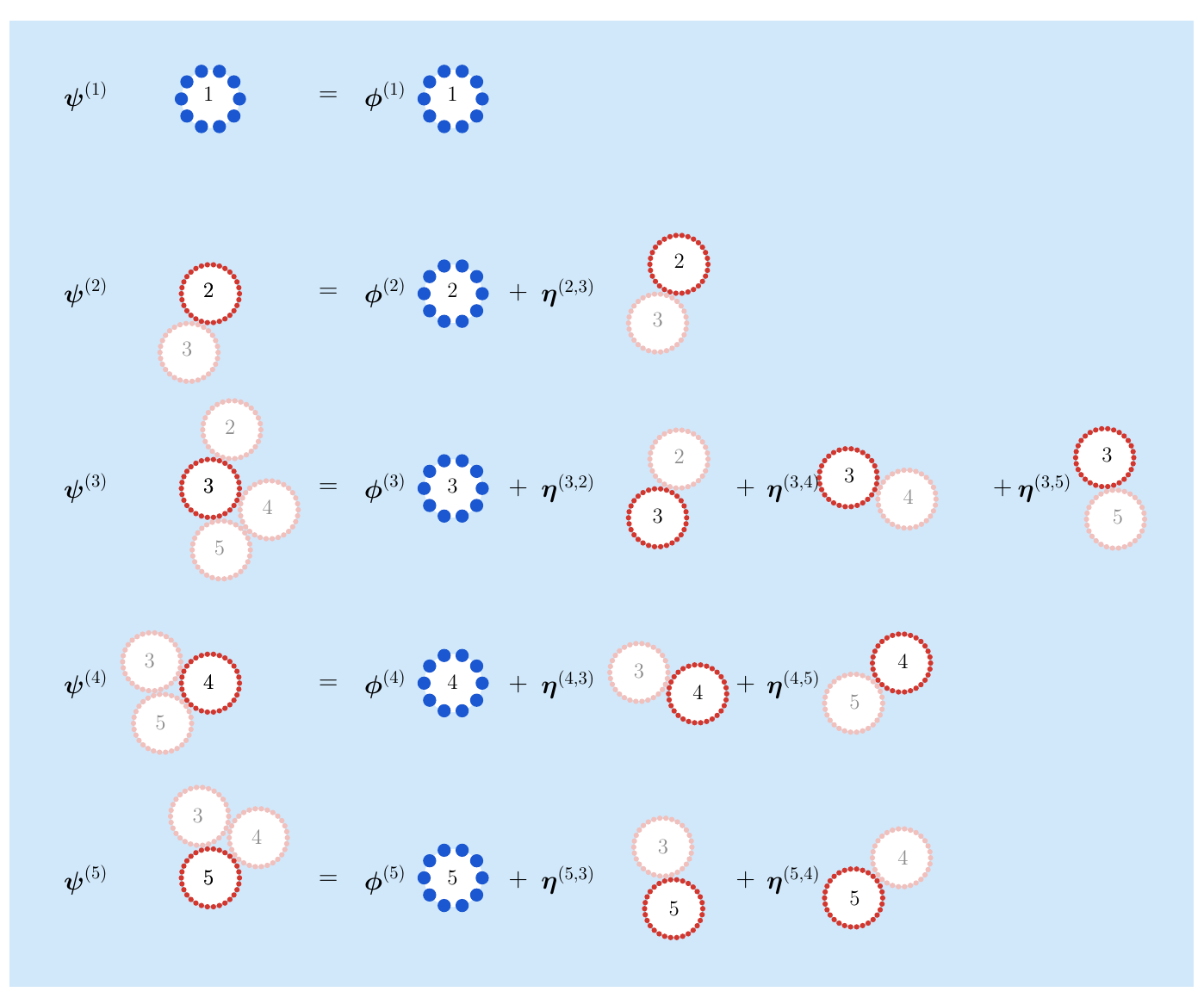}
      \caption{}
      \label{fig:c}
    \end{subfigure}
  \end{minipage}

  \caption{For the five particles shown in panel (a), the
    construction of their two-body bases
    $\vec\psi^{(1)}, \dots,\vec\psi^{(5)}$ is diagrammed in panel (c). For bodies in near contact, the one-body basis is augmented with pairwise corrections $\vec\eta^{(i,j)}$,
    each computed by solving a BVP with a single neighbor as in panel (b), using a
    fine discretization (red). For well-separated particles, no correction is needed: the one-body basis alone suffices, and the coarse discretization (blue) resolves all interactions.}
  % ***ahb something about compression here maybe?
  \label{fig:multicontacts}
\end{figure}

If body 1 has multiple near contacts, their %corresponding
pairwise corrections are simply summed: % to get the new basis for body 1:
\begin{equation}\label{eta_sum}
  \vec \psi^{(1)} = \vec \phi^{(1)} + \sum_{k\in \mathcal C^{(1)}}\vec \eta^{(1,k)}
%    \vec \psi^{(1)}[\mm^{(1)}](\x) = \vec \phi^{(1)}[\mm^{(1)}](\x) + \sum_{k\in \mathcal C^{(1)}}\vec \eta^{(1,k)} [\mm^{(1)}](\x)
%  \qquad \x \in \R^2 \setminus \overline{\Omega}
  \qquad \mbox{ (two-body basis for body 1), }
\end{equation}
with $\mathcal C^{(1)}$ the set of neighbor indices of body $1$.
This is repeated for the rest of the bodies $i=2,\dots,P$.
The resulting multi-contact construction of two-body bases for $P=5$ particles is illustrated in Figure~\ref{fig:multicontacts}.
This idea avoids the large system sizes associated with fully coupled multi-body treatments for near contacts,
while maintaining low GMRES iteration counts, which we demonstrate in Section~\ref{sec:num}.

%%%%%%%%%%%%%%

% \begin{figure}[htb]
%   \centering
%   \subcaptionbox{}{ %
%     \vspace*{30ex}
%     \includegraphics[scale=0.6\paircorrScale]{figures/fine_group.pdf}%
%   }
%  % \hspace{0.01em}
%   \subcaptionbox{}{%
%     \includegraphics[scale=0.55\paircorrScale]{figures/fine_group_basis_rot.pdf}%
%   }
%     \hspace{0.06em}
%   \subcaptionbox{}{%
%       \vspace*{30ex}  
%      \includegraphics[scale=0.9\paircorrScale]{figures/pair_problem.pdf}%
%   }
%   \caption{Panel (b) constructs the two-body basis $\{\vec\psi^{(i)}\}_{i=1}^P$ for the $P=7$ particles in panel (a). For bodies in near contact, the one-body basis is augmented with pairwise corrections, each computed by isolating the interaction with a single neighbor. The fine discretization (red) is required only at the level of individual pairs. For well-separated particles, no correction is needed: the one-body basis alone suffices, and the coarse discretization (blue) resolves all interactions. Panel (c) illustrates the pair BVP for body $i$ in near contact with body $j$.}
%   \label{fig:multicontacts}
% \end{figure}
%%%%%%%%%%%%5fffffffffffffffffffffffffffffffffff

 \begin{remark}
   The above generalizes a method of Cheng \& Greengard for the iterative solution of Laplace Dirichlet BVPs with multiple disks \cite{Cheng1998} or spheres \cite{Cheng2000}; they used an analytic multipole image sum for the pairwise basis, giving a slow matrix-vector apply.
   Our generalization is to arbitrary boundary-based linear PDE solvers in multi-body geometries,
   given only a black-box pairwise solver, and allows for accelerated matrix-vector applies.
   %
%   Interpreting image sums as local BVP solutions allows for accelerated matrix-vector applies.
%   
   We refer to \cite{Cheng1998} for a discussion of why pairwise bases are
   sufficient even in the case of multiple clusters, such as a
   triangle of three near-touching disks.
   % AHB: I think this has already been said in Rmk 1:
%In this broader setting, the two-body basis functions of \cite{Cheng1998} are reinterpreted as local BVP solutions, enabling acceleration.  
%Our scheme constitutes a partial factorization, or hybrid method and may be viewed as a close-touching analog of the corner-adapted bases in \cite{Bremer2012b,Hoskins2019,Hoskins2020} or 1-level RCIP \cite{HelsingRCIP}.
\end{remark}
%%%%%%%%%%%%%%%%%%%%%%%%%%%

 \subsection{The mobility problem}
 \label{s:mob_gen}
% [Anna: modified after disc with AK]
% AHB: I have trouble with this section - the main issue is that the RBM needs to be described as some linear functional of \mu. It may be that there is *no* "local mobility BVP" interpretation of the S(I-L)(S_L)^+ rep that our one-body mobility precond uses!
% ***AHB (6/6/26): conclusion now is: there *is* a local mobility solve formulation of the 1-body basis (in an MFS-indep way), but it is a little odd - I have explained it. The difficulty is that it's not possible to explain the PM*PM lin sys having Id diag blocks in the mobility case, because of needing the ansatz for the RBM = B K^T.mu itself :(
% The 2-body basis I am not confident about and we can touch up at the end...

For mobility, the one-body flow field representation is
\begin{equation}\label{mob_compflow}
\uu(\x)
=
\sum_{i=1}^P
\vec \phi^{(i)}[\mm^{(i)}](\x)
+
\sum_{i=1}^P
\vec C^{(i)}\left[\vec f^{(i)},t^{(i)}\right](\vec x),
% AHB: note I made C completions body-specific not global. It is more flexible.
\quad \x \in \R^2 \setminus \overline{\Omega}
\quad \mbox{(mobility one-body representation).}
\end{equation}
The first term is a sum of mobility-specific one-body basis functions, each producing zero net force and torque, while $\vec C^{(i)}$ is a known \emph{completion flow}
\cite{PowerMiranda,Pozrikidis1992,Bagge2021,malhotra2023,Broms2024c}, that is,
an explicit exterior Stokes solution
carrying (in the sense of \eqref{eq:forces}) the prescribed force and torque for body $i$, but in general not a rigid body motion on any body.
The construction of each above one-body function differs from the resistance setting.
%here we use a choice generalizing a mobility method due to the authors \cite{Broms2024c}.
% Anna: this is the general description. Do we really want to say it's MFS solver specific? AHB: actually I think it's pretty specific, but edited it.
Rather than the plain reproducing property, we have a rank-3 perturbation of it:
for any surface velocity vector $\mm^{(i)}$ on body $i$,
\begin{equation}\label{mobility_bc}
  % AHB: note I added missing mu... Anna: I don't understand. Do you mean \vec v^{(i)}[\vec\mu^{(i)}] and \omega^{(i)}[\vec\mu^{(i)}]? I'm a bit confused now. Don't we want to say that \vec v and \omega are unknown? The ansatz for \vec v and \omega is MFS specific.
  % AHB: actually the ansatz for v and omega *has* to be part of the 1-body basis defn. Without trying this out with BIE, it's hard to say it's MFS-specific or not. But it does match our MFS-mob choice.
\vec \phi^{(i)}[\mm^{(i)}](\vec x_j^{(i)})
=
\mm^{(i)}_j
+ ({\cal R}^{(i)}\mm^{(i)})_j,      % ***AHB: yes, it's +
%\vec v^{(i)} + \omega^{(i)}(\vec x_j^{(i)}-\vec c^{(i)})^{\perp},
\qquad j=1,\dots,M,
\end{equation}
where ${\cal R}^{(i)}$ is a
%method-specific
$2M\times 2M$ matrix of rank 3 whose range is the subspace of rigid body velocity data, namely
$\{\vec v + \omega(\vec x_j^{(i)}-\vec c^{(i)})^{\perp}\}_{j=1}^M$
for all $\vec v \in \R^{2}, \omega \in \R$.
In practice, ${\cal R}^{(i)}$ is fixed by an ansatz
made for the specific local BVP solver used, that couples rigid body motions to unknowns
(e.g.\ for BIE see \cite[Eq.~(38)]{malhotra2023}, while for MFS see \cite[\S3.2]{Broms2024c}).
%where $\vec v^{(i)}$ and $\omega^{(i)}$ are certain linear functionals of $\mm^{(i)}$.
To solve the full mobility problem, one treats the one-body coefficients $\{\mm^{(i)}\}_{i=1}^P$
as unknown, and uses collocation to enforce that the representation \eqref{mob_compflow} equals
the ansatz ${\cal R}^{(k)}\mm^{(k)}$ at the coarse nodes of the $k$th body, for all $k=1,\dots,P$.
As with the resistance problem,
the resulting system matrix has size $2PM\times 2PM$, with identity blocks along the diagonal,
and enables an accelerated apply.
However, now the right-hand side becomes the surface data of the completion flow
from \eqref{mob_compflow} with a negative sign.
Full details are given in Section~\ref{sec:mob}.

Applying pair-corrections to the mobility one-body basis functions we similarly get
\begin{equation}\label{rep2mob}
  \uu(\x) = \sum_{i=1}^P\vec \psi^{(i)}[\mm^{(i)}](\x) +
  \sum_{i=1}^P
  \vec C^{(i)}\left[\vec f^{(i)},t^{(i)}\right](\vec x),
  \quad \mbox{(mobility two-body representation).}
\end{equation}
with $\vec \psi^{(i)}$ again a two-body basis constructed with additive corrections to the one-body basis $\vec \phi^{(i)}$ as in \eqref{eta_sum}, but now with the additional constraints on the corrections $\vec\eta^{(i,k)}$ that they produce no net force or torque on the bodies within the pair.
With body 1 and body 2 sufficiently far away from any other neighbors, so that no other corrections are needed, the pair corrected basis functions for the two bodies satisfy
\begin{equation}\label{corr_syst_mob}
\begin{aligned}
  \vec\psi^{(1)}[\vec\mu^{(1)}](\vec x_j^{(k)}) = & \begin{cases}
    \vec\mu^{(1)}_j + ({\cal R}^{(1)}\mm^{(1)})_j,
    %\vec v^{(1)}-\omega^{(1)}(\vec x_j^{(1)}-\vec c^{(1)})^{\perp}
      & k = 1,\\
        \vec 0, & k = 2,\\
    \end{cases}\\
    \vec\psi^{(2)}[\vec\mu^{(2)}](\vec x_j^{(k)}) = & \begin{cases}
        \vec 0, & k = 1,\\
        \vec\mu^{(2)}_j + ({\cal R}^{(2)}\mm^{(2)})_j,
        %- \vec v^{(2)}-\omega^{(2)}(\vec x_j^{(2)}-\vec c^{(2)})^{\perp},
        & k = 2,
    \end{cases}
    \end{aligned}
    \qquad j = 1,\dots,M.
\end{equation}
% ***AHB: Still I'm not sure I believe this 2-body case (are the \vec 0 cases not also RBMs?
Details on how these local mobility problems are solved using MFS are given in Section \ref{sec:mob}.

A summary of the functions used to represent the flow field in both the resistance and mobility settings is given in Table \ref{tab:functions}.

\begin{table}[h]
\centering
\begin{tabular}{ l l }
\textit{Function} & \textit{Description} \\
\hline
\\[-2mm]
$\displaystyle \vec\phi^{(\partone)}(\vec x)$ & One-body velocity solution operator for particle $\partone$. \\
$\displaystyle \vec \phi^{(\partone)} [\mm^{(\partone)}](\x)$ & One-body basis representation for body $i$ with coefficient vector $\vec\mu^{(i)}$. \\ &
For MFS-specific expressions, see: \eqref{b1rep} for resistance and \eqref{1body_mob} for mobility. \\
$\displaystyle \vec \eta^{(\partone,\parttwo)}[\mm^{(\partone)}](\x)$
& Pair-correction to one-body basis. For MFS-specific expressions, see: \eqref{chi} for \\ & resistance  and \eqref{chi_mob} for mobility, where the  function is evaluated using fine sources \\ & $\displaystyle \vec\beta^{(\partone,\parttwo)}$, that in turn depend on $\mm^{(\partone)}$. 
For efficiency, $\vec\eta^{(i,k)}$ is ``peanut compressed'' \\ & together with $\vec\eta^{(k,i)}$ and evaluated using the coarse correction vector $\tilde{\vec{\lambda}}^{(\partone\text{-}\parttwo)}$.  \\ & For expressions, see \eqref{peanut} for resistance and \eqref{peanut_mob} for mobility. \\
$\vec \psi^{(\partone)}[\mm^{(\partone)}](\x)$ & Two-body basis for particle $\partone$: 
%$\displaystyle
%\vec \psi^{(\partone)}[\mm^{(\partone)}](\x) = \vec \phi^{(\partone)}[\mm^{(\partone)}](\x) + \sum_{k\in \mathcal C^{(\partone)}}\vec \eta^{(\partone,k)} [\mm^{(\partone)}](\x)$.
$\displaystyle \vec \psi^{(\partone)} = \vec \phi^{(\partone)} + \sum_{k\in \mathcal C^{(\partone)}}\vec \eta^{(\partone,k)}$.
\\
$\vec C^{(i)}\left[\vec f^{(i)},t^{(i)}\right](\vec x)$ & Completion flow for the mobility problem, supplying the $i$th body force and torque $\vec f^{(i)},t^{(i)}$.\\
\\[-2mm]
\hline
\end{tabular}
\caption{Functions appearing in the flow field representations. The same notation is used to describe functions in both the resistance and mobility problems, although their construction is based on different boundary value problems. Here, $\mathcal C^{(i)}$ denotes the set of indices of bodies in near contact with body~$i$.
}
\label{tab:functions}
\end{table}

 %sssssssssssssssssssssssssssssssssssssssss
\section{Solving the resistance problem via MFS}\label{sec:res}

We now apply the MFS to the approach outlined in Section~\ref{sec:general_res}, specializing
to the case where $\Omega^{(i)}$ are unit-radius (monodisperse) disks. We start with the MFS itself, followed by one-body preconditioning (recapping prior work). Recasting the latter in the new basis-function framework then provides a natural route to the new two-body algorithm.
% We start with the MFS itself, then one-body preconditioning (recapping prior work),
% which we then present in the new basis interpretation, in order then
% to give a complete presentation of the new two-body algorithm.

% Old notation comment from Alex:
% $2N_c$ is the number of coarse degrees of freedom on a body ($N_c$ coarse sources),
% $M_c$ the number of coarse collocation nodes ($2M_c$ conditions).
% We ignore getting such factors of 2 correct for now since
% it complicates bookkeeping.

The fundamental solution to the Stokes equations---the 2D \emph{Stokeslet}---describes the velocity field induced by a point force. It is given by the $2 \times 2$ matrix-valued kernel
\begin{equation}\label{stokeslet2d}
\mathbb S(\vec x, \vec y) = \frac{1}{4\pi\mu} \left[ -\log \|\vec x - \vec y\|\vec I_2 + \frac{(\vec x - \vec y)(\vec x - \vec y)^T}{\|\vec x - \vec y\|^2} \right],
\end{equation}
with $\vec I_2$ the $2\times 2$ identity matrix.
 %This kernel satisfies the Stokes equations in the fluid domain $\mathbb{R}^2 \setminus \{\vec y\}$.
A general velocity field $\vec u$ solving the Stokes equations in the multiply-connected exterior can be approximated as a linear superposition of Stokeslets:
\begin{equation}\label{stokes_sum}
\vec u(\vec x) = \sum_{i=1}^P \sum_{j=1}^N \mathbb S(\vec x, \vec y^{(i)}_j) \vec \lambda^{(i)}_j, \qquad \vec x \in \mathbb R^2\setminus\overline{\Omega},
\end{equation}
where $\vec y^{(i)}_j$, $j=1,\dots,N$, are source points
inside disk $i$, and the associated vector source strengths $\vec \lambda^{(i)}_j \in \mathbb{R}^2$ are to be determined.
For well-separated disks, source points may be chosen equispaced on the concentric circle of radius $R_c<1$, with typically $0.6 \le R_c \le 0.8$.
% ***AHB: felt like the reader shouldn't wait until p. 26 for this.
However, for lubrication driven problems with near-touching disks, more elaborate choices
are needed (as in Section~\ref{sec:fine}). % with refinement in close-to-touching regions. Such source refinement also goes with refinement of the boundary collocation nodes.  % Anna: great! But also the collocation points are refined in this setting and no longer correspond to the coarse colloc points... Tried to add but notation gets a little messy. $\vec Y^{(i)}$ are later said to be coarse and the same for \vec\lambda.. .

For ease of notation, let $\vec Y^{(i)}  = \{ \vec y_j^{(i)} \}_{j=1}^N$ and $\vec\lambda^{(i)} = \lbrace \vec\lambda_j^{(i)}\rbrace_{j=1}^N$ denote the stacked source points and strengths.  %(sometimes referred to as proxy points). 
We can then write the representation \eqref{stokes_sum} compactly as
\begin{equation}\label{stokes_operator}
\vec u(\vec x) = \sum_{i = 1}^P \mathbb{S}(\vec x,\vec Y^{(i)}) \vec \lambda^{(i)}.
\end{equation}
% where~$\mathcal{S}^{(i)}$ denotes the single-layer Stokes operator associated with sources on particle $i$.

The source strengths vectors
$\vec\lambda^{(i)}$, $i = 1,\dots,P$, are determined  by enforcing boundary conditions in the least-squares sense at all sets of target coarse collocation points $\vec X^{(i)} \coloneqq 
\{ \vec x^{(i)}_j \}_{j=1}^{M}$ 
on the $i$th particle boundary. 
It has been found that setting $M$ slightly larger than $N$ increases MFS accuracy \cite{Barnett2007}.
% We will later refer to this basic discretization as the \emph{coarse} grid.
% Anna: this gets a little messy... g is notation for bc at coarse nodes. Can't use same notation if we want to write this for the more general setting with adapted source and collocation sets.
Letting $\vec g$ again denote the stacked boundary data at all such coarse collocation nodes, and $\vec \lambda$ the stacked source coefficients, we obtain an
overdetermined least-squares linear system
\begin{equation}\label{stokes_system}
\vec S \vec \lambda = \vec g,
\end{equation}
where $\vec S$ is a dense $2MP \times 2NP$ matrix whose blocks are given by evaluating $\mathbb S(\vec x, \vec y)$ between all pairs of target and source points. As already alluded to, the matrix $\vec S$ is exponentially ill-conditioned. In addition, it is typically too large for a dense backward-stable solve
if $P\gg 1$. Hence, efficient preconditioning is needed to transform \eqref{stokes_system} into a well-conditioned square system amenable for an iterative solve.

%%%%%%%%%%%%%%%%%%%%%%%%%%%%%%
\subsection{One-body preconditioning}
\label{1B_precond}

We first review the one-body preconditioning approach of \cite{Liu2016, Broms2025, Broms2024c}, and then formulate it in the general framework of one-body basis functions of Section \ref{sec:general}. 
% The one-body framework also plays a crucial role in what follows: it serves as the foundation for the two-body preconditioner of the least-squares problem \eqref{stokes_system} developed in this paper, much more effective in resolving the strong hydrodynamic interactions that arise at small particle separations. 
We will sometimes refer to the point sets $\lbrace\vec X^{(i)}\rbrace_{i=1}^P$, $\lbrace\vec Y^{(i)}\rbrace_{i=1}^P$ respectively as the coarse collocation and source points.
%\begin{remark}[Source vector notation]
We will use the bar notation $\bar{\vec\lambda}$ to indicate the
vector $\vec\lambda$ obtained via one-body preconditioning alone.
% This equality will not hold with ``peanut compressed'' two-body preconditioning, where we will add a correction to $\bar{\vec\lambda}$ from each close pair to form the vector $\vec\lambda$ to be used in the representation \eqref{stokes_sum}; more on this in Section \ref{sec:peanut}. % Anna: vague! How to formulate?
%\end{remark}

\subsubsection{Standard formulation}\label{sec:standard}
 The global target-from-source MFS matrix $\vec S \in \mathbb{R}^{2MP \times 2NP}$ has the block structure,
\begin{equation}\label{Gmat}
	\vec S = 
	\begin{bmatrix} 
		\vec S^{(11)} & \vec S^{(12)} & \cdots & \vec S^{(1P)} \\
		\vec S^{(21)} & \vec S^{(22)} & \cdots & \cdots \\
		\vdots & \vdots & \ddots & \vdots \\
		\vec S^{(P1)} & \cdots & \cdots & \vec S^{(PP)} 
	\end{bmatrix},
\end{equation}
with the block 
\be
\vec S^{(ik)}_{mn} = \mathbb S(\x^{(i)}_m,\y^{(k)}_n), \quad m=1,\dots,M, \quad n=1,\dots,N
\ee
mapping source points on particle $k$ to target points on particle $i$.

Following \cite{Liu2016, Broms2025, Broms2024c}, we precondition the ill-conditioned overdetermined system $\vec S\bar{\vec\lambda} = \vec g$ from the right using a block-diagonal matrix with entries ${\vec S^{(ii)}}^{+}$, the pseudoinverses of the one-body self-interaction blocks.
Since all particles are identically shaped and sized, the self-interaction block ${\vec S^{(ii)}}$ is shared across bodies, so we
denote it simply by $\vec P$. For each particle, define the preconditioned unknown vector $\vec \mu^{(i)} = \vec P \bar{\vec \lambda}^{(i)} \in \mathbb{R}^{2M}$, so that $\bar{\vec \lambda}^{(i)} = \vec P^+ \vec \mu^{(i)}$. The pseudoinverse $\vec P^+$ is determined e.g.\ using the SVD, $\vec P=\vec U \vec \Sigma \vec V^{T}$, where $\vec\Sigma$ is the diagonal matrix containing the singular values $\sigma_1^{(i)}\ge \sigma_2^{(i)}\ge \dots \sigma_N^{(i)}$.
To ensure numerical stability, one has to apply the pseudoinverse in a backward-stable fashion \cite{TrefethenBau,Lai2015,Malhotra2015,Stein2022,Parolin2022}, so that 
	\begin{equation}\label{howG+}
		\bar{\vec \lambda}^{(i)} = \vec V \vec \Sigma^{+} \left(\vec U^T \vec \mu^{(i)}\right), \quad i = 1, \dots, P,
	\end{equation}
	where $\vec \Sigma^+$ denotes the truncated pseudoinverse of $\vec\Sigma$; its diagonal entries are set to $1/\sigma_j^{(i)}$ when $\sigma_j^{(i)} > \sigma_1^{(i)}\epsilon_{\text{trunc}}$, or zero otherwise. The truncation level $\epsilon_{\text{trunc}}$  is typically set smaller than the desired error, but somewhat larger than machine precision. 
    
    Since $N<M$, each diagonal block $\vec P\vec P^+$ of the resulting preconditioned target-from-source matrix has a nontrivial nullspace: it annihilates some $(2M - 2N)$-dimensional subspace of inputs. %has a subspace of dimension $2M-2N$ that maps $\vec\mu^{(i)}$ to vanishing velocity  values. %This motivates replacing 
%Since $N<M$, by a rank consideration,
%there must always be a subspace of dimension
%$M-N$ of vectors $\vec\mu^{(j)}$ that map to vanishing source strengths.
Thus, a well-conditioned matrix can only be
achieved by replacing the diagonal blocks by $\vec I$, which $\vec P\vec P^+$ approximates %whose action they approximate 
for vectors smooth on each body.
    The preconditioned system is now square and %therefore results in the following \emph{square} linear system, now 
    of size $2MP \times 2MP$:
\begin{equation}
\begin{bmatrix}
	\vec I & \vec S^{(12)} \vec P^+ & \cdots & \vec S^{(1P)} \vec P^+ \\
	\vec S^{(21)} \vec P^+ & \vec I & \cdots & \cdots \\
	\vdots & \vdots & \ddots & \vdots \\
	\vec S^{(P1)} \vec P^+ & \cdots & \cdots & \vec I
\end{bmatrix}
\begin{bmatrix}
	\vec \mu^{(1)} \\ \vec \mu^{(2)} \\ \vdots \\ \vec \mu^{(P)}
\end{bmatrix}
=
\begin{bmatrix}
	\vec g^{(1)} \\ \vec g^{(2)} \\ \vdots \\ \vec g^{(P)}
\end{bmatrix}.
\label{resistance}
\end{equation}
% after regularizing the diagonal matrix products ${\vec S^{(ii)}}{\vec P}^+$ by the $2M\times2M$ identity matrix. 
Once this has been solved iteratively for $\mm := \{\mm^{(i)}\}_{i=1}^P$,
the MFS strengths are recovered via \eqref{howG+}.

% \begin{remark}
% Note that Since $N<M$, by a rank consideration,
% there must always be a subspace of dimension
% $M-N$ of vectors $\vec\mu^{(j)}$ that map to vanishing source strengths.
% Thus, the well-conditioned nature of \eqref{resistance} is only
% achieved by replacing $\vec P\vec P^+$ by $\vec I_{M\times M}$ in all diagonal blocks, whose action they approximate for
% vectors smooth on each body The latter cannot be well-conditioned because $N<M$,
% thus in the actual system \eqref{resistance}
% they are replaced by $\vec I$ whose action they approximate for
% vectors smooth on each body.
% \end{remark}

\begin{remark}[New unknowns] One-body preconditioning reformulates the problem so that the unknowns, $\vec\mu$, are collocation velocity values at the boundaries rather than interior source strengths. Although this increases the number of unknowns slightly (typically, for the coarse grid, 
we choose $M=1.2N)$, the resulting system is square and much better conditioned, at least for sufficiently well-separated particles. For close-to-touching configurations, however, both $N$ and $M$ must be increased substantially, often with $M\gg N$, leading to much larger linear systems. This provides another indication that one-body preconditioning alone is insufficient for large-scale simulations of dense suspensions. By contrast, the two-body basis will resolve the near-contact interactions locally while retaining a coarse global discretization. The resulting reduction in the number of unknowns is quantified in Section~\ref{mob_res}. % Can we polish? Move around?
\end{remark}

\subsubsection{One-body basis reformulation and accelerated matrix-vector product}

We now recast \eqref{resistance} in terms of one-body basis functions.
It is easy to check that if we define
the MFS-solved one-body basis function for the $i$th body as
\be
   \vec \phi^{(i)} [\mm^{(i)}](\x) := \sum_{j=1}^M \vec\mu_j^{(i)} \vec\phi^{(i)}_j(\x)
   = \begin{cases} \vec\mu_m^{(i)},\quad &\vec x = \vec x_m^{(i)},\\
   \sum_{n=1}^N \mathbb S(\x,\y^{(i)}_n) (\vec P^+\mm^{(i)})_n,\quad&\text{otherwise,}
\end{cases}
\label{b1rep}
\ee
then the one-body basis collocation procedure described at the beginning of Section \ref{sec:general_res} gives precisely the preconditioned MFS linear system \eqref{resistance}.
Recall that for numerical stability, $\vec P^{+}$ must be applied as in \eqref{howG+}.
%The reproducing property gives the diagonal identity blocks.
Converting the resulting solution vector $\mm$ to $\vec\lambda$ using \eqref{howG+},
the MFS evaluation of the flow $\vec u$ at new targets using \eqref{stokes_sum} is
identical to the one-body basis sum \eqref{evaluate}.
In short, one-body (block-diagonal) preconditioning is equivalent to \emph{using
one-body basis functions that solve isolated-particle BVPs}.

%Given a surface coefficient vector
%$\vec\mu^{(i)}$ on the boundary of body $i$,
%evaluating its MFS representation
%$\vec u(\x) = \sum_{j=1}^N \mathbb S(\vec x, \vec y^{(i)}_j) \vec \lambda^{(i)}_j$

%This is a one-body flow field that solves (to the MFS solution accuracy)
%the BVP outlined in the beginning of Section \ref{sec:general_res}.
% keeping G for consistency.
%The global flow field representation is of the form \eqref{evaluate}.
% \begin{equation}\label{evaluate}
%   \uu(\x) = \sum_{i=1}^P[\vec \phi^{(i)}\mm^{(i)}](\x),
%   \hspace{2in}\mbox{(one-body representation)}
% \end{equation}
% which is how the flow field is reconstructed after solution for the 
% coefficients $\lbrace \vec\mu^{(k)}\rbrace_{k=1}^P$. 
% Note the rep *is* the evaluation.
%Evaluating \eqref{evaluate} at the coarse collocation nodes and enforcing the boundary data
%nearly 
%gives the linear system \eqref{resistance}. %; the difference is that
%it gives $\vec P \vec P^+$ as diagonal blocks.
% The latter cannot be well-conditioned because $N<M$,
% thus in the actual system \eqref{resistance}
% they are replaced by $\vec I$ whose action they approximate for
% vectors smooth on each body.

A practical aspect is to apply the huge matrix in \eqref{resistance} using a fast algorithm.
This matrix-vector apply follows the general algorithm outlined in Algorithm \ref{alg:general}. We will refer back to this basic algorithm later, when two-body corrections are added.
\begin{algorithm}
 \algrenewcommand\algorithmiccomment[2][\footnotesize]{{#1\hfill\(\triangleright\) #2}}  
 \algnewcommand{\LeftComment}[1]{\Statex \hspace*{1.6ex}\(\triangleright\) #1}
  \algnewcommand{\LeftLeftComment}[1]{\Statex \hspace*{3.6ex}\(\triangleright\) #1}
\caption{General matrix-vector apply for a Stokes BVP using a basis representation of the flow field}\label{alg:general}
\begin{algorithmic}
\Function{matvec}{$\vec \mu$}
    \State \textbf{Input:} Stacked velocity data vector $\vec \mu = \{ \vec\mu^{(i)}\}_{i=1}^P$ at all coarse collocation nodes $\lbrace \vec X^{(i)}\rbrace_{i=1}^P$ %Coarse collocation data $\vec \mu$ for all particles
    \State \textbf{Output:} Stacked surface velocity vector $\vec u = \{ \vec u^{(i)} \}_{i=1}^P$ at these same nodes $\lbrace \vec X^{(i)}\rbrace_{i=1}^P$ %Surface velocity values $\vec u$ at all coarse collocation points
    \LeftComment{Step 1: Map collocation data to source strengths internal to the basis}
    \LeftComment{Step 2: Evaluate velocity field at coarse collocation points using basis representation:}
    \LeftLeftComment{a. Compute global sum of all internal sources at all coarse collocation nodes via fast summation}
    \LeftLeftComment{b. Correct locally %within each body 
    to account for diagonal blocks}
      % ***AHB is that ok to add, or does it mess up something later? Anna: The corrections are not only per body, but also per pair (subtract off wrong velocity contribution, add back in what the fine grid should produce (rhs in lsq problem). 
    \State \Return $\vec u$
\EndFunction
\end{algorithmic}
\end{algorithm} 

Specifically, Step 1 and Step 2 in the linear-scaling matrix-vector apply in the one-body basis case is
\cite[Alg.~1]{Broms2024c}:
\begin{enumerate}
	\item[1.] Apply the self-interaction pseudoinverse backward-stably to each particle as in \eqref{howG+}: $\bar{\vec \lambda}^{(i)}  = \vec P^+\vec\mu^{(i)}$.
	\item[2a.] Evaluate the total velocity field from all sources  at  all collocation points (via, e.g., a Stokes FMM):
    \begin{equation}\label{getall}
   \vec u_{\text{All}} = \sum_{i = 1}^P \mathbb S\left(\vec X^{(i)},\vec Y^{(i)}\right) \bar{\vec \lambda}^{(i)}.
	\end{equation}
	\item[2b.] Subtract the self-interaction field and add the identity correction for each particle:
\begin{equation}\label{onesub}
\vec u^{(i)} = \vec u_{\text{All}}^{(i)} - \vec P \bar{\vec \lambda}^{(i)} + \vec \mu^{(i)}, \quad i = 1, \dots, P.
\end{equation}
This enforces the diagonal identity action in \eqref{resistance} (and the equivalent interpolation case in  \eqref{b1rep}).
\end{enumerate}

Once GMRES has converged, the final source coefficients $\bar{\vec \lambda}^{(i)}$ are recovered from $\vec \mu^{(i)}$, $i = 1,\dots,P$, using the pseudoinverse relation in \eqref{howG+}. The solution field can then be evaluated via the representation in \eqref{stokes_sum}, using $\vec\lambda = \bar{\vec\lambda}$. 
\subsection{Two-body preconditioning}\label{sec:2body}
% As rigid particles move relative to each other close distances apart, the number of GMRES iterations required for convergence with one-body preconditioning grows significantly, limiting the  efficiency. To address this, we now introduce two-body preconditioning, inspired by Cheng \& Greengard \cite{Cheng1998}.  The idea is to modify the basis in \eqref{basis} for close-to-touching particles. By a ``contact'', we hereafter refer to any pairwise interaction for which 

%*** Anna: Very nice! Easier to read than before :)

We now introduce pairwise corrections to the MFS one-body representation. Pair interactions are resolved locally per pair on a fine grid of $N_f$ source points and $M_f$ collocation points per particle (their precise location choices are deferred to Section~\ref{sec:fine}).
We use superscripts $(i,k)$ to denote quantities associated with $\vec\eta^{(i,k)}$, the basis correction to particle $i$ due to near contact with particle $k$, but $(i\text{-}k)$ to denote quantities needed \emph{jointly} for $\vec\eta^{(i,k)}$ and $\vec\eta^{(k,i)}$. The stacked vector of all $2N_f$ fine source points for the $(i\text{-}k)$ pair is denoted by $\vec{\mathcal{Y}}^{(i\text{-}k)}$, and the corresponding $2M_f$ collocation points by $\vec{\mathcal{X}}^{(i\text{-}k)}$. To distinguish fine and coarse discretizations, let $N_c$ and $M_c$ now denote the number of coarse source and collocation points per particle.
%recall for the $i$th particle their coordinates are stacked in $\vec Y^{(i)}$ and $\vec X^{(i)}$ respectively.
% Anna: Is the pair notation ok? AHB: yes!

% \begin{remark}[The fine grid]
% For simplicity, one may think of the fine grid as a refined distribution of Stokeslets placed in the interior of each of the two circles—finer than the coarse distribution used to construct the one-body basis—along with a denser sampling of collocation points. It is, however, often more efficient to include additional source types to better capture the near-singular behavior of lubrication forces, as described in the introduction. To accommodate this generality, we henceforth let $\vec G(\vec x,\vec y)$ denote a generic fundamental solution centered at $\vec y$ and evaluated at $\vec x$, not necessarily limited to the Stokeslet. While the precise composition of the fine grid is not essential at this stage, it will be fully described in Section~\ref{sec:fine}. 
% \end{remark}
% Unfortunately, we still expect iteration count to grow with the number of
%bodies; the only way around that would be fast direct solvers.

Our goal is to use the fine MFS to derive formulae for the pair correction basis $\vec \eta^{(1,2)}$ (for concreteness we pick $i=1$, $k=2$). This requires two ingredients: an algorithm to evaluate $\vec \eta^{(1,2)}[\mm^{(1)}](\x)$, and a procedure to compute the associated internal fine source strengths for the pair, denoted by $\vec\beta^{(1,2)} \coloneqq \{\vec\beta_q^{(1,2)}\}_{q=1}^{2N_f}$. %, with $\vec\beta_q^{(1,2)}\in\mathbb R^2$. 
The detailed steps outlined here are specific to the resistance problem, but the corresponding corrections to the one-body basis functions appearing in the mobility problem (see Section~\ref{sec:mob_two}) follow the same structure.
% AHB: good!

The first ingredient is simple:
once %a stacked vector
\(\vec\beta^{(1,2)}\) is known, we define \(\vec \eta^{(1,2)}\) via a
``fine'' sum of Stokeslets: %evaluation formula:
\begin{equation}
   \vec \eta^{(1,2)}[\mm^{(1)}](\x) = \sum_{q=1}^{2N_f} \mathbb S(\x,\vec{\mathcal Y}_q^{(1\text{-}2)}) \vec\beta_q^{(1,2)},
   \label{chi}
\end{equation}
with $\vec \beta^{(1,2)}$ some linear function of $\vec\mu^{(1)}$.
With \(\vec\beta^{(1,2)}\) known, the global representation \eqref{rep2} can be
evaluated with an FMM that includes
the fine source points \(\vec{\mathcal Y}^{(1\text{-}2)}\) and strengths \(\vec\beta^{(1,2)}\).

\begin{figure}[t]
    \centering
            	\includegraphics[trim={3cm 11cm 2cm 12.5cm},clip,width=\textwidth]{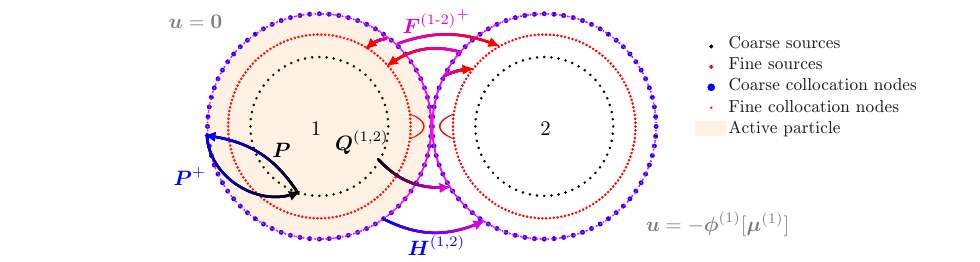}	
                \caption{Illustration of the linear maps used to construct the correction $\vec\eta^{(1,2)}$ to the basis of particle~1 (here called the active particle) due to near-contacting particle~2, as in Section~\ref{sec:2body}.}
                  %The fine source strengths $\vec\beta^{(1,2)}$ defining the correction basis are recovered from the coarse collocation data $\vec\mu^{(1)}$ on particle~1; see \eqref{sys2b}. This involves the maps ${\vec F^{(1\text{-}2)}}^+$ and $-\vec H^{(1,2)}$, with $\vec H^{(1,2)}=\vec Q^{(1,2)}\vec P^+$. The correction basis $\vec\eta^{(1,2)}$ cancels the velocity field generated by the one-body basis of body~1 on body~2 and vanishes on body~1 itself (shaded boundary data, enforced at the fine collocation nodes). }
                % ***AHB: as you see, I think this caption should be cut to the first sentence. The rest repeats the text. Anna: Ok! A-K had some comments about the figure being hard to follow. Perhaps we modified the text enough by now for the Figure to be clear enough?
   \label{fine_map}
\end{figure}

It remains to define the second ingredient: the linear map from $\mm^{(1)}$ to $\bb^{(1,2)}$.
For this we use the fine MFS with strength vector $\vec\beta^{(1,2)}$ to solve
the near-contact pairwise BVP as in \eqref{local_BVP}.
Enforcing its boundary data at the fine collocation nodes $\vec{\mathcal X}^{(1,2)}$
gives the overdetermined least-squares
$4M_f \times 4N_f$ system (again recalling the Stokes vector character where
each $\vec \beta_q^{(1,2)}$ is a 2-vector),
\[
\sum_{q=1}^{2N_f} \mathbb S(\vec{\mathcal X}^{(1\text{-}2)}_\ell,\vec{\mathcal Y}^{(1\text{-}2)}_q) \vec\beta_q^{(1,2)}
\;=\;
\left\{\begin{array}{ll}
\vec 0, & \ell = 1,\dots,M_f \qquad\qquad \mbox{(on body 1)}\\
-\vec \phi^{(1)}[\mm^{(1)}](\vec{\mathcal X}^{(1\text{-}2)}_\ell), & \ell = M_f+1,\dots,2M_f  \qquad \mbox{(on body 2)},\\
\end{array}\right.
\]
where we assume indexing of the fine collocation points for the pair,
$\vec{\mathcal X}^{(1\text{-}2)}_\ell$, such that the first $M_f$ lie on $\partial\Omega^{(1)}$ and the rest on $\partial\Omega^{(2)}$. 
Compactly, we write this as
\be
\vec F^{(1\text{-}2)} \bb^{(1,2)} = \begin{bmatrix} \vec 0 \\ -\vec H^{(1,2)}\end{bmatrix} \mm^{(1)},
\label{sys2b}
\ee
where $\vec H^{(1,2)}$ is a $2M_f\times 2M_c$ one-body basis evaluation matrix from the ``active'' body 1 to its neighbor.
% Its $2\times 2$ tensor valued 
% entries are $\vec H^{(1,2)}_{\ell j} = \vec\phi^{(1)}_j(\vec{\mathcal X}^{(1\text{-}2)}_{l+M_f})$, for $l=1,\dots,M_f$
% and $j=1,\dots,M_c$. 
Via the definition of the one-body basis in \eqref{b1rep}, $\vec H^{(1,2)}$ has its own factorization as $\vec H^{(1,2)} = \vec Q^{(1,2)} \vec P^+$, where $\vec Q_{\ell n}^{(1,2)} \coloneqq \mathbb S(\vec{\mathcal X}^{(1\text{-}2)}_{l+M_f},\y^{(1)}_n)$, $l = 1,\dots,M_f$, $n = 1,\dots,N_c$. The actions of the four matrices ${\vec F^{(1\text{-}2)}}^+$, $\vec H^{(1,2)}$, $\vec Q^{(1,2)}$ and $\vec P^+$ are illustrated in Figure \ref{fine_map}. % by applying the first factor of $\vec S^+$ leftward
%first.

In practice the stable solution of the ill-conditioned linear system \eqref{sys2b}
has two stages:
in a precomputation the matrices are filled, densely factorized via $\vec F^{(1\text{-}2)} = \vec U\vec\Sigma \vec V^T$,
and the factors $\vec V \vec\Sigma^+$
and $\vec U^T [\vec 0;-\vec H^{(1,2)}]$
stored, noting that
\be
\bb^{(1,2)} = {\vec F^{(1\text{-}2)}}^+ [\vec 0;-\vec H^{(1,2)}] \mm^{(1)} = \vec V\vec \Sigma^+ \left((\vec U^T [\vec 0;-\vec H^{(1,2)}]) \mm^{(1)}\right).
\label{bfm}
\ee
In subsequent ``apply'' stages (in each GMRES iteration), one uses the final above expression
as a pair of dense matrix-vector multiplies. This is simply a ``fine'' analog
of \eqref{howG+}.
This completes the recipe
to map $\mm^{(1)}$ to $\bb^{(1,2)}$. To evaluate the correction basis $\vec\eta^{(1,2)}$, the resulting $\vec\beta^{(1,2)}$ is inserted into \eqref{chi}.

Now equipped with the two-body basis, we evaluate the two-body flow field \eqref{rep2}
in the framework of Algorithm \ref{alg:general}:
Step 1 maps the coefficients $\mm$ to both the coarse source strengths $\bar{\vec\lambda}$ and the fine correction sources $\vec\beta$,
Step 2a inputs all fine and coarse source points to, e.g., a Stokes FMM.
The final step, 2b, proceeds exactly as in the one-body case (see \eqref{onesub}), since the two-body basis is constructed as an additive correction to the former.

% replace step 1 and 2 of the matrix-vector multiply outlined in Section \ref{sec:standard}(which uses the representation of the flow field in \eqref{evaluate} based on a one-body basis). Instead, the flow field at the coarse collocation nodes of all particles is evaluated via the representation in \eqref{rep2}: $\vec u_{\text{All}} = \sum_{i=1}^P[\vec \psi^{(i)}\mm^{(i)}]\left(\vec{\mathcal X}^{(j)}\right),\medspace j = 1,\dots,P.$
% As with one-body preconditioning, $\uu_\text{All}$ is not the matrix-vector multiply used and step 3 %Algorithm \ref{alg1} 
% in Section \ref{sec:standard} is left unchanged:
% As the two-body basis is based on corrections to the one-body basis, one must subtract the one-body diagonal blocks and add in the $\vec I$ contributions.

\begin{remark}[Symmetrizing pairs]
  It is faster to solve for $\vec \eta^{(i,k)}$ and $\vec \eta^{(k,i)}$ simultaneously, since both involve the same system matrix $\vec F^{(i\text{-}k)}$. For example, to compute $\vec \beta^{(2,1)}$ needed to evaluate $\vec \eta^{(2,1)}$, the right-hand side in \eqref{sys2b} becomes $[-\vec H^{(2,1)}; \vec 0] \mm^{(2)}$. The factorization of $\vec F^{(1\text{-}2)}$ from computing $\vec \beta^{(1,2)}$ can thus be reused, effectively halving the setup time. 
    % *** Anna: Do we want to have something like the box also here? E.g.
    One then computes the \emph{total} fine strength vector for the pair $\vec\beta^{(1\text{-}2)}\coloneqq\vec\beta^{(2,1)}+\vec\beta^{(1,2)}$ in the apply stage. A schematic of this apply stage is then, using arrows to denote linear maps:
\begin{center}
\setlength{\fboxsep}{8pt}
\framebox{%
\begin{minipage}{0.95\textwidth}
\centering
\makebox[1.3cm][c]{%
  \shortstack{%
    {\small Fine}\\
    {\small sources}\\
    $\vec{\beta}^{(1\text{-}2)}$%
  }%
}
$\xleftarrow[
  {\shortstack{\vspace{0.3ex}\\[-0.2ex]\scriptsize $\mathbb R^{4N_f\times 4M_f}$}}
]{
  \shortstack{\scriptsize ${\vec F^{(1\text{-}2)}}^+$\\[-0.2ex]\vspace{0.3ex}}
}$
\makebox[2.7cm][c]{%
  \shortstack{%
    \hspace*{-9ex}{\small Fine}\\
    \hspace*{-3ex}{\small collocation}\\%
    \hspace*{-4ex}{\small velocity}
  }%
}
$\xleftarrow[
  {\shortstack{\vspace{0.5ex}\\[-0.2ex]\scriptsize $\mathbb R^{4M_f\times 4N_c}$}}
]{%
  \makebox[1.0cm][c]{%
    \scriptsize
    $\mathord{-\left[\!\!
      \begin{array}{cc}
        \arraycolsep=1pt
        \vec 0 & \! \vec Q^{(2,1)} \\
        \vec Q^{(1,2)} & \! \vec 0
      \end{array}
    \!\!\right]}$%
  }%
}$
% \makebox[1.9cm][c]{\hspace*{2.5ex}\shortstack{\small{Coarse}\\\small{sources}}}
\makebox[4cm][c]{\hspace*{-5ex}\shortstack{\small{Coarse}\\ \small{sources}\\$\bar{\vec\lambda}^{(1)};\bar{\vec\lambda}^{(2)}$}} \hspace*{-10ex}
$\xleftarrow[
  {\shortstack{\vspace{0.5ex}\\[-0.2ex]\scriptsize $\mathbb R^{4N_c\times 4M_c}$}}
]{%
  \makebox[0.7cm][c]{%
    \scriptsize
    $\mathord{-\left[\!\!
      \begin{array}{cc}
        \arraycolsep=1pt
        \vec P^+ & \! \vec 0 \\
        \vec 0 & \! \vec P^+
      \end{array}
    \!\!\right]}$%
  }%
}$
\makebox[3.3cm][c]{\shortstack{\small{Coarse}\\\small{collocation data}\\$\vec\mu^{(1)};\vec\mu^{(2)}$}}
\end{minipage}
}% end framebox
\end{center}
\end{remark}
% *** AHB: I don't know if this box is needed. But it's nice that it sets up the reader for the peanut box later, so let's keep it.

%%%%5sssssssssssssssssssssssssssssssssssssssssssssssssssssssss
\subsection{Peanut compression}\label{sec:peanut}

% \alert{Anna: see what you think. The section is modified after discussion with A-K.}
The above two-body representation uses finely resolved MFS source points for each contact pair. To accelerate evaluation of the correction bases in the far field, we now introduce a recompression step that replaces these fine sources with effective strengths back at the original coarse source points. Fine sources are only retained for evaluations in the vicinity of the particle pair.
%exactly when this is needed will soon be made more precise. % how very vague...

Consider again the $(1\text{-}2)$ pair.
To evaluate $\vec\eta^{(1,2)}[\mm^{(1)}]$ in \eqref{chi} away from the pair, we replace the fine source set $\vec{\mathcal Y}^{(1\text{-}2)}$, with strengths $\vec\beta^{(1,2)} \in \mathbb R^{4N_f}$, by equivalent strengths at the original coarse source points $\vec Y^{(1)}$ and $\vec Y^{(2)}$, whose union we denote by $\vec Y^{(1\text{-}2)}$. These equivalent strengths $\tilde{\vec\lambda}^{(1,2)}$ are added to the
source strengths $\bar{\vec\lambda}^{(1)}$ and $\bar{\vec\lambda}^{(2)}$ stemming from the one-body basis, leaving an FMM cost no more than that of the one-body basis.
As above, it is convenient to combine the source corrections from
$\vec\eta^{(1,2)}$ and $\vec\eta^{(2,1)}$ into
\[
\tilde{\vec{\lambda}}^{(1\text{-}2)}
\coloneqq
\tilde{\vec{\lambda}}^{(1,2)}
+
\tilde{\vec{\lambda}}^{(2,1)},
\]
which will be obtained by applying a precomputed pair correction matrix $\vec A^{(1\text{-}2)}$
of size  $4N_c\times 4N_c$, so that
\begin{equation}\label{correction}
\tilde{\vec{\lambda}}^{(1\text{-}2)}
=
\vec A^{(1\text{-}2)}
\begin{bmatrix}
\bar{\vec \lambda}^{(1)} \\
\bar{\vec \lambda}^{(2)}
\end{bmatrix}.
\end{equation}
%Despite its small size, it is by design fully capturing the lubrication forces for the pair. % *** Not true, only when seen by everyone else. But can we find another good formulation here?
% *** Anna: Research question: Does the map have to be dense? Is it enough to use a subset of the coarse sources? Or, as an alternative, is it best to first build the coarse-to-coarse map A but then only keep its most ``information heavy'' rows/cols via e.g. ID?
% AHB: good ideas, but since Nc~50, not worth it. Or for 3D, a very coarse set could be used for all far-field, but delta* would have to be big O(1). There would be a trade off with finer near-evals. Really you're designing a FDS here... and maybe trying an existing 3D FDS instead would be better (eg FLAM).
% Anna: The thinking here was to see if one can do something to have fewer values at an interpolation step (particularly in 3D). In Dhairya's and Dan's work, the number of entries to be interpolated is much smaller. 
The one-body coarse source vector
$\bar{\vec\lambda}^{(1)}$
for particle 1 is then corrected to give
\begin{equation}\label{total_corr}
\vec\lambda^{(1)}_j
=
\bar{\vec{\lambda}}^{(1)}_j
+
\sum_{k\in \mathcal C^{(1)}}
\tilde{\vec{\lambda}}^{(1\text{-}k)}_j, \quad j = 1,\dots, N_c.
\end{equation}
Once all coarse strength vectors have been corrected in this manner,
the two-body representation \eqref{rep2} is evaluated simply by plain coarse summation
\eqref{stokes_sum}.

 % ***AHB: I replaced this with one sentence above, and clearer, I think... Anna: Ok!
% \begin{remark}[Coarse source strengths]
%Compare \eqref{total_corr} with the one-body source strengths in \eqref{total_lambda}. Away from body~1 and its near neighbors, evaluating $\vec\psi^{(1)}(\vec x)$ via compression no longer requires summation over fine Stokeslets as in \eqref{chi}, but reduces to the coarse evaluation
%\[
%\vec\psi^{(1)}(\vec x)$
%\sum_{j=1}^{N_c}
%\mathbb S(\vec x,\vec y_j^{(1)})
%\vec\lambda_j^{(1)}.
%\]
%\end{remark}

The rest of the subsection is devoted to assembling
the ``coarse-to-coarse'' correction matrix $\vec A^{(1\text{-}2)}$ via a sequence of precomputed linear maps during the setup stage.
This matrix accurately compresses
all of the near-contact lubrication BVP solution information, and is thus
%that needed fine discretization to solve accurately.
analogous to a scattering matrix in the FMPS or FDS frameworks (Remark~\ref{r:hybrid}).
As a reminder for the rest of the paper, all point sets and source types are collected in Tables \ref{tab:sets} and \ref{tab:sources}.

\begin{table}[t] % ttttttttttttttttttttttttttttttttttttttttttttttttttttttttttttttttttt
% ***AHB: figs and tables are clearest at [t] not [h]
  \centering
\begin{tabular}{ l l }
\textit{Point set} & \textit{Description} \\
\hline
\\[-2mm]
$\displaystyle \vec Y^{(\partone)}  = \{ \vec y_\csi^{(\partone)} \}_{\csi=1}^{N_c}$
& Set of coarse source points for particle $\partone$\\
$\displaystyle \vec X^{(\partone)}  = \{ \vec x_\cci^{(\partone)} \}_{\cci=1}^{M_c}$ & Set of coarse collocation points for particle $\partone$\\
$\displaystyle \vec Z^{(\partone\text{-}\parttwo)}=\vec Z^{(\partone)}\bigcup \vec {Z}^{(\parttwo)}$ & Union of point sets for particles $\partone$ and $\parttwo$, where $\vec Z$ can be $\vec Y$, $\vec X$\\
$\displaystyle \vec {\mathcal Y}^{(\partone\text{-}\parttwo)}$  
%= \{ \vec {\mathcal y}_\csi^{(\partone)} \}_{\csi=1}^{N_f}$
& Set of fine source points for particles $\partone$ and $\parttwo$ in near contact\\
$\displaystyle \vec {\mathcal X}^{(\partone\text{-}\parttwo)}$ 
%= \{ \vec {\mathcal x}_\csi^{(\partone)} \}_{\cci=1}^{M_f}$
& Set of fine collocation points for particle $\partone$ and $\parttwo$ in near contact\\
\\[-2mm]
\hline
%$\displaystyle \{\vec X^{(\partone)}\}_{\partone=1}^{P}$ & Full set of coarse source points for all $P$ particles. 
\end{tabular}
\caption{The four different point sets defined on each particle. We also denote $\displaystyle \{\vec X^{(\partone)}\}_{\partone=1}^{P}$ as the full set of coarse source points for all $P$ particles, and similarly for the coarse collocation points. The number of fine source and collocation points used to set up a pair-correction will depend on the particle distance within each pair, see Section \ref{sec:fine}.}
\label{tab:sets}
\end{table}

\begin{table}[t] % tttttttttttttttttttttttttttttttttttttttttttttttttttttttttttttttttt
\centering
\begin{tabular}{ l l }
\textit{Quantity} & \textit{Description} \\
\hline
\\[-2mm]
$\displaystyle \vec\lambda^{(\partone)}\in\mathbb R^{2N_c}$ %= \lbrace \vec\lambda_\csi^{(\partone)}\rbrace_{\csi=1}^{N_c}$
& Coarse source strengths associated with $\vec Y^{(\partone)}$
\\
$\displaystyle \mm^{(\partone)} \in\mathbb R^{2M_c}$ %= \lbrace \mm_\cci^{(\partone)}\rbrace_{\cci=1}^{M_c} 
& Transformed coefficient vector associated with $\vec X^{(\partone)}$ \\
$\displaystyle \bar{\vec\lambda}^{(\partone)} \in\mathbb R^{2N_c}$ % = \lbrace \bar{\vec\lambda}_\csi^{(\partone)}\rbrace_{\csi=1}^{N_c}$
& Obtained from $\vec P^+\vec\mu^{(i)}$, applied as in \eqref{howG+}, associated with $\vec Y^{(\partone)}$
\\
$\displaystyle \vec\beta^{(\partone\text{-}\parttwo)}\in\mathbb R^{4N_f}$ % = \{\vec\beta_\fsi^{(\partone,\parttwo)}\}_{\fsi=1}^{2N_f}$
& Fine source strengths for pair-correction, associated with $\vec {\mathcal Y}^{(\partone\text{-}\parttwo)}$\\
$\displaystyle\tilde{\vec{\lambda}}^{(\partone\text{-}\parttwo)}\in\mathbb R^{4N_c}$ & Correction vector to $[\vec\lambda^{(\partone)};\vec\lambda^{(\parttwo)}]$ on $\vec Y^{(\partone\text{-}\parttwo)}$\\%$=\vec Y^{(\partone)}\bigcup \vec Y^{(\parttwo)}$, \\
& 
$\rightarrow$ see \eqref{peanut} for resistance and \eqref{peanut_mob} for mobility \\
\\[-2mm]
\hline
\end{tabular}
\caption{Coefficient vectors in solution procedures.}
\label{tab:sources}
% \textcolor{red}{A-K: Note: j is used both for point index and particle index as in rest of paper. Change in command def, but also in rest of paper. This thing with $\vec \lambda$ and $\vec {\bar{\lambda}}$ should be sorted out. If we should make a distinction, motivate it clearly where it is introduced and be consistent. Anna: I have tried to change. Is it consistent?}}
\end{table}

\begin{figure}[t] % ffffffffffffffffffffffffffffffffffffffffff
\centering
\includegraphics[trim={5.1cm 15cm 1.5cm 6.7cm},clip,width=0.85\textwidth]{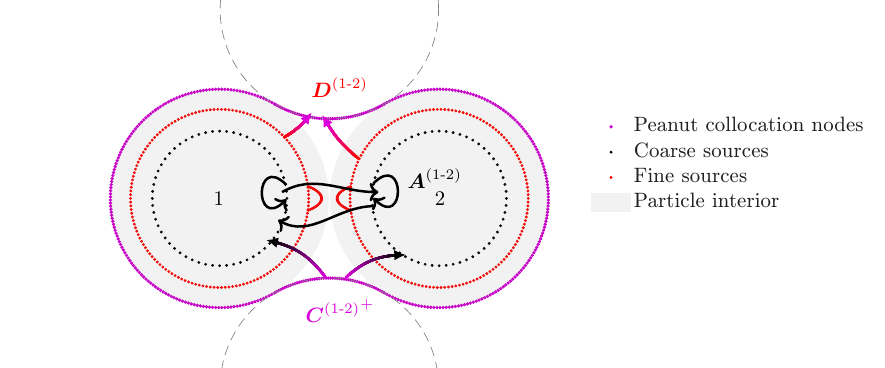}
% ***AHB: is the aspect ratio a bit off? (stretched vertically?) Anna: yes, you really have a sharp eye! Corrected.
  \caption{The peanut %, discretized with magenta nodes, is the union of four radius-1 circular arcs. It 
    (magenta) forms a separation or proxy surface; no other disk can get closer to the ($1\text{-}2$) contacting pair. In Section~\ref{sec:peanut} this is used, along with matrices $\vec C^{(1\text{-}2)}$ and $\vec D^{(1\text{-}2)}$,
    to fill $\vec A^{(1\text{-}2)}$ (black arrows), the coarse-to-coarse matrix that compresses the true near-contact pairwise lubrication interaction.}
  % ***AHB: suggest to cut, defer to the text... Anna: it was an attempt to help the reader, also after disc w A-K...
  %    correction bases $\vec\eta^{(1,2)}$ and $\vec\eta^{(2,1)}$, constructed so that particles outside the peanut see the correct flow field generated by the pair.}
    %Internally, we build $\vec A^{(1\text{-}2)}$ by first determining a set of fine sources by solving two local BVPs---the problem illustrated in Figure \ref{fine_map} and its mirror where particle 2 is the active particle---and then compressing these fine sources to equivalent coarse ones. This last step is done by matching flow fields at collocation nodes on the peanut boundary. Algebraically this is achieved via the matrix product ${\vec C^{(1\text{-}2)}}^+\vec D^{(1\text{-}2)}$, which first maps fine sources to the peanut and then maps the result to the coarse proxy surfaces.}
  \label{peanut_fig}
\end{figure}

The idea of the replacement is to ensure that the coarse sources reproduce the same velocity field as the fine ones on (and therefore exterior to) a special ``peanut'' proxy surface.
This is the separation boundary formed by rolling a unit circle around the fixed pair,
comprising four arcs,
which defines the minimal clearance from any third body; see Figure \ref{peanut_fig}.
  % The replacement is done via a linear map that can be determined in a setup stage for each pair by
  % enforcing that the coarse sources generate the same potential $\uu$
  % at a new set of collocation nodes as the fine sources. The new collocation nodes live on a ``peanut'' separation surface defined by rolling a body around the fixed pair; it is the closest any other particle can get to the pair, see Figure \ref{peanut_fig}.
Let $\z_p^{(1\text{-}2)}$, $p=1,\dots,M_p$ be peanut collocation nodes, uniformly sampled in arc length.
%
%with $2N_c<M_p$, %<2N_f$, % *** Anna: why this constraint: $2N_c<M_p<2N_f$? Has always been here...
%***AHB: Suggest to postpone M_p choice to numer expts.
 As usual with the MFS (or other proxy point applications), we solve for the effective strengths via least-squares velocity collocation at the peanut nodes. That is, given fine strengths $\vec\beta^{(1\text{-}2)}$, we solve for $\tilde{\vec{\lambda}}^{(1\text{-}2)}$ such that
% As usual with the MFS, the matching problem is solved in the least-squares sense. Given $\bb^{(1,2)}$, we solve
% for $\tilde{\vec{\lambda}}^{(1,2)}$ as
\be
\sum_{j=1}^{2N_c} \mathbb S(\z_p^{(1\text{-}2)},\vec Y^{(1\text{-}2)}_j) \tilde{\vec{\lambda}}^{(1\text{-}2)}_j =
\sum_{q=1}^{2N_f} \mathbb S(\z_p^{(1\text{-}2)},\vec{\mathcal Y}^{(1\text{-}2)}_q) \vec\beta_q^{(1\text{-}2)},
\quad p=1,\dots,M_p
\quad \mbox{(peanut compression)},
\label{peanut}
\ee
where, as before, we suppress the $2$-vector nature of each component for notational convenience. Compactly, this is summarized as
$\vec C^{(1\text{-}2)} \tilde{\vec{\lambda}}^{(1\text{-}2)} = \vec D^{(1\text{-}2)}  \bb^{(1\text{-}2)}$, with the linear maps illustrated in Figure \ref{peanut_fig}.
% Other conditions should not be needed;
%   certainly for fine Stokeslets, representation by coarse Stokeslets is
%   consistent (eg see Thm. 5 and comment after in Stein--Barnett QFS paper).
If the residual of this linear solve is small for any data $\bb^{(1\text{-}2)}$, then the
far-field flow generated by the fine sources is faithfully captured by the smaller number of coarse Stokeslets.
% ***AHB: don't think needed:
%No further conditions are required and %: %for instance, if both fine and coarse sources are Stokeslets, then 
%the representation is consistent; see Theorem~5 and the surrounding discussion in the work by Stein \& Barnett \cite{Stein2022}. 

To this end, one precomputes $\vec C^{(1\text{-}2)} = \vec U\vec \Sigma \vec V^T$, then stores the
factors $\vec V \vec \Sigma^+$ and $\vec U^T \vec D^{(1\text{-}2)} $, so that $\tilde{\vec{\lambda}}^{(1\text{-}2)} = {\vec C^{(1\text{-}2)}}^+ \vec D^{(1\text{-}2)} \bb^{(1\text{-}2)} $
is performed via $\tilde{\vec{\lambda}}^{(1\text{-}2)} = \vec V\vec\Sigma^+ \left(\vec U^T \vec D^{(1\text{-}2)} \bb^{(1\text{-}2)}\right)$.
We then combine this factorization with the two-body strength factorization
\eqref{bfm} to give (at a high level, without breaking pseudoinverses into their stable application recipes),
\begin{equation*}
    \begin{aligned}
     \tilde{\vec{\lambda}}^{(1\text{-}2)} & = {\vec C^{(1\text{-}2)}}^+ \vec D^{(1\text{-}2)} \bb^{(1\text{-}2)} = {\vec C^{(1\text{-}2)}}^+ \vec D^{(1\text{-}2)} {\vec F^{(1\text{-}2)}}^+ [-\vec H^{(2,1)} \mm^{(2)};-\vec H^{(1,2)} \mm^{(1)}] = \\ &=
{\vec C^{(1\text{-}2)}}^+ \vec D^{(1\text{-}2)} {\vec F^{(1\text{-}2)}}^+ [-\vec Q^{(2,1)} \vec P^+ \mm^{(2)};-\vec Q^{(1,2)} \vec P^+ \mm^{(1)}] = \\ & =  {\vec C^{(1\text{-}2)}}^+ \vec D^{(1\text{-}2)} {\vec F^{(1\text{-}2)}}^+ [-\vec Q^{(2,1)} \bar{\vec\lambda}^{(2)};-\vec Q^{(1,2)} \bar{\vec\lambda}^{(1)}].   
    \end{aligned}
\end{equation*}
This serves as the correction vector to $[\vec\lambda^{(1)};\vec\lambda^{(2)}]$ due to
$\vec\eta^{(1,2)}$ and $\vec\eta^{(2,1)}$, to be added to $[\bar{\vec\lambda}^{(1)};\bar{\vec\lambda}^{(2)}]$ as per \eqref{total_corr}.
%\coloneqq[\vec P^+ \mm^{(1)}; \vec P^+ \mm^{(2)}]$. 
%In practice, the sum of correction sources $\tilde{\vec{\lambda}}^{(1\text{-}2)}\coloneqq\tilde{\vec{\lambda}}^{(1,2)}+\tilde{\vec{\lambda}}^{(2,1)}$ due to $\vec\eta^{(1,2)}$ and $\vec\eta^{(2,1)}$ is computed jointly.
% The corrected $\vec\lambda$ vectors are then sent to FMM. As
% in Section \ref{sec:standard}, the diagonal one-body contributions are then subtracted off and $\mm$
% added for the identity mapping of the one-body basis.
% The diagonal blocks in the preconditioned system matrix will however no longer be strictly $\vec I$ due to the $\vec\eta$ correction).
Its correction matrix stack extends the stack used to determine $\vec\beta^{(1\text{-}2)}$ (presented at the end of Section \ref{sec:2body}) and has the interpretation (going right-to-left):
\begin{center}
\setlength{\fboxsep}{8pt}
\framebox{%
\begin{minipage}{0.97\textwidth}
\centering \hspace*{-6ex}
\makebox[2cm][c]{%
  \shortstack{%
    {\small Coarse source}\\
    {\small pair correction}\\
    $\tilde{\vec{\lambda}}^{(1\text{-}2)}$%
  }%
}
$\xleftarrow[
  {\shortstack{\vspace{0.3ex}\\[-0.2ex]\scriptsize $\mathbb R^{4N_c\times 2M_p}$}}
]{
  \shortstack{\scriptsize ${\vec C^{(1\text{-}2)}}^+$\\[-0.2ex]\vspace{0.3ex}}
}$
\makebox[1.3cm][c]{\shortstack{\small{Peanut}\\\small{collocation}\\\small{velocity}\\\small{}}}
%$\xleftarrow{\makebox[0.7cm]{\scriptsize$\vec D^{(1\text{-}2)}$}}$
$\xleftarrow[
  {\shortstack{\vspace{0.3ex}\\[-0.2ex]\scriptsize $\mathbb R^{2M_p\times 4N_f}$}}
]{
  \shortstack{\scriptsize ${\vec D^{(1\text{-}2)}}$\\[-0.2ex]\vspace{0.3ex}}
}$
%\makebox[1.3cm][c]{\shortstack{\small{Fine}\\\small{sources}}}
\makebox[1.3cm][c]{\shortstack{\small{Fine}\\ \small{sources}\\$\vec{\beta}^{(1\text{-}2)}$}}
%$\xleftarrow{\makebox[0.7cm]{\scriptsize${\vec F^{(1\text{-}2)}}^+$}}$
$\xleftarrow[
  {\shortstack{\vspace{0.3ex}\\[-0.2ex]\scriptsize $\mathbb R^{4N_f\times 4M_f}$}}
]{
  \shortstack{\scriptsize ${\vec F^{(1\text{-}2)}}^+$\\[-0.2ex]\vspace{0.3ex}}
}$
\makebox[2.2cm][c]{\shortstack{\hspace*{-4ex}\small{Fine}\\ \hspace*{-1ex}\small{collocation}\\\small{velocity}}} \hspace*{1ex}
$\xleftarrow[
  {\shortstack{\vspace{0.5ex}\\[-0.2ex]\scriptsize $\mathbb R^{4M_f\times 4N_c}$}}
]{%
  \makebox[1.0cm][c]{%
    \scriptsize
    $\mathord{-\left[\!\!
      \begin{array}{cc}
        \arraycolsep=1pt
        \vec 0 & \! \vec Q^{(2,1)} \\
        \vec Q^{(1,2)} & \! \vec 0
      \end{array}
    \!\!\right]}$%
  }%
}$
% \makebox[1.9cm][c]{\hspace*{2.5ex}\shortstack{\small{Coarse}\\\small{sources}}} 
\hspace*{2ex}
\makebox[2.0cm][c]{\shortstack{\small{Coarse}\\ \small{sources}\\$\bar{\vec\lambda}^{(1)};\bar{\vec\lambda}^{(2)}$}} \hspace*{-10ex}
% *** Anna: I'm leaving this last bit out.
% $\xleftarrow{\makebox[1.0cm]{%
% \scriptsize$\mathord{\left[\!\begin{array}{cc}\arraycolsep=1pt
% \vec P^+ & \vec 0 \\
% \vec 0 &  \vec P^+
% \end{array}\!\right]}$}}$
% \makebox[1.9cm][c]{\shortstack{\small{Coarse}\\\small{colloc.}\\$\vec\mu^{(1)};\vec\mu^{(2)}$}}
\end{minipage}
}% end framebox
\end{center}

%Each pseudoinverse here should generally be split into two factors for stable application,
%but that would require excessive storage. We have found empirically that one can safely group all factors together into the matrix $\vec A^{(1\text{-}2)}$, which we are now finally ready to define as

The above stack of matrices gives the coarse-to-coarse correction matrix formula
\begin{equation}\label{Amat}
    \vec A^{(1\text{-}2)} = -{\vec C^{(1\text{-}2)}}^+\vec D^{(1\text{-}2)}{\vec F^{(1\text{-}2)}}^+\begin{bmatrix} \vec 0 & \vec Q^{(2,1)} \\
 \vec Q^{(1,2)} &  \vec 0 \end{bmatrix}.
\end{equation}
As expected, its form resembles a block Schur complement when eliminating the fine source degrees of freedom.
Since $\vec A^{(1\text{-}2)}$ is well-conditioned, once is has been filled
(using backward-stable applications of pseudoinverses along the stack), the
intermediate matrices may be discarded, avoiding excessive storage.
% ***AHB: what do you think?
% Anna: Good!

%  such that the additive update to the coarse source vector due to the contacting pair $(1\text{-}2)$  is given by \begin{equation}\label{correction}
% \tilde{\vec{\lambda}}^{(1\text{-}2)} = 
%     \vec A^{(1\text{-}2)}\begin{bmatrix}\bar{\vec \lambda}^{(1)} \\ \bar{\vec \lambda}^{(2)} \end{bmatrix}.
% \end{equation} 
We have now described how the corrections are done for a single contact pair. In the global matrix-vector apply, it remains to loop over all contacts to correct the entire source vector and store the result in $\vec\lambda$; see \eqref{total_corr}.

There are two clear benefits of peanut compression, compared to doing pair corrections as in Section~\ref{sec:2body}: the number of sources to be fed in to an FMM is smaller and the required amount of storage is reduced. We will in Section \ref{sec:res_solve} describe how the compressions are handled within the resistance solve. 

\begin{remark}[Proxy compression]
  Compression using collocation points that hug a particle pair has independently been developed in work to appear by Martinez Aguilar, Malhotra \& Fortunato on solving elliptic PDEs in the exterior of close-to-touching disks. The idea of compressing numerous unknowns to fewer equivalent proxy sources is common to fast direct solvers
  \cite{Martinsson2005,Martinsson2019}
  % AHB: stein+barnett did not invent this ... hence I cite earlier... Anna: thanks.
  and kernel-independent FMMs \cite{Malhotra2015}.
\end{remark}

%%%ssssssssssssssssssssssssssssssssssssssssssss
\subsection{Extracting forces and torques}
In a resistance problem, the final output is particle net forces and torques, $\lbrace\vec f^{(i)},\,t^{(i)}\rbrace_{i=1}^P$. 
Following Gauss' law for Stokes, all Stokeslet sources in the interior of a particle contribute to these quantities. For particle 1, the force and torque pair $(\vec f^{(1)},\,t^{(1)})$ can be determined from the computed coarse sources $\bar{\vec \lambda}^{(1)}$ stemming from the one-body basis, together with the fine sources $\vec\beta^{(1,k)}$ associated with all close neighbors $k$ of particle 1: 
	\begin{equation}\label{force_eq}
		\vec f^{(1)} =\sum_{i = 1}^{N_c}\bar{\vec\lambda}_i^{(1)}+\sum_{k\in\mathcal C^{(1)}}\sum_{i=1}^{N_f}\vec\beta^{(1,k)}_i,\quad t^{(1)} = \sum_{i = 1}^{N_c} \bar{\vec\lambda}_i^{(1)}\cdot (\vec y_i^{(1)}-\vec c^{(1)})^{\perp} + \sum_{k\in\mathcal C^{(1)}}\sum_{i=1}^{N_f}\vec\beta_i^{(1,k)}\cdot(\vec{\mathcal  Y}_i^{(1\text{-}k)}-\vec c^{(1)})^{\perp}.
	\end{equation} 
 With more compact notation, \eqref{force_eq} can be written as
 \begin{equation}\label{force_syst}
		\left[{\vec f^{(1)}}^T,t^{(1)}\right]^T  = \vec K^T\bar{\vec\lambda}^{(1)} + \sum_{k\in\mathcal C^{(1)}}{\vec{\mathcal K}^{(1,k)}}^T\vec \beta^{(1,k)},
        \end{equation}
        with the rigid body matrices
        \begin{equation}\label{Kmat}
        \vec K = \begin{bmatrix} \vec I_2 & (\vec y^{(1)}_1-\vec c^{(1)})^{\perp}\\ \vec I_2 & (\vec y^{(1)}_2-\vec c^{(1)})^{\perp} \\ \vdots & \vdots \\ \vec I_2 & (\vec y^{(1)}_{N_c}-\vec c^{(1)})^{\perp}
        \end{bmatrix}\in \mathbb R^{2N_c\times 3}\text{ and } \vec{\mathcal K}^{(1,k)} = \begin{bmatrix} \vec I_2 & (\vec{\mathcal Y}^{(1\text{-}k)}_1-\vec c^{(1)})^{\perp}\\ \vec I_2 & (\vec{\mathcal Y}^{(1\text{-}k)}_2-\vec c^{(1)})^{\perp} \\ \vdots & \vdots \\ \vec I_2 & (\vec{\mathcal Y}^{(1\text{-}k)}_{N_f}-\vec c^{(1)})^{\perp}
        \end{bmatrix}\in \mathbb R^{2N_f\times 3}.
	\end{equation}     
        The second term in \eqref{force_syst} can be viewed as a correction to the net forcing due to all near contacts. We assume that the coarse matrix $\vec K$ is equal for every disk,
        %(the system is monodisperse and the same resolution),
        and hence add no superscript. The fine matrices $\vec{\mathcal K}^{(1,k)}$ however depend on $\vec{\mathcal Y}^{(1\text{-}k)}$, $k\in\mathcal C^{(1)}$. To avoid storing the submatrices needed to stably reconstruct the fine sources $\vec\beta^{(1,k)}$, $k\in\mathcal C^{(1)}$, we instead store the much smaller matrices $\vec E_{ft}^{(1,k)}\in\mathbb R^{6\times 4N_c}$ that for the $(1\text{-}k)$ pair extracts the correction to the net forces and torques, given the coarse one-body sources only. This means for the $(1\text{-}2)$ near contact with no other close neighbors that 
\begin{equation}\label{two_force}
    \begin{bmatrix}
        \vec f^{(1)}\\ t^{(1)} \\ \vec f^{(2)}\\ t^{(2)}
    \end{bmatrix} = \begin{bmatrix}
        \vec K & \vec 0 \\ \vec 0 & \vec K
    \end{bmatrix}\begin{bmatrix}
    \bar{\vec \lambda}^{(1)} \\ \bar{\vec\lambda}^{(2)}\end{bmatrix} +\vec E_{ft}^{(1\text{-}2)}\begin{bmatrix} \bar{\vec\lambda}^{(1)} \\ \bar{\vec\lambda}^{(2)}
    \end{bmatrix},
\end{equation}
with
\begin{equation}\label{force_ext}
\vec E_{ft}^{(1\text{-}2)} = \begin{bmatrix} \vec{\mathcal K}^{(1,2)} & \vec 0 \\ \vec 0 & \vec{\mathcal K}^{(2,1)} \end{bmatrix}{\vec F^{(1\text{-}2)}}^+\begin{bmatrix} \vec 0 & \vec Q^{(2,1)} \\
 \vec Q^{(1,2)} &  \vec 0 \end{bmatrix}.
 \end{equation}
% The representation in \eqref{force_syst} is used to determine forces and torques in the absence of peanut compression. When peanut compression is applied, there is a linear map...
%Forces and torques are returned as output from the solve outlined in the pseudocode in Algorithm \ref{alg:peanut}.
\begin{remark}[Fine sources are needed]
  It is not possible to recover the net forces and torques on the individual particles from the
  coarse sources computed via peanut compression. Although these sources reproduce the same far-field flow and the same total force for the pair, they do \emph{not} in general preserve the forcing on each particle \emph{separately}.
  % AHB: this must be that forces/torques are split differently between the 2 bodies,
  % since mathematically the flow in the far field can only carry one net force and net torque.
  % Do you agree? Clarify for reader?
  % Anna: I modified the description. At infinity, only the total force should be the same for the pair system (not the torque). I tested  numerically and indeed this seems to be what happens; the coarse equivalent sources give the correct total force for the pair.
  % AHB: cool!
\end{remark}

%%%%%%%%5ssssssssssssssssssssssssssssssssssssssssss
\subsection{The full resistance algorithm}\label{sec:res_solve}
Pseudocode summarizing the full resistance solve is given in Algorithm~\ref{alg:peanut}. It assumes a uniform discretization of coarse proxy and collocation point sets per body, takes as input the stacked vector of rigid body velocities evaluated at all coarse collocation nodes, and returns net forces and torques on all bodies. More specifically, we solve the
two-body collocation system
\begin{equation}
    \sum_{i=1}^P \vec \psi^{(i)}[\mm^{(i)}](\x)=\vec g(\x),
    \qquad \vec x \in \{\vec X^{(i)}\}_{i=1}^{P},
\end{equation}
and use internal one-body coarse %and fine
% ***AHB correct? Anna: No. I modified, but we can probably write it shorter.  Only coarse sources are needed explicitly.
sources to recover forces and torques (via the force/torque correction matrix for each close pair).
In the pseudocode, a boolean flag \texttt{flowEval} indicates whether the full sets of coarse and fine sources, 
\begin{equation}\label{source_types}
(\vec\lambda,\bar{\vec\lambda},\lbrace\vec\beta^{(i\text{-}k)}, k\in\mathcal C^{(i)},i = 1,\dots,P\rbrace),
\end{equation}
are returned, to enable flow field evaluations both far from and near to the particles in a separate routine. Returning the triplet in \eqref{source_types} instead of merely $\vec\mu$ enables direct flow field evaluations, e.g.~via the FMM.
% ***AHB: I think the reader has seen Q already:
%The action of the matrices $\vec Q^{(i,k)}$, needed to  retrieve the fine sources, is exemplified in Figure \ref{fine_map}. Anna: ok

The detailed matrix--vector apply, needed in each GMRES iteration, is outlined separately in Algorithm~\ref{alg:matvec}. In each call, the flow field due to the corrected sources is evaluated globally via the FMM, and then modified locally on a per-pair basis, following the general structure of Algorithm \ref{alg:general}. For the $(i\text{-}k)$ pair, the coarse source correction $\tilde{\vec\lambda}^{(i\text{-}k)}$ cannot replace the fine sources when evaluating the field on the pair itself. Hence, the contribution from $\tilde{\vec\lambda}^{(i\text{-}k)}$ should in principle be subtracted off and the flow field due to the fine sources $\vec\beta^{(i\text{-}k)}$ added back in. This would however mean that fine sources would need to be retrieved, at extra cost. We can avoid this by replacing the flow field due to fine sources by the evaluation of the one-body basis on the neighboring particle, which is precisely the quantity that the fine sources are constructed to match; see \eqref{local_BVP} and Step~2b of Algorithm~\ref{alg:matvec}. For a reminder of the action of the matrices $\vec H^{(i,k)}$ and $\vec P$, see Figure \ref{fine_map}. As a last step, the one-body blocks must be corrected, as with the one-body basis. This procedure ensures that the $4N_c\times 4N_c$ two-body system matrix obtained for an isolated contact pair %when evaluating the two-body basis at the coarse collocation nodes
is approximately the identity.

 This completes the interpretation of the peanut-compressed two-body MFS scheme as a
1-level fast direct solver coupled to an iterative solver for the coarse unknowns.

% The replacement of fine sources by coarse ones cannot be used to evaluate the field on the pair itself. At the solve stage, we can elevate the issue by equating the evaluation of fine sources by the one-body basis on the neighboring particle that the fine sources are constructed to match; see \eqref{local_BVP} and Step 2b in Algorithm \ref{alg:matvec}.

% The pair corrections are constructed so that the preconditioned two-body system matrix for an isolated contact pair is approximately the identity. With formulae, letting $\vec T^{(1\text{-}2)}$ denote the evaluation of the two-body basis at the coarse collocation nodes for the pair in isolation, we have
% \[
%     \vec T^{(1\text{-}2)}\begin{bmatrix}
%         \vec \mu^{(1)} \\ \vec \mu^{(2)} 
%     \end{bmatrix}
%     =
%     \begin{bmatrix} \vec b^{(1)} \\ \vec b^{(2)} \end{bmatrix},
% \]
% with $\vec T^{(1\text{-}2)} \approx \vec I$. To avoid excessive storage, we therefore simply replace $\vec T^{(1\text{-}2)}$ by $\vec I$ when accounting for the self-interaction of the contact pair.  see Step 2b in Algorithm \ref{alg:matvec}. Practically, this is done by equating the evaluation of fine sources by the one-body basis on the neighboring particle that the fine sources are constructed to match; see \eqref{local_BVP}.

\begin{remark}[Fast pair corrections] 
In practice, the loop over close pairs in Algorithm \ref{alg:matvec} can be replaced by two applications of precomputed sparse matrices, with nonzero entries only for coarse sources associated with close pairs. The first matrix maps one-body coarse source strengths $\bar{\vec\lambda}$ to corrected coarse source strengths $\vec\lambda$, replacing the loop over pairs in Step~1b. The second matrix determines the corresponding correction to the flow field, ensuring block-diagonal identity contributions for each pair, and replaces the loop over pairs in Step 2b. 
\end{remark}

\begin{algorithm}
  % custom comment defs
 \algrenewcommand\algorithmiccomment[2][\footnotesize]{{#1\hfill\(\triangleright\) #2}}  
 \algnewcommand{\LeftComment}[1]{\Statex \hspace*{1.6ex}\(\triangleright\) #1}
\caption{Fast resistance solve with two-body preconditioning and peanut compression}\label{alg:peanut}
\textbf{Global data:} Proxy point sets $\{ \vec X^{(i)} \}_{i=1}^P$, collocation point sets $\{ \vec Y^{(i)} \}_{i=1}^P$, boolean flag \texttt{flowEval}
\begin{algorithmic}

\Function{solve}{$\vec g$}
    \State \textbf{Input:} Stacked right-hand side surface velocity data vector $\vec g$, 
    \State \textbf{Output:} Particle net forces and torques $\lbrace \vec f^{(i)},t^{(i)}\rbrace_{i=1}^P$, if \texttt{flowEval} then stacked coarse source \hspace*{3.5ex}strengths (force density) vectors $\vec\lambda$, $\bar{\vec\lambda}$, fine source strengths $\lbrace \vec \beta^{(i\text{-}k)}\rbrace$, fine source points $\lbrace \vec{\mathcal Y}^{(i\text{-}k)}\rbrace$, $k\in\mathcal C^{(i)}$,\, \\ \quad\hspace*{0.5ex} $i =1,\dots,P$
    %\If{\texttt{flowEwal}} stacked coarse source strengths (force density) vector $\vec\lambda$, fine source strengths $\lbrace \vec \beta_{ij}\rbrace$, fine source points $\lbrace \vec{\mathcal Y}^{(i\text{-}k)}\rbrace$ \EndIf
\LeftComment{\emph{Precomputations:} Identify all $C$ contacts (particle pairs with interparticle distance $<\delta_c$). }
\For{$c = 1$ to $C$}
\State Get particles $i,k$ in contact $c$, set fine discretization  $\vec{\mathcal X}^{(i\text{-}k)}$, $\vec{\mathcal Y}^{(i\text{-}k)}$ and peanut nodes $\lbrace\vec z^{(i\text{-}k)}_q\rbrace_{q=1}^{M_p}$
\State Compute correction matrices and store:
\If{\texttt{flowEval}}\\ \vspace*{-4ex}\State Coarse-to-fine mapping via SVD of $\vec F^{(i\text{-}k)}$: $\vec V^{(i\text{-}k)}{\vec\Sigma^{(i\text{-}k)}}^+$ and ${\vec U^{(i\text{-}k)}}^T\begin{bmatrix} \!\vec 0 &\!\! \vec Q^{(k,i)} \\ 
 \vec Q^{(i,k)} & \!\! \vec 0\! \end{bmatrix}$,\\
 \Comment{for local post-processing; see eqn.~\eqref{bfm}} 
 \EndIf
 \vspace*{-1ex}
 \State Coarse-to-forcing correction matrices $\vec E_{ft}^{(i\text{-}k)}$
\Comment{see eqn.~\eqref{force_ext}}
\State Coarse-to-coarse correction matrices $\vec A^{(i\text{-}k)}$
\Comment{peanut compression; see eqn.~\eqref{Amat}}
\EndFor

\LeftComment{\emph{Solve} for stacked surface values $\vec\mu$ using matrix-vector multiply function defined in Algorithm \ref{alg:matvec}, which\\ \hspace*{3ex} gets access to all local variables} %of the \texttt{SOLVE} function
    \State $\vec\mu \gets \mathrm{GMRES}(\texttt{MATVEC}, \vec g)$
    \LeftComment{\emph{Postprocessing:} Determine net forces and torques and if \texttt{flowEval} then prepare flow field evaluation:\\ \hspace*{3.5ex}Recover coarse strengths for each body via local pseudoinverse and determine forcing contribution:}
    \For{$i = 1$ to $P$}
     
        \State $\bar{\vec\lambda}^{(i)} \gets %({\vec S^{(kk)}})^{+} \bga^{(k)}$ \Comment{By \eqref{reconbal}}
        \vec P^+ %\vecmc V\vec \Sigma^+ (\vecmc U^T 
        \vec\mu^{(i)}$
        \Comment{apply pseudoinverse backward-stably as in eqn.~\eqref{howG+}}
        \State $\begin{bmatrix} \vec f^{(i)} \\ t^{(i)}\end{bmatrix} \gets \vec K^T\bar{\vec\lambda}^{(i)}$ \Comment{see first term in eqn.~\eqref{two_force}}
    \EndFor
    Determine fine force/torque correction for each contact and if \texttt{flowEval} then recover fine source vectors
    \For{$c=1$ to $C$}
    \State Get particle indices $i,k$ for contact $c$
    \State  $\begin{bmatrix} {\vec f^{(i)}}^T & t^{(i)} & {\vec f^{(k)}}^{T} & t^{(k)}  \end{bmatrix}^T\gets \begin{bmatrix} {\vec f^{(i)}}^T & t^{(i)} & {\vec f^{(k)}}^{T} & t^{(k)}  \end{bmatrix}^T+\vec E_{ft}^{(i\text{-}k)}\begin{bmatrix}\bar{\vec\lambda}^{(i)} \\ \bar{\vec\lambda}^{(k)}
    \end{bmatrix}$ \Comment{see second term in \\ \vspace*{-2.6ex} \hspace*{124.5ex} eqn.~\eqref{two_force}}
    \If{\texttt{flowEval}}\State
    $\vec\beta^{(i\text{-}k)} = \vec V^{(i\text{-}k)}{\vec\Sigma^{(i\text{-}k)}}^+\left({\vec U^{(i\text{-}k)}}^T\begin{bmatrix} \vec 0 & \vec Q^{(k,i)} \\
 \vec Q^{(i,k)} &  \vec 0 \end{bmatrix}\begin{bmatrix} \bar{\vec \lambda}^{(i)} \\ \bar{\vec \lambda}^{(k)} \end{bmatrix}\right)$\Comment{see eqn.~\eqref{bfm}} 
  \EndIf
    \EndFor
    \If{\texttt{flowEval}}
    \State Call \texttt{MATVEC}$(\vec\mu)$ in Algorithm \ref{alg:matvec} to recover $\vec\lambda$ \EndIf
    \State \Return $\lbrace \vec f^{(i)},t^{(i)}\rbrace_{i=1}^P$, \textbf{if} \texttt{flowEval} \textbf{then} $\vec\lambda$, $\bar{\vec\lambda}$, $\lbrace 
     \vec\beta^{(i\text{-}k)}\rbrace$, $\lbrace \vec{\mathcal Y}^{(i\text{-}k)}\rbrace$, $k\in\mathcal C^{(i)},\, i =1,\dots,P$ 
\EndFunction
\end{algorithmic}
\end{algorithm}

%%%%%%%%%%%%%%%%%%%%%%%%%%%%%%%%%%%%%%%%%%%%%%%%%%%%%%%%%
\begin{algorithm}[th]
  % custom comment defs
 \algrenewcommand\algorithmiccomment[2][\footnotesize]{{#1\hfill\(\triangleright\) #2}}  
 \algnewcommand{\LeftComment}[1]{\Statex \hspace*{1.6ex}\(\triangleright\) #1}
\caption{Matrix-vector apply for use in the fast resistance solve of Algorithm \ref{alg:peanut}}\label{alg:matvec}
\begin{algorithmic}
\Function{matvec}{$\vec \mu$}
    \State \textbf{Input:} Stacked coarse collocation data vector $\vec \mu = \{ \vec\mu^{(i)}\}_{i=1}^P$, coarse-to-coarse correction matrices \hspace*{3.7ex}$\lbrace\vec A^{(i\text{-}k)}\rbrace$ for all $C$ close pairs of particles $(i\text{-}k)$
    \State \textbf{Output:} Stacked surface velocity vector $\vec u = \{ \vec u^{(i)} \}_{i=1}^P$
    \LeftComment{Step 1a: Recover proxy source strengths for each body via local pseudoinverse apply and add to total:}
    \For{$i = 1$ to $P$}
      %  \State Apply local pseudoinverse $\vec S^{(kk)}$ to get proxy source strengths:
        \State $\bar{\vec\lambda}^{(i)} \gets %({\vec S^{(kk)}})^{+} \bga^{(k)}$ \Comment{By \eqref{reconbal}}
        \vec P^+\vec\mu^{(i)}$ 
        \State  $\vec\lambda^{(i)} \gets \bar{\vec\lambda}^{(i)}$ \Comment{leave $\bar{\vec\lambda}$ for correction at the end}
    \EndFor
    \LeftComment{Step 1b: Apply two-body corrections for all close pairs:}
    \For{$c = 1$ to $C$}
    \State $\tilde{\vec{\lambda}}^{(i\text{-}k)} \gets \vec A^{(i\text{-}k)}\begin{bmatrix}
        \bar{\vec\lambda}^{(i)}\\
        \bar{\vec\lambda}^{(k)} 
    \end{bmatrix}$, \quad
$\begin{bmatrix}
        \vec\lambda^{(i)}\\
        \vec\lambda^{(k)} 
    \end{bmatrix}\gets \begin{bmatrix}
        \vec\lambda^{(i)}\\
        \vec\lambda^{(k)} 
    \end{bmatrix}+\tilde{\vec{\lambda}}^{(i\text{-}k)}$, for particles $i,k$ in contact $c$ \Comment{see eqn.~\eqref{correction}}
    \EndFor
    %\State
    \LeftComment{Step 2a: Fast potential evaluation at all $PM$ targets from all $PN$ sources:}
%with 
% $\vec u^{k}_i = \sum\limits_{k'=1}^P\sum\limits_{j=1}^N\mathbb S\left(\vec x_i^{(k)},\vec y_j^{(k')}\right)\alpha_j^{(k')}, \quad i =1,\dots,N
    \State $\vec u \gets \texttt{FMM\_evaluate}(\{ \vec X^{(i)} \}_{i=1}^P,\{ \vec Y^{(i)} \}_{i=1}^P, \lbrace\vec\lambda^{(i)}\rbrace_{i=1}^P) $
    \Comment{applies bare MFS matrix in \eqref{Gmat}}
    %\State
    \LeftComment{Step 2b: Locally correct to convert two-body-to-two-body blocks to the identity:}
    \For{$c = 1$ to $C$}
        \State  $\begin{bmatrix}\vec u^{(i)} \\ \vec u^{(k)} \end{bmatrix}\gets \begin{bmatrix}\vec u^{(i)} \\ \vec u^{(k)} \end{bmatrix} - \begin{bmatrix} \vec P & \vec S^{(ik)} \\ \vec S^{(ki)} & \vec P\end{bmatrix}\tilde{\vec\lambda}^{(i\text{-}k)}
        -\begin{bmatrix}
            \vec 0 & \vec H^{(k,i)} \\ \vec H^{(i,k)} & \vec 0
        \end{bmatrix}\begin{bmatrix}
            \vec\mu^{(i)} \\ \vec\mu^{(k)}
        \end{bmatrix}$  
        \Comment{subtract local pair-contribution\\ \vspace*{-2.6ex} \hspace*{93.5ex} from compressed coarse sources and\\ \vspace*{0ex} \hspace*{93.5ex} add back right hand side from local\\ 
        %\vspace*{-2.6ex}
        \hspace*{93.5ex} BVPs}
    \EndFor
    \vspace*{-2ex}
    \For{$i=1$ to $P$} 
    \State $\vec u^{(i)} \gets  \vec u^{(i)}-\vec P\bar{\vec \lambda}^{(i)}+\vec\mu^{(i)}$ 
        \Comment{correct one-body identities; see eqn.~\eqref{onesub}}
    \EndFor
    \State \Return $\vec u$
\EndFunction
\end{algorithmic}
\end{algorithm}

%***AHB: I find it mysterious that 2-body then 1-body Id block corrections are needed! Anna: You mean the loop over bodies at the end? Two-body is the one-body + corrections and we want the one-body part to produce identities. (where corrections should be zero),

%mmmmmmmmmmmmmmmmmmmmmmmmmmmmmmmmmmmmmmmmmmmmmmmmmmmmmmmmmmmmmmmm
\section{Solving the mobility problem via MFS}\label{sec:mob}

We now show how the same type of basis representation can be used to efficiently solve a mobility problem. As a reminder, the forces and torques, $\lbrace \vec f^{(i)},t^{(i)}\rbrace_{i=1}^P$, are prescribed and 
% unknown source strengths are related to the prescribed  via (compare \eqref{two_force})
% \begin{equation}\label{force_torque}
%     \left[{\vec f^{(i)}}^T,t^{(i)}\right]^T = \vec K^T\vec\lambda^{(i)},\medspace i=1,\dots, P.
% \end{equation}
particle velocities, $\lbrace\vec v^{(i)}, \omega^{(i)}\rbrace_{i=1}^P$, are unknown in this setting. The coarse and fine point sets and their respective sources follow the conventions introduced in Tables \ref{tab:sets} and Table \ref{tab:sources}.

%%%sssssssssssssssssssssssssssssssssssssssssssssssssssssssssss
\subsection{One-body representation}
\label{s:mob_1b}
We begin with a dilute suspension, where a set of coarse sources per body suffices to resolve all interactions, and specify the representation in \eqref{mob_compflow}. In Section \ref{sec:one_mob}, we will see that also the mobility one-body basis functions lead to a system equal to that stemming from one-body preconditioning, as was the case for the one-body basis functions used for resistance.

We first focus on particle $i$ and construct the corresponding one-body basis functions. The slip boundary condition \eqref{mobility_bc} will by design be satisfied exactly at the collocation nodes. %***Anna: word choice: "exactly" is probably not what I want to say 
We therefore now first determine $\vec\phi^{(i)}[\vec\mu^{(i)}](\vec x)$, $\vec x\notin \vec X^{(i)}$ and soon return to the case $\vec x\in\vec X^{(i)}$.
The discretized slip boundary condition of \eqref{mobility_bc} can in the one-body problem be written as ${\vec s}^{(i)} = \vec B\begin{bmatrix} \vec v^{(i)} \\ \omega^{(i)} \end{bmatrix} +\vec\mu^{(i)}$. Here, the matrix $\vec B$ maps rigid body velocities to boundary velocity. It is defined analogously to $\vec K$ in \eqref{Kmat}, but has size $2M_c \times 3$, and is formed by replacing $\vec y_j^{(1)}$, $j=1,\dots,N_c$, with $\vec x_j^{(1)}$, $j=1,\dots,M_c$. Using the representation in \eqref{stokes_sum} with $P=1$, we formulate the constrained least-squares problem $\vec S\hat{\vec\lambda} = \vec s$ subject to $\vec f^{(i)} = \vec 0$, $t^{(i)} = 0$. To avoid enforcing these constraints explicitly, we use a ``recompleted'' formulation as in \cite{Broms2024c}. The idea is to construct $\vec\phi^{(i)}$ so that the constraints are automatically satisfied. The force and torque constraints then read $\vec K^T \hat{\vec\lambda}^{(i)} = \vec 0$. This is satisfied by choosing $\hat{\vec\lambda}^{(i)} = (\vec I - \vec L)\bar{\vec\lambda}^{(i)}$, where $\vec L = \vec K(\vec K^T \vec K)^{-1}\vec K^T \in \mathbb R^{2N_c \times 2N_c}$ projects onto rigid body motions and where $\bar{\vec\lambda}$ emphasizes that these are the one-body source strengths, as in the resistance setting. We thus represent $\vec\phi^{(i)}$ as 
\begin{equation}\label{phi_eq}
    \vec\phi^{(i)}[\vec\mu^{(i)}](\vec x)  = \mathbb S(\vec x,\vec Y^{(i)})\left(\vec I-\vec L\right)\bar{\vec\lambda}^{(i)},\quad \vec x\notin \vec X^{(i)},
\end{equation}
where the dependence on $\vec\mu^{(i)}$ is encoded in the coefficient vector $\bar{\vec\lambda}^{(i)}$, whose explicit form is derived below.
% ***AHB: writing \phi without [\mm] is a bit weird, since \phi alone has been defined in Sec 2 as
% a solution *operator*,
% not a basis function (flow field). After arxiv, you might want to rewrite this section less in the "telling a story" and more as "define the thing used and comment as to why later" (the latter mode is clearer!).
%
An unused subspace can be exploited as an ansatz for the unknown rigid body velocities, allowing us to express them as
\begin{equation}\label{rbm_ansatz}
\begin{bmatrix} \vec v^{(i)} \\ \omega^{(i)} \end{bmatrix} =  -\vec K^T\bar{\vec\lambda}^{(i)}.
\end{equation}
 The unknown boundary data at the collocation nodes can then be written as
 \begin{equation}\label{bc_data}
 \vec s^{(i)} = -\vec B \vec K^T\bar{\vec\lambda}^{(i)}+\vec\mu^{(i)}.
 \end{equation}
 One solves for $\bar{\vec\lambda}^{(i)}$ by evaluating \eqref{phi_eq} at $\vec X^{(i)}$ and  matching the result with \eqref{bc_data} such that 
 \begin{equation}\label{SL}
 \bar{\vec\lambda}^{(i)} = \vec S_L^+\vec\mu^{(i)},
 \end{equation}
  where $\vec S_L = \vec P(\vec I-\vec L)+\vec B\vec K^T$.
 % \begin{equation}
 %     \vec\phi^{(i)}(\vec x)  = \mathbb S(\vec x,\vec Y^{(i)})\left(\vec I-\vec L\right)\vec S_L^{+}\vec\mu^{(i)},
 % \end{equation}
 At the collocation nodes, we 
 %define $\vec\phi^{(i)}$ so that the BVP \eqref{mobility_bc} is satisfied exactly and 
 express the unknown rigid body velocity pair $(\vec v^{(i)},\omega^{(i)})$ as linear functionals of $\vec\mu^{(i)}$ via \eqref{rbm_ansatz} and \eqref{SL}.  Hence, $\vec\phi^{(i)}$ can be expressed as
\begin{equation}\label{1body_mob}
    \vec\phi^{(i)}\coloneqq \vec\phi^{(i)} [\mm^{(i)}](\x) := \sum_{j=1}^{M_c} \vec\mu_j^{(i)} \vec\phi^{(i)}_j(\x)
   = \begin{cases} \vec \mu_m^{(i)} - \left(\vec B\vec K^T \vec S_L^+\vec\mu^{(i)}\right)_m, \quad \vec x=\vec x_m^{(i)}, \\
   \sum\limits_{n=1}^{N_c} \vec S(\x,\y^{(i)}_n)\left[\left(\vec I-\vec L\right) (\vec S_L^{+}\mm^{(i)})\right]_n, \quad \text{otherwise.}\end{cases}
\end{equation}
This specifies the rigid-body matrices in \eqref{mobility_bc} as
${\cal R}^{(i)} = -\vec B\vec K^T \vec S_L^+$, for all $i$.

By construction, $\vec\phi^{(i)}$ carries zero net force and torque on particle $i$.
To represent flows with non-zero prescribed force and torque, we add a known completion flow to the sum of one-body basis functions, as in \eqref{mob_compflow}.
In the present MFS setting, we construct it using Stokeslet sources at the coarse source points,
\begin{equation}
    \vec C^{(i)}\left[\vec f^{(i)},t^{(i)}\right](\vec x) = \mathbb S(\vec x,\vec Y^{(i)})\vec\lambda_0^{(i)}, \qquad i=1,\dots,P.
\end{equation}
Its source strengths $\vec\lambda_0^{(i)}$ are set per particle to sum exactly to the prescribed force and torques:  $\vec K^T\vec\lambda_0^{(i)} = \vec F$ (compare \eqref{force_syst}), determined via the ansatz $\vec\lambda_0 = \vec K\vec a$, for some rigid body velocity components stacked in $\vec a\in\mathbb R^3$. This completes the description of the representation in \eqref{mob_compflow}.

% This formulation will also be the foundation for the mobility choice of one-body and two-body basis functions. The basic idea is to split the source vectors $\vec\lambda^{(i)}$, $i=1,\dots,P$, each into two terms, one in the null-space of the constraint matrix $\vec K^T$, and one in the orthogonal complement to this space, i.e.~in the image space of $\vec K$. The first term is obtained via a projection, by multiplying with the projection matrix $\vec I-\vec L$, , while the second term is obtained by combining the ansatz $ \vec\lambda_0^{(i)} = \vec K\vec F$ with \eqref{force_torque}.
 % \begin{equation}
 % \begin{aligned}
 %     \vec\lambda_0^{(i)} &= \vec K_c\vec F,\\
 %     \vec K_c^T\vec\lambda_0^{(i)} &=\begin{bmatrix} \vec f^{(i)} \\ t^{(i)} \end{bmatrix}.
 % \end{aligned}
 % \end{equation}

%%%%%%%%%%%%%%%%%%%%%%%%%%%%%%%%%%%%%%%%%%%%%%%%%%
\subsection{Equivalence with one-body preconditioning}
\label{sec:one_mob}
This subsection explains the equivalence of the above to our prior one-body preconditioning method for mobility \cite{Broms2024c}, and connects it to the general formulation of
Section~\ref{s:mob_gen}.
Written in terms of internal proxy source strengths $\lbrace \bar{\vec\lambda}^{(i)}\rbrace_{i=1}^P$ instead of boundary unknowns $\lbrace \vec\mu^{(i)}\rbrace_{i=1}^P$, the representation in \eqref{mob_compflow} may be expressed as
\begin{equation}\label{mob}
     \vec u(\vec x) = \sum_{i=1}^P\mathbb S(\x,\vec Y^{(i)})\left[\left(\vec I-\vec L\right)\bar{\vec\lambda}^{(i)}+\vec\lambda_0^{(i)}\right].
 \end{equation}
This was the representation of the flow field in \cite{Broms2024c}, upon a split of the source vector into one component  in the null-space of the constraint matrix $\vec K$ and one component that sums to the net force and torque on each particle. Evaluating \eqref{mob} at the collocation nodes $\lbrace\vec X^{(i)}\rbrace_{i=1}^P$, expressing the boundary data as $\vec g^{(i)} = -\vec B\vec K^T\bar{\vec\lambda}^{(i)},\,i=1,\dots, P$, and reordering terms, the system to solve for the mobility problem takes the form 
\begin{equation}\label{tosolve}
    \begin{bmatrix}
         \vec S^{(11)} \left(\vec I-\vec L\right)+\vec B\vec K^T & \vec S^{(12)} \left(\vec I-\vec L\right)& \dots & \vec S^{(1P)} \left(\vec I-\vec L\right)  \\
        \vec S^{(21)} \left(\vec I-\vec L\right) & \vec S^{(22)} \left(\vec I-\vec L\right)+\vec B\vec K^T  & \dots & \dots \\
        \vdots & \vdots & \ddots & \vdots \\
        \vec S^{(P1)} \left(\vec I-\vec L\right) & \dots & \dots & \vec S^{(PP)} \left(\vec I-\vec L\right)+\vec B\vec K^T  
    \end{bmatrix}
    \begin{bmatrix}
        \bar{\vec \lambda}^{(1)} \\  \bar{\vec \lambda}^{(2)} \\ \vdots \\ \bar{\vec \lambda}^{(P)}
        \end{bmatrix} = 
        \begin{bmatrix}
        -\vec u_0^{(1)} \\  -\vec u_0^{(2)} \\ \vdots \\ -\vec u_0^{(P)}
    \end{bmatrix},
\end{equation}
with the completion flow surface velocity data
\begin{equation}\label{u0_def}
    \vec u_0^{(i)} = \sum_{k=1}^P\mathbb S(\vec X^{(i)},\vec Y^{(k)})\vec\lambda_0^{(k)}.
\end{equation} This large unconstrained least-squares problem can be solved by applying the same type of one-body preconditioning from the right as was applied for the resistance problem in Section \ref{sec:standard}. The result is
%preconditioned system becomes
    \begin{equation}\label{2esysprec}
    \begin{bmatrix}
        \vec I & \vec S^{(12)} \left(\vec I-\vec L\right)\vec S_L^+ & \dots & \vec S^{(1P)} \left(\vec I-\vec L\right)\vec S_L^+  \\
        \vec S^{(21)} \left(\vec I-\vec L\right)\vec S_L^+ & \vec I & \dots & \dots \\
        \vdots & \vdots & \ddots & \vdots \\
        \vec S^{(P1)} \left(\vec I-\vec L\right)\vec S_L^+ & \dots & \dots & \vec I
    \end{bmatrix}
    \begin{bmatrix}
        \vec \mu^{(1)} \\  \vec \mu^{(2)} \\ \vdots \\ \vec \mu^{(P)}
        \end{bmatrix} = 
        \begin{bmatrix}
        -\vec u_0^{(1)} \\  -\vec u_0^{(2)} \\ \vdots \\ -\vec u_0^{(P)}
    \end{bmatrix},
\end{equation}
    where the diagonal blocks have been replaced by $\vec I$ to regularize the system, as in the resistance setting.
    The above linear system is precisely the one described in Section \ref{s:mob_gen} when
    using one-body basis functions that obey \eqref{mobility_bc}, recalling the above
    definition of ${\cal R}^{(i)}$.
% ***AHB I don't think the cancellation needs to be explained any more.

%    It is easy to check that 
%This is the same system as obtained by solving for $\lbrace\vec\mu^{(i)}\rbrace_{i=1}^P$ from the representation in \eqref{mob_compflow}, by matching the flow field at the coarse collocation nodes to the unknown rigid body motion, and using the ansatz \eqref{rbm_ansatz} based on internal source strengths for the unknown rigid body velocities. (The ansatz decouples and can therefore be used also in this multi-particle setting.) 

% The regularization manifests as a modification of $\vec\phi^{(i)}$ so that 
%Following the definition of $\vec\phi^{(i)}$ in \eqref{1body_mob}, the term $\left(-\vec B\vec K^T \vec S_L^+\vec\mu^{(i)}\right)$ is cancelled by the expression for the unknown rigid body data at the collocation nodes and therefore leaves identities along the diagonal of the system matrix.

    % We recover \eqref{2esysprec} from \eqref{evaluate_mob}, up to reordering of terms, by collocating at the coarse nodes and imposing the boundary data $\vec b^{(i)} = -\vec B\vec K^T\vec\lambda^{(i)}$, where $\vec\lambda^{(i)} = \vec S_L^{+}\mm^{(i)}$.

    %%%%%%%sssssssssssssssssssssssssssssssssssssssssssss
\subsection{Two-body preconditioning}\label{sec:mob_two}
% The flow field in \eqref{evaluate_mob}-\eqref{1body_mob} is written as the sum of two terms: one that contributes no net force or torque, and a second that accounts for the net forcing through a completion flow, following Remark \ref{mob_remark}. 

%as in the resistance setting.
% Consider first a system consisting of two circles only. The flow field can be expressed as 
% \begin{equation}\label{2rep_mob}
%     \vec u(\vec x) = \left[\vec \phi^{(1)}\vec\mu^{(1)}\right](\vec x) + \left[\vec \eta^{(1,2)}\vec\mu^{(1)}\right](\vec x) +\left[\vec \phi^{(2)}\vec\mu^{(2)}\right](\vec x) + \left[\vec \eta^{(2,1)}\vec\mu^{(2)}\right](\vec x) + \SR ^{(1)}\vec\lambda_0^{(1)} + \SR ^{(2)}\vec\lambda_0^{(2)}.
% \end{equation}
% The two-body correction $\vec \eta^{(1,2)}$ %to the one-body basis of particle 1 due to the near contact with body 2 
% is constructed such that its internal source strengths do not contribute to the net force and torque on body 1. 
% We need to assure that the internal source strengths of the two-body correction $\vec\eta^{(1,2)}$ exert no net force or torque on either body 1 or 2.
The two-body correction $\vec\eta^{(1,2)}$ is constructed so that its internal sources exert no net force or torque on either body 1 or 2.
As in the one-body case \eqref{1body_mob}, this is achieved by a projection---now applied to the fine source strengths on both particles in the pair. The correction basis takes the form
\be
   \vec \eta^{(1,2)}[\mm^{(1)}](\x) = \sum_{q=1}^{2N_f} \mathbb S(\x,\vec{\mathcal Y}_q^{(1\text{-}2)})\left[(\vec I-\vec{\mathcal L}^{(1\text{-}2)}) \vec \beta^{(1,2)}\right]_q,
   \label{chi_mob}
\ee
with the fine pair projection matrix 
\begin{equation}
    \vec{\mathcal L}^{(1\text{-}2)} = \begin{bmatrix} \vec{\mathcal L}^{(1,2)} & \vec 0 \\ \vec 0 & \vec{\mathcal L}^{(2,1)} \end{bmatrix},\text{ where } \vec{\mathcal L}^{(i,k)}  \coloneqq  \vec{\mathcal K}^{(i,k)} \left({\vec{\mathcal K}^{(i,k)}}^T\vec{\mathcal K}^{(i,k)}\right)^{-1}{\vec{\mathcal K}^{(i,k)}}^T
\end{equation}
and $\vec{\mathcal K}^{(i,k)}$ is the fine rigid body matrix of \eqref{Kmat}.

Equation \eqref{chi_mob} leaves an unused image space of ${\vec{\mathcal K}^{(1,2)}}^T$ and ${\vec{\mathcal K}^{(2,1)}}^T$, which we now exploit to close the system---together with the corresponding image space of $\vec K^T$ from the one-body basis.
% In \eqref{chi_mob}, there is hence an unused image space of the constraint matrix $\vec{\mathcal K}^T$ that may be used to close the system, together with the unused image space of $\vec K^T$ appearing in the one-body basis. 
For the contact pair $(1\text{-}2)$, assumed to be isolated from other particles, we make the ansatz 
% \begin{equation}
%     \begin{aligned}
%          = -\vec K_c^T\vec\lambda^{(1)} -\vec{\mathcal K}^T\left(\vec\beta^{(1,2)}_1+\vec\beta^{(2,1)}_1\right),\\
%         \vec U^{(2)} = -\vec K_c^T\vec\lambda^{(2)} -\vec{\mathcal K}^T\left(\vec\beta^{(1,2)}_2+\vec\beta^{(2,1)}_2\right),      
%     \end{aligned}
% \end{equation} 
\begin{equation}\label{rbm2}
    \begin{aligned}
         \begin{bmatrix} \vec v^{(1)} \\ \omega^{(1)} \\  \vec v^{(2)} \\ \omega^{(2)} \end{bmatrix}  = -\begin{bmatrix} {\vec{\mathcal K}^{(1,2)}}^T &\vec 0 \\ \vec 0 & {\vec{\mathcal K}^{(2,1)}}^T \end{bmatrix}\left(\vec\beta^{(1,2)}+\vec\beta^{(2,1)}\right)-\begin{bmatrix} \vec K^T &\vec 0 \\ \vec 0 & \vec K^T \end{bmatrix} \begin{bmatrix} \bar{\vec\lambda}^{(1)} \\ \bar{\vec\lambda}^{(2)}\end{bmatrix}      
    \end{aligned}
\end{equation} 
for the unknown rigid body velocities.
At the particle boundaries, the no-slip velocity field at the coarse collocation nodes is, as before, expressed as  
\[
\vec g^{(1)} = \vec B 
\begin{bmatrix} \vec v^{(1)} \\ \omega^{(1)} \end{bmatrix}, \qquad
\vec g^{(2)} = \vec B 
\begin{bmatrix} \vec v^{(2)} \\ \omega^{(2)} \end{bmatrix},
\]
which gives the combined form:
\begin{equation}
\begin{bmatrix} \vec g^{(1)} \\ \vec g^{(2)} \end{bmatrix} 
= - 
\begin{bmatrix}
    \vec B{\vec{\mathcal K}^{(1,2)}}^T & \vec 0 \\
    \vec 0 & \vec B{\vec{\mathcal K}^{(2,1)}}^T
\end{bmatrix}
\left(\vec\beta^{(1,2)} + \vec\beta^{(2,1)}\right)
-
\begin{bmatrix}\label{fine_bc}
    \vec B\vec K^T \vec S_L^+ & \vec 0 \\
    \vec 0 & \vec B\vec K^T \vec S_L^+
\end{bmatrix}
\begin{bmatrix} \vec\mu^{(1)} \\ \vec\mu^{(2)} \end{bmatrix}.
\end{equation}

We now match the velocity field produced by \eqref{rep2mob} with $P=2$ to the unknown no-slip boundary data in \eqref{fine_bc}. After reordering terms, the resulting equation to be satisfied at the  coarse collocation nodes $\vec X^{(1\text{-}2)}$ becomes:

% \begin{equation}\label{2body_re}
% \begin{aligned}
% \sum_{k=1}^{N_c} 
% &\left[ 
% \mathbb S(\vec x,\y^{(1)}_k)(\vec I - \vec L) + \vec B_c\vec K_c^T 
% \right] 
% (\vec S_L^+ \vec\mu^{(1)})_k \\
% + \sum_{k=1}^{N_c} 
% &\left[ 
% \mathbb S(\vec x,\y^{(2)}_k)(\vec I - \vec L) + \vec B_c\vec K_c^T 
% \right] 
% (\vec S_L^+ \vec\mu^{(2)})_k \\
% + \sum_{q=1}^{2N_f} 
% &\left[ 
% \G(\vec x,\vec{\mathcal Y}_q^{(1\text{-}2)})(\vec I - \vec L_f^{\text{pair}}) 
% + 
% \begin{bmatrix}
% \vec B_c\vec{\mathcal K}^T & \vec 0 \\
% \vec 0 & \vec B_c\vec{\mathcal K}^T
% \end{bmatrix}
% \right] 
% \left( \vec\beta_q^{(1,2)} + \vec\beta_q^{(2,1)} \right) \\
% &= - \sum_{k=1}^{N_c} \mathbb S(\vec x,\y^{(1)}_k) \vec\lambda_0^{(1)}_k
%    - \sum_{k=1}^{N_c} \mathbb S(\vec x,\y^{(2)}_k) \vec\lambda_0^{(2)}_k.
% \end{aligned}
% \end{equation}

\begin{equation}\label{2body_re}
\begin{aligned}
  & 
  \sum\limits_{i=1}^2\vec\phi^{(i)}[\vec\mu^{(i)}]\left(\vec X_m^{(1\text{-}2)}\right)+\sum_{j=1}^{M_c}\begin{bmatrix}
        \vec B\vec{K}^T\vec S_L^+ & \vec 0 \\ \vec 0 & \vec B\vec{K}^T\vec S_L^+
    \end{bmatrix}_{mj}\begin{bmatrix} \vec\mu^{(1)} \\ \vec\mu^{(2)}\end{bmatrix}_j\\
  % \sum_{k=1}^{N_c}\left(\sum_{j=1}^{N_c}\mathbb S(\vec X^{(1\text{-}2)}_i,\y^{(1)}_j)\left(\vec I-\vec L\right)_{jk}+\left[\vec B\vec K^T\right]_{ik} \right)(\vec S_L^{+}\mm^{(1)})_k \\ 
  &+    
  % \sum_{k=1}^{N_c}\left(\sum_{j=1}^{N_c}\mathbb S(\vec X^{(1\text{-}2)}_i,\y^{(2)}_j)\left(\vec I-\vec L\right)_{jk}+\left[\vec B\vec K^T\right]_{ik} \right)(\vec S_L^{+}\mm^{(2)})_k \\&+
   \sum_{n=1}^{2N_f} \left(\sum_{j=1}^{2N_f}\mathbb S(\vec X^{(1\text{-}2)}_m,\vec{\mathcal Y}_j^{(1\text{-}2)})(\vec I-\vec{\mathcal L}^{(1\text{-}2)})_{jn}+\begin{bmatrix}
        \vec B{\vec{\mathcal K}^{(1,2)}}^T & \vec 0 \\ \vec 0 & \vec B{\vec{\mathcal K}^{(2,1)}}^T
    \end{bmatrix}_{mn}\right) \vec\beta_n^{(1,2)} \\ &+
   \sum_{n=1}^{2N_f} \left(\sum_{j=1}^{2N_f}\mathbb S(\vec X^{(1\text{-}2)}_m,\vec{\mathcal Y}_j^{(1\text{-}2)})(\vec I-\vec{\mathcal L}^{(1\text{-}2)})_{jn}+\begin{bmatrix}
        \vec B{\vec{\mathcal K}^{(1,2)}}^T & \vec 0 \\ \vec 0 & \vec B{\vec{\mathcal K}^{(2,1)}}^T
    \end{bmatrix}_{mn}\right) \vec\beta_n^{(2,1)}  = \\ &  \hspace*{60ex} = - \sum_{i=1}^2\vec C^{(i)} \left[\vec f^{(i)},t^{(i)}\right](\vec X^{(1\text{-}2)}_m),  m = 1,\dots, 2M_c.
   \end{aligned}
\end{equation}
% The first two terms on the left-hand side arise from the one-body bases $\vec\phi^{(1)}$ and $\vec\phi^{(2)}$, along with their associated contributions to the rigid body velocities at the particle surfaces. 
The block-diagonal matrix appearing in the second term cancels the same-body contribution $\vec B\vec K^T\vec S_L^+\vec\mu^{(i)}$ appearing in $\vec\phi^{(i)}[\vec \mu^{(i)}]$; see \eqref{1body_mob}.
Terms three and four represent the two-body correction bases $\vec\eta^{(1,2)}$ and $\vec\eta^{(2,1)}$ and their associated contributions to the no-slip boundary velocity.
% The first two terms in the left-hand side arise from $\vec\phi^{(1)}$ and $\vec\phi^{(2)}$, written together with the one-body contributions to the unknown rigid body velocities. The third and fourth terms correspond to the correction bases and the associated contributions to the rigid body velocities for the interacting pair. 
It remains to express $\vec \beta^{(1,2)}$ and $\vec\beta^{(2,1)}$ as the solutions to least-squares problems involving $\vec\mu^{(1)}$ and $\vec\mu^{(2)}$.  We do so for  the sum $\vec \beta^{(1\text{-}2)}\coloneqq \vec\beta^{(2,1)}+\vec\beta^{(1,2)}$ jointly by thinking about the two-body system in \eqref{2body_re} in the form
\begin{equation}
    \vec T^{(1\text{-}2)}\begin{bmatrix}
        \vec \mu^{(1)} \\ \vec \mu^{(2)} 
    \end{bmatrix} = \begin{bmatrix} -\vec u_0^{(1)} \\ -\vec u_0^{(2)} \end{bmatrix}, 
\end{equation}
with $\vec u_0^{(k)} = \sum_{i=1}^2\vec C^{(i)}\left[\vec f^{(i)}, t^{(i)}\right](\vec X^{(k)})$. To ensure that the system matrix fulfills $\vec T^{(1\text{-}2)}\approx\vec I$---so that the system for more than two particles is efficiently preconditioned---the one-body contributions on the neighboring particle have to be cancelled out by the two-body corrections, mirroring the approach taken for the resistance problem. Hence, the sum
$\vec\beta^{(1\text{-}2)}$ is chosen to satisfy
\begin{equation}
\begin{aligned}
\sum_{n=1}^{2N_f} \left(\sum_{j=1}^{2N_f}\mathbb S(\vec{\mathcal X}^{(1\text{-}2)}_m,\vec{\mathcal Y}^{(1\text{-}2)}_j)(\vec I-\vec{\mathcal L}^{(1\text{-}2)})_{jn}+\begin{bmatrix}
        \vec{\mathcal B}{\vec{\mathcal K}^{(1,2)}}^T & \vec 0 \\ \vec 0 & \vec{\mathcal B}{\vec{\mathcal K}^{(2,1)}}^T
    \end{bmatrix}_{mn}\right)\vec\beta_n^{(1\text{-}2)}%\left(\vec\beta_k^{(1\text{-}2)}+\vec\beta_k^{(2,1)}\right)
\;=\\\;
=\left\{\begin{array}{ll}
-\vec \phi^{(2)}[\mm^{(2)}](\vec{\mathcal X}^{(1\text{-}2)}_m), & m = 1,\dots,M_f \qquad \mbox{(on body 1)}\\
-\vec \phi^{(1)}[\mm^{(1)}](\vec{\mathcal X}^{(1\text{-}2)}_m), & m = M_f+1,\dots,2M_f  \qquad \mbox{(on body 2)}.\\
\end{array}\right.
\end{aligned}
\end{equation}
The second term on the left-hand side is a correction to the no-slip boundary velocities, where $\vec{\mathcal B}$ is the fine counterpart of $\vec B$, constructed using fine collocation points.
With this relation for $\vec\beta^{(1\text{-}2)}$, the two-body basis functions $\vec\psi^{(1)}$ and $\vec\psi^{(2)}$ solve \eqref{corr_syst_mob}. 
%sssssssssssssssssssssssssssssssssssssssssssssssssssssssssssssssss
\subsection{Peanut compression}

It remains to express the correction fields $\vec\eta^{(1,2)}$ and $\vec\eta^{(2,1)}$ using only coarse sources. As in the resistance case \eqref{peanut}, the fine source strengths for a near-contact pair are replaced by an equivalent coarse representation obtained through a least-squares match of the fine and coarse velocity fields on the peanut separation surface, now with both representations constrained to produce zero total force and torque on the pair.
% To accelerate this step, we compress the fine source sum $\vec\beta^{(1\text{-}2)}$—which accounts for both $\vec\eta^{(1,2)}$ and $\vec\eta^{(2,1)}$—into a single coarse source vector $\tilde{\vec\lambda}^{(1\text{-}2)}$. 
Specifically, we solve
\begin{equation}
\sum_{j=1}^{2N_c} \mathbb S(\z_p^{(1\text{-}2)},\tilde{\vec Y}^{(1\text{-}2)}_j)\left[\left(\vec I-\vec L^{\text{pair}}\right)\tilde{\vec{\lambda}}^{(1\text{-}2)}\right]_j =
\sum_{q=1}^{2N_f} \mathbb S(\z_p^{(1\text{-}2)},\vec{\mathcal Y}^{(1\text{-}2)}_q)\left[(\vec I-\vec{\mathcal L}^{(1\text{-}2)}) \vec\beta^{(1\text{-}2)}\right]_q,
\quad p=1,\dots,M_p,
\label{peanut_mob}
\end{equation}
where $\vec L^{\text{pair}}$ is block-diagonal with two copies of $\vec L$ on the diagonal. The full source vector is formed as in \eqref{total_corr} and the flow evaluated in the far field as (compare \eqref{mob})
\begin{equation}
     \vec u(\vec x) = \sum_{i=1}^P\mathbb S(\x,\vec Y^{(i)})\left[\left(\vec I-\vec L\right)\vec\lambda^{(i)}+\vec\lambda_0^{(i)}\right].
 \end{equation}

The solve procedure mirrors that of the resistance problem (Algorithms~\ref{alg:peanut}–\ref{alg:matvec}), with the following key differences:
\begin{itemize}
    \item The input right-hand side is a surface velocity field representing the completion flow $-\vec u_0^{(i)}$ on particle $i$, $i =1,\dots,P$.
   \item The output consists of translational and angular velocities $\{\vec v^{(i)}, \omega^{(i)}\}_{i=1}^P$ instead of net forces and torques, which are returned alongside the source strengths used for flow evaluation. The velocities are computed in a post-processing step using a procedure closely resembling that used to determine the net forcing (compare \eqref{rbm2} and \eqref{two_force}). 
    \item The linear maps $\vec H^{(1,2)}$, $\vec Q^{(1,2)}$, $\vec F^{(1\text{-}2)}$, $\vec C^{(1\text{-}2)}$, $\vec D^{(1\text{-}2)}$ and $\vec A^{(1\text{-}2)}$ illustrated in Figures \ref{fine_map} and \ref{peanut_fig} are defined via projected evaluations: 
    \begin{equation}
    \begin{aligned}
        \vec Q_{mn}^{(1,2)} & \coloneqq \sum\limits_{j=1}^{N_c}\mathbb  S(\vec{\mathcal X}^{(1\text{-}2)}_{m+M_f}, \y^{(1)}_j)(\vec I - \vec L)_{jn},&\in\mathbb R^{2M_f\times 2N_c}, \\
        \vec H^{(1,2)} & \coloneqq \vec Q^{(1,2)} \vec S_L^+,&\in\mathbb R^{2M_f\times 2M_c},\\ \vec F^{(1\text{-}2)}_{mn} &\coloneqq \left(\sum\limits_{j=1}^{4N_f}\mathbb S(\vec{\mathcal X}^{(1\text{-}2)}_m,\vec{\mathcal Y}^{(1\text{-}2)}_j)\left(\vec I-\vec{\mathcal L}^{(1\text{-}2)}_{jn}\right)+\begin{bmatrix}
        \vec B{\vec{\mathcal K}^{(1,2)}}^T & \vec 0 \\ \vec 0 & \vec B{\vec{\mathcal K}^{(2,1)}}^T \end{bmatrix}\right)_{mn},&\in\mathbb R^{4M_f\times 4N_f},\\
        \vec C^{(1\text{-}2)}_{pn} &\coloneqq \sum\limits_{j=1}^{2N_c}\mathbb S(\vec{z}^{(1\text{-}2)}_p,\vec{\tilde Y}^{(1\text{-}2)}_j)(\vec I-\vec{L}^{\text{pair}})_{jn},&\in\mathbb R^{2M_p\times 4N_c},\\
    \vec D^{(1\text{-}2)}_{pn} & \coloneqq \sum\limits_{j=1}^{4N_f}\mathbb S(\vec{z}^{(1\text{-}2)}_p,\vec{\mathcal Y}^{(1\text{-}2)}_j)(\vec I-\vec{\mathcal L}^{(1\text{-}2)})_{jn},&\in\mathbb R^{2M_p\times 4N_f},\\ \vec A^{(1\text{-}2)} & \coloneqq -{\vec C^{(1\text{-}2)}}^+\vec D^{(1\text{-}2)}{\vec F^{(1\text{-}2)}}^+\begin{bmatrix} \vec 0 & \vec Q^{(2,1)} \\
 \vec Q^{(1,2)} &  \vec 0 \end{bmatrix}, & \in\mathbb R^{4N_c\times 4N_c}.
    \end{aligned}
    \end{equation}
\end{itemize}

%sssssssssssssssssssssssssss
\section{The fine MFS representation for near-contact pairs}\label{sec:fine}
Recall that the $N_c$ coarse MFS sources are equispaced on the concentric circle of radius $R_c$; the $M_c$ coarse collocation nodes are equispaced on the disk boundary, and we will
fix $M_c= 1.2 N_c$ in this paper.
The rest of this subsection describes
the MFS sources and collocation nodes used for the local fine pairwise BVPs.
This is inspired by the method of images of Cheng \& Greengard \cite{Cheng1998} and extends our earlier 3D approach for spheres \cite{Broms2025} to 2D disks.

\begin{figure}[t] % fffffffffffffffffffffffffffffffffffffffffffffffffffffffffffffffffffff
    \centering
            	\includegraphics[trim={3cm 11.3cm 2cm 12.6cm},clip,width=1.05\textwidth]{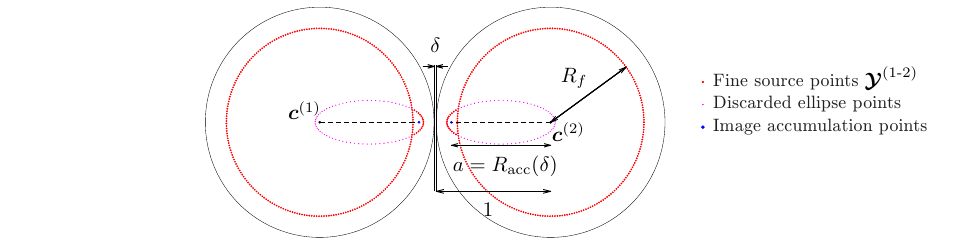}	
   \caption{The proxy source points used for the local fine BVPs are enhanced with sources on ellipse segments that shield singularities located up to a distance $R_{\text{acc}}(\delta)$ from each particle's center.}
   \label{fine_disc}
\end{figure}

%Similarly to \cite{Broms2025},
The pairwise fine source sets $\vec{\mathcal Y}^{(i\text{-}k)}$ include a fixed set of $N_{fp}$ fine ``proxy'' points per disk, equispaced on the concentric circle of radius $R_f$.
Close enough pairs additionally carry sources adapted to the pair contact separation $\delta$, lying on two arcs designed to enclose %analytically
known
singularities in the analytic continuation of the exterior BVP solution into the disk interiors.
See Figure~\ref{fine_disc}. % ***AHB: note pointing the reader here asap is best.
These singularities arise from successive image
reflections of the centers through the disk boundaries
(in the Laplace case see, e.g., \cite{Cheng1998}), giving an infinite series
along the line connecting the two centers, that accumulates at radius
%As in \cite{Broms2025}, these adapted point sets are based on an image construction. Specifically, we utilize that successive reflection across the boundaries of two neighboring unit disks results in an image accumulation point within each body at a radius \cite{Cheng1998,Broms2025}
% For a point $\vec x$ in the exterior of a unit circle centered at $\vec c^{(i)}$, the reflection of $\vec x$ in the circle gives the \emph{image point}
% \begin{equation}\label{reflect}
% 	\vec x' = \vec c^{(i)}+ \frac{(\vec x-\vec c^{(i)})}{\|\vec x-\vec c^{(i)}\|^2}.
% \end{equation}
% For two disjoint unit circles, centered at $\vec c^{(1)},\vec c^{(2)}$ and separated by a small distance $\delta$, successive reflections in one circle followed by the other converge to fixed points:
% \begin{equation}\label{acc_points}
% \begin{aligned}
% 	\vec x_{\text{acc}}^{(1)} &=\vec c^{(1)}+R_{\text{acc}}\frac{\vec c^{(2)}-\vec c^{(1)}}{\|\vec c^{(2)}-\vec c^{(1)}\|},\\ 
% 	\vec x_{\text{acc}}^{(2)} &=\vec c^{(2)}+R_{\text{acc}}\frac{\vec c^{(1)}-\vec c^{(2)}}{\|\vec c^{(2)}-\vec c^{(1)}\|},
% \end{aligned}
% \end{equation}
% with
\begin{equation}
	R_\text{acc}(\delta) = 1+ \delta/2 - \sqrt{\delta + \delta^2/4}.
\end{equation}
The formula is the same for 2D and 3D. In the Stokes case
(see \cite{Broms2025}) it is
expected that the entire line from center to radius $R_\text{acc}$ is singular.
%This accumulation radius is not merely a geometric curiosity; it marks the location of a singularity in the analytic continuation of the PDE solution for two close-to-touching particles.
In \cite{Broms2025}, to handle this singularity we distributed a mixture of fundamental solution types along the outer pieces of the lines, clustered toward the accumulation points.
Instead here we present a scheme based solely on Stokeslets:
for pairs close enough that $R_{\text{acc}}(\delta) > R_f$, the additional Stokeslets lie on
elliptical arcs that ``shield'' the line singularities. The shielding property is necessary for strengths to remain of bounded magnitude
as the MFS solution converges \cite[Thm~2.4,~Ch.~4]{Doicu2000}
(also see \cite[Conj.~11]{Barnett2007}), required for numerical stability.

The shielding elliptical arcs are shown in Figure~\ref{fine_disc}, where both centers lie
on the $x$-axis.
In body-local coordinates, the ellipse associated with particle~1 has foci at $x=0$ and $x=a$, while its rightmost point is located at $x=a+\gamma(1-a)$, where $0<\gamma<1$ is a tunable parameter. % ***AHB gamma is the eccentricity? Anna: I don't think that's the same thing? 
%
%Associating the origin of the complex plane with the center of disk 1.
In complex notation centered at the center of disk 1,
% ***AHB this is still a bit vague. Associating (x,y) with x+iy, ... maybe?
the ellipse is parametrized by
\[
z(\theta) = \frac{a}{2} [1+\cos(\theta+i\alpha)], \qquad
%c+s\cos(\theta+i\alpha),
%\qquad
%c=s=\frac a2, \qquad     % ***AHB no need for c,s since never used later? Anna: Ok!
\mbox{ where } \quad \alpha=\operatorname{arcosh}
\left(1+\frac{2\gamma(1-a)}{a}\right).
\]
We sample the ellipse on the half-shifted angular grid
$$
\theta_j=\frac{(j+1/2)\pi}{N_{\rm ell}},
\qquad
j=0,\dots,2N_{\rm ell}-1,
$$
which clusters nodes toward the tips of the ellipse. We then retain only the nodes satisfying $|z(\theta_j)|>R_f$. This gives an ellipse segment lying outside the fine proxy circle with a high density of nodes in a neighborhood of the image accumulation points. If
$R_{\rm acc}(\delta)\le R_f$, we do not add any enhancing nodes, since
singularities are already shielded by the existing fine proxy circle.

For the pair $(1,2)$ with unit center--center direction \(e^{i\varphi}\), the
enhancing source sets in complex notation are
\[
\vec{\mathcal Y}^{(1,2)}_{\rm enh}
=
\left\{
\vec c^{(1)}+\begin{bmatrix}\Re\lbrace e^{i\varphi}z(\theta_j)\rbrace \\ \Im\lbrace e^{i\varphi}z(\theta_j)\rbrace \end{bmatrix} :\ |z(\theta_j)|>R_f
\right\},
\qquad
\vec{\mathcal Y}^{(2,1)}_{\rm enh}
=
\left\{
\vec c^{(2)}-\begin{bmatrix}\Re\lbrace e^{i\varphi}z(\theta_j)\rbrace \\ \Im\lbrace e^{i\varphi}z(\theta_j)\rbrace \end{bmatrix} :\ |z(\theta_j)|>R_f
\right\}.
\]
This pair's full fine source set is then
\[
\vec{\mathcal Y}^{(1\text{-}2)}
=
\left\{\vec c^{(1)}+R_f\begin{bmatrix} \cos(t_i) \\ \sin(t_i) \end{bmatrix}_{i = 1}^{N_{fp}}\right\}
\cup
\vec{\mathcal Y}^{(1,2)}_{\rm enh}
\cup
\left\{
\vec c^{(2)}+R_f\begin{bmatrix} \cos(t_i) \\ \sin(t_i) \end{bmatrix}_{i = 1}^{N_{fp}}\right\}
\cup
\vec{\mathcal Y}^{(2,1)}_{\rm enh}.
\]
%and is exemplified in Figures \ref{fine_map}, \ref{peanut_fig} and \ref{fine_disc}. 
% In Figure \ref{fine_disc}, the ellipses and image accumulation points are also indicated.
Since at most $2N_{\text{ell}}$ ellipse nodes are sampled, one always has
$
N_f \le N_{fp} + 2N_{\text{ell}}
$. In practice, $R_f$ is chosen so that $
N_f \le N_{fp} + N_{\text{ell}}
$.
Based on empirical performance, we set $\gamma = 0.3$.

%***AHB note I found it unclear here, rewrote...
The local two-body MFS solve requires a set of $2M_f$ boundary collocation nodes, $\vec{\mathcal{X}}^{(i\text{-}k)}$. For this, we first place $M_f = \lceil a_f N_{fp} \rceil$ equispaced nodes on each disk boundary.
%***AHB I deleted a_c since it came out of the blue and this section is about fine not coarse. Anna: ok!
This grid is then augmented by additional pair-adapted collocation nodes obtained through a fixed Möbius angle map. Specifically, we sample $M_{\text{ell}} = a_{\text{ell}}N_{\text{ell}}$ angles $\theta_j$ uniformly and reparameterize them according to $ \tilde{\theta}_j = \pi+\arg\left(\frac{e^{i\theta_j}+r}{1+re^{i\theta_j}}\right). $ This concentrates nodes where the boundary data vary most rapidly, without refining the entire geometry. Throughout the paper, we use $a_f = 1.2$, $r = 0.6$ and $a_{\text{ell}} = 5$.
% ***AHB: N_ell is not given ?
% Also a_ell=5 seems huge (>1e3 fine colloc nodes?) %Anna: this is what I used...

% *** Anna: Remove, right? % The performance of the fine discretization for an isolated pair is demonstrated in Figure \ref{fig:fine_test1}.
The above fine MFS parameters were chosen via convergence studies, that, since they are not
directly relevant to the proposed compressed two-body preconditioning method, will be reported elsewhere.
Recall that in this method, the fine MFS matrices are used only in precomputations;
the main iterative solution involves only coarse surface unknowns ($2M_c$ unknowns per body).

\begin{remark}[Unknowns with one-body preconditioning]
With one-body preconditioning, each body must carry the full fine boundary discretization required to resolve all of its near-contact interactions, leading to very large linear systems. As in our earlier work on spheres \cite{Broms2025}, this remains true even when using an image-enhanced discretization, essential for resolving near-contact pair interactions efficiently. Without such enhancement, dramatically more source and collocation points would be required, especially in three dimensions.
\end{remark}

\begin{remark}[Combining fine sources] 
        In the evaluation of the flow field, where the fine grid is needed, we combine the fine proxy source strengths together from multiple contacts to speed up the computations.
  This is easy since the fine proxy grids do not rotate; they will be coincident. The fine
grid from ellipse segments, however, cannot be combined. 
\end{remark}

% \begin{remark}[Column scaling]
% We do column scaling as in the literature on the Stokes lightning method. In that work, the analytic Goursat functions are expressed in a suitable basis with unknown coefficients, such as a composition of a Laurent series and a rational basis, leading to a severely ill-conditioned system matrix upon collocation at the boundary of the fluid domain. With row and column weighting, highly accurate results can despite the ill-conditioning be achieved for multiply connected domains.
% \end{remark}
% VERY small effect. not worth mentioning

%%%%%%%%%%%%%%%%%%%%%%%%%%
\section{Numerical results}\label{sec:num}
We now
demonstrate numerically that the proposed peanut-compressed two-body preconditioner dramatically accelerates
% ***AHB note you need to sell it more!                    ^             ^
and stabilizes the solve even for gaps on the order of $\delta = 10^{-3}$, for unit-radius disks.

Two bodies $i$ and $j$ are considered a near contact if $\|\vec c^{(i)}-\vec c^{(j)}\|\leq 2+\delta_{\text{c}}$, where we make the choice $\delta_{\text{c}} = 0.2$\footnote{The peanut self-intersects for separation $D \ge \sqrt{3}-1$;
  and $\delta_{\text{c}}$ should therefore be chosen certainly closer than that.}. The coarse sources are sampled uniformly on a circle of radius $R_c$, with $R_c = 1-(1/N_c)\log(\epsilon^{-1})$ \cite[Alg.~1]{Stein2022}, with $\epsilon = 10^{-12}$. The fine proxy radius $R_f$ is set analogously, with $N_{fp}$ replacing $N_c$. % *** Anna: not sure where to report common parameters so that they don't get in the way. 
For fast flow evaluation, we use the Stokes FMM in the FMM2D library \cite{fmm2d} with tolerance set to $10^{-9}$.

\begin{remark}[Code availability and hardware]\label{r:cpu} % rrrrrrrrrrrrrrrrrrrrrrrrrrrrrrrrrrrrrrrrrrr
  MATLAB code implementing the proposed methods
  is available at \url{https://github.com/annabroms/StokesMFS2D}. All timed examples were performed on a single AMD Genoa node of the Rusty cluster at the Flatiron Institute, using up to 96 CPU cores.
\end{remark}

%   There is little downside to $M_p$ being too large since they are discarded later. 
%%%sssssssssssssssssssssssssssss
\subsection{Resistance and mobility: two-body vs.~one-body preconditioning}\label{mob_res}
We compare the proposed scheme---peanut-compressed
% ***AHB needed, right? Or we need an acronym PC2BP. It's unclear that "two-body" is supposed to mean "peanut-compressed two-body" with Sec 3 as written.
  two-body preconditioning---against the one-body preconditioner, for geometries where small inter-particle distances can be controlled systematically.
With one-body preconditioning, to resolve all near-contacts the source and collocation points are formed by taking the union of the fine discretizations associated with all near contact pairs.
% Anna: Is the above better to understand what's going on? AHB:yes
  We generate random clusters so that each disk has separation exactly $\delta$ to at least one neighbor, and vary $\delta$ while fixing $P$ (the number of disks). For the resistance problem, we take $P = 20$, while for mobility, $P = 50$. The solves use input quantities (forces/torques or translational/angular velocities) sampled from a standard normal
% AHB you mean normal? (there is no "standard uniform distn"!). Anna: Yes, thanks.
  distribution, with forces subsequently shifted
  %***AHB "rescaling" would multiply only. But, you mean a constant force added to each disk? Anna: yes, is "shifted" a better word?
  so that the total force in the system sums to zero.

% *** Anna: moved here from Sec 5. Needed in Sec 5 if we want to display any results there.
To assess the accuracy of the solution, we define the \emph{pointwise residual}
\begin{equation}\label{residual_x}
 \epsilon_{\text{res}}(\vec x)\coloneqq \|\vec u(\vec x)-\vec g^{(i)}(\vec x)\|_{2}%{\|\vec g^{(i)}(\vec x)\|_{2}}
 ,\quad \vec x\in\partial\Omega^{(i)}.
\end{equation}
For a mobility problem, this measures the local discrepancy between the computed boundary velocity $\vec u(\vec x)$ and the no-slip boundary data $\vec g^{(i)}(\vec x)$ determined by the computed velocity pair $\vec v^{(i)},\omega^{(i)}$.  
We primarily report the \emph{relative boundary residual}, defined as
\begin{equation}\label{residual}
  \epsilon_{\text{res}}^{\max}\coloneqq \max\limits_{\vec x\in \partial\Omega} \epsilon_{\text{res}}(\vec x) /
  % ***AHB why not use the simpler max over $\pO$ here too:? Anna: you mean \max \|\vec g(\vec x)\|? What is $\vec g(\vec x)$?
  \left(\max\limits_{i = 1,\dots,P} \|\vec g^{(i)}(\vec x)\|\big|_{\vec x\in\partial\Omega^{(i)}}\right).
\end{equation}

\begin{figure}[t] % fffffffffffffffffffffffffffffffffffffffffffffffffffffffffffffffffffffffff
    \centering
     \begin{subfigure}[t!]{0.24\textwidth}
        \includegraphics[trim={0cm 0cm 0cm 0cm},clip,width = \textwidth]{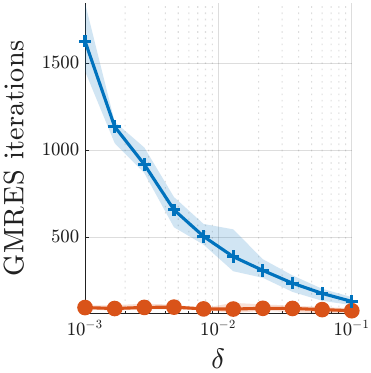}
        \caption{Resistance iterations}
        \label{fig:compare_iters_res}
      \end{subfigure}
    \begin{subfigure}[t!]{0.24\textwidth}
        \includegraphics[trim={0cm 0cm 0cm 0cm},clip,width = \textwidth]{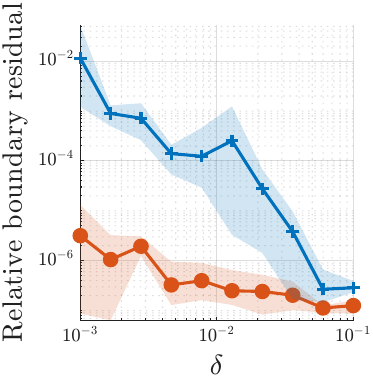}
        \caption{Resistance accuracy}
        \label{fig:cluster_err_res}
    \end{subfigure}
    \begin{subfigure}[t!]{0.24\textwidth}
        \includegraphics[trim={0cm 0cm 0cm 0cm},clip,width = \textwidth]{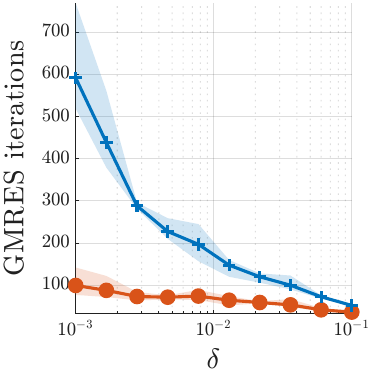}
        \caption{Mobility iterations}
        \label{fig:compare_iters_mob}
      \end{subfigure}
    \begin{subfigure}[t!]{0.24\textwidth}
        \includegraphics[trim={0cm 0cm 0cm 0cm},clip,width = \textwidth]{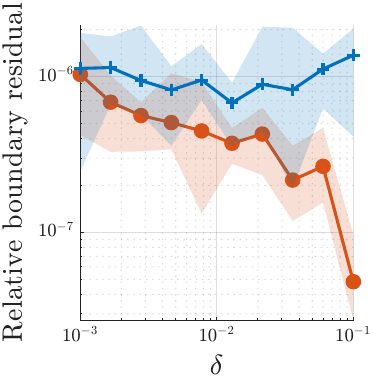}
        \caption{Mobility accuracy}
        \label{fig:cluster_err_mob}
    \end{subfigure}
    \caption{Comparison of proposed peanut-compressed two-body preconditioning method
      (red circles), against one-body preconditioning (blue + signs), for resistance and
      mobility BVPs.
      % ***AHB note clearer.
      The geometries are random clusters of unit disks with each disk having separation $\delta$ from at least one neighbor.
      % ***AHB: save this for the text (it's visible!):
      %A dramatic reduction in iteration count is observed, especially for small $\delta$.
      Curves show the mean (and shaded color the range) over five randomized runs at each $\delta$.}
   \label{fig:compare_iters}
\end{figure}

With one-body preconditioning, the iteration count grows rapidly as $\delta$ decreases in both settings, whereas it remains low with the two-body basis; see Figure~\ref{fig:compare_iters_res} for the resistance results and~\ref{fig:compare_iters_mob} for mobility. Resistance is harder to resolve, and we use $N_c = 150$, $N_{fp} = 120$ and $N_{\text{ell}} = 100$. In this setting, the proposed two-body solver achieves up to 4 digits more accuracy than the one-body solver; see Figure~\ref{fig:cluster_err_res}. For the mobility solves, we keep $N_{\text{ell}} = 100$, but use $N_c = 80$ and $N_{fp} = 60$,
%***AHB not N_p, right?  But it's weird that Nc > Nfp in both cases... really?
obtaining relative boundary residuals below $10^{-6}$; see Figure~\ref{fig:cluster_err_mob}.
Two-body preconditioning leads to a substantial reduction in the number of unknowns: at small $\delta$, the proposed mobility solver achieves an average 12.5x reduction while resolving a mean of 2.16 near contacts per body.
This combines with the iteration count reduction to give a speed-up in solve time of a factor of 30. For the resistance case, the speed-up for $\delta = 10^{-3}$ is a factor of 65.

We have also investigated the two-way error, where the output of, for example, the resistance problem is used as input to the mobility problem, and the result is compared with the original resistance input. As in earlier work \cite{Broms2024c}, we found that the error is small and proportional to the relative boundary residual; we do not discuss it further.

\begin{remark}[Increasing $N_c$] %Anna : perhaps remove?
The underlying recompression
assumption is that the two-body correction $\vec\eta^{(i,k)}[\mm^{(i)}](\x)$
  is smooth enough outside of the peanut for the $i\text{-}k$ pair to be accurately represented there by
  the union of the coarse sources for the contacting pair. 
  %This may not hold; 
  %the first this is to plot $\eta$ to test smoothness.
  Resistance problems require an %slight 
  increase in coarse nodes to get the residuals
  small in the two least-squares problems \eqref{sys2b} and \eqref{peanut}, explaining the larger $N_c$.  %It might fail because the union of two coarse
  %sources does not ``connect'' the middle part of the peanut enough (see picture).
  %***AHB true - but this is future work.
  \end{remark}
  
%%sssssssssssssssssssssssssssssssssssss
\subsection{Mobility for random disk packings}
We next demonstrate the proposed peanut-compressed two-body preconditioned
mobility solver on large random disk packings, experiments that would be infeasible using one-body preconditioned MFS or
standard boundary-based iterative methods.
We will study the effect of area fraction upon iteration count.

%***AHB this needed to be up-front...
We use a geometry generation method inspired by \cite{Helsing1995}, based on \cite{Hansen1986,Metropolis1953}. To generate the particle configuration, we first place $P$ unit disks in a square domain whose side length is chosen to match the prescribed area fraction $\varphi$. The initial configuration is taken to be a square lattice satisfying the non-overlap constraint. We then randomize this configuration by performing several sweeps over the particles. In each sweep, the disks are visited sequentially; for each disk, a random trial displacement is proposed and accepted only if the resulting configuration remains admissible, i.e., the disk stays within the domain and maintains the prescribed minimum separation from all other disks. Repeating these accept/reject moves produces a disordered packing at the desired area fraction while preserving a controlled minimum inter-particle distance.

We first discuss the large-scale $P=10{,}000$ experiment shown in Figure~\ref{large_ex}: the solve using a single core takes 9.6\,min, is 95.7\% FMM-dominated, and requires 72.4\,GiB of RAM
(120\,GiB if fine factorizations are kept for near-particle flow evaluation).
%Precompute time 23 min. 120GiB RAM if fine compressions are stored for velocity evaluation on the boundary
Using all 96 cores, the solve time drops to 36\,s, with a precomputation time of 2.9\,min. %larger RAM!
Here, parallel solve-time speed-up is limited by the use of MATLAB's single-threaded sparse matrix-vector multiply.
%AHB quote the parallel efficiency and say mostly due to FMM2D ... or SVDs??? The speedup in solve time is due to FMM2D and limited by a factor of 16 due to the application of the sparse correction maps. % Anna: ..why? The speedup in precompute time is obtained by batching the SVDs in parallel. 
In the first panel of Figure~\ref{large_ex}, particles are colored by a measure of overall boundary speed,
\begin{equation}\label{speed}
  s^{(i)} \coloneqq \|\vec v^{(i)}\|_2 + |\omega^{(i)}|,
  % ***AHB surely this is a absval of omega? Otherwise s<=0 is possible. I corrected.
\end{equation}
where $\vec v^{(i)}$ and $\omega^{(i)}$ denote the computed translational and angular velocity of particle $i$, obtained with $N_c = 60$. The second panel shows the maximum relative residual for each body, as defined in \eqref{residual_x}--\eqref{residual}. The rightmost panel shows the relative rigid-body error, computed against a reference solution obtained on a finer coarse grid with $N_c = 120$. Specifically, the absolute error for body $i$ is defined by
\[
e^{(i)} =
\left\|
\vec v^{(i)} - \vec v_{\mathrm{ref}}^{(i)}
\right\|_2
+
\left|
\omega^{(i)} - \omega_{\mathrm{ref}}^{(i)}
\right|,
\]
with corresponding relative error
\begin{equation}\label{rel_vel_err}
E^{(i)} = e^{(i)}/\max\limits_{j=1,\dots,P} s^{(j)}_{\mathrm{ref}}.
\end{equation}

\begin{figure}[t] % fffffffffffffffffffffffffffffffffffffffffffffffffffffffffffffffffff
    \centering
    \begin{subfigure}[b!]{0.24\textwidth}
        \includegraphics[trim={0cm 0cm 0cm 0cm},clip,width = 1.02\textwidth]{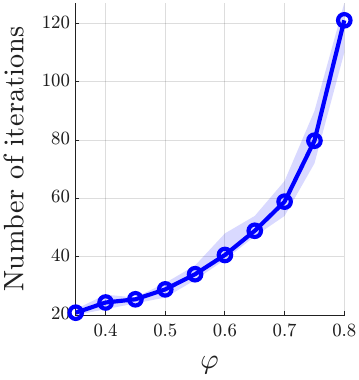}
        \caption{}
        \label{random_iters}
      \end{subfigure}
    % \begin{subfigure}[t!]{0.32\textwidth}       \includegraphics[trim={0cm 0cm 0cm 0cm},clip,width = \textwidth]{figures/phi_sweep_res.pdf}
    %     \caption{}
    %     \label{random_residual}
    % \end{subfigure}
       \begin{subfigure}[b!]{0.24\textwidth}
      \includegraphics[trim={0cm 0cm 0cm 0cm},clip,width = \textwidth]{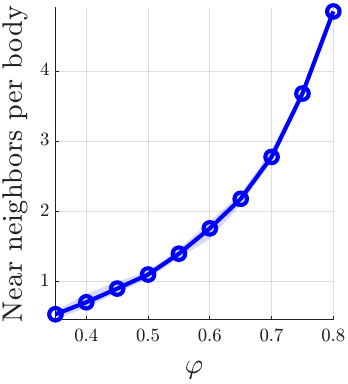}
        \caption{}
        \label{random_neigh}
    \end{subfigure}
    \begin{subfigure}[b!]{0.24\textwidth}
        \includegraphics[trim={0cm 0cm 0cm 0cm},clip,width = \textwidth]{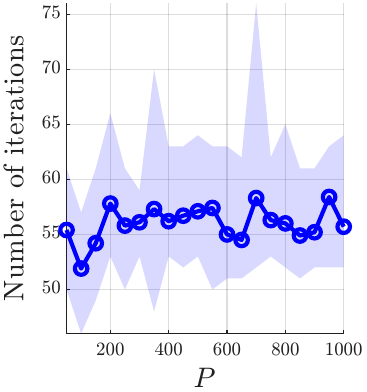}
        \caption{}
        \label{p_random_iters}
      \end{subfigure}
    % \begin{subfigure}[t!]{0.32\textwidth}       \includegraphics[trim={0cm 0cm 0cm 0cm},clip,width = \textwidth]{figures/p_sweep_res.pdf}
    %     \caption{}
    %     \label{p_random_residual}
    % \end{subfigure}
    \begin{subfigure}[b!]{0.24\textwidth}    
      \includegraphics[trim={0cm 0cm 0cm 0cm},clip,width = 1.02\textwidth]{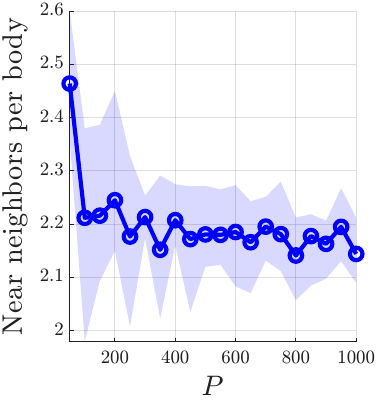}
        \caption{}
        \label{p_random_neigh}
           \end{subfigure}
    \caption{The mobility problem is solved for a random packing of $P$ particles at packing density $\phi$. In panels (a) and (b), $P = 500$ is fixed while $\phi$ is varied; in panels (c) and (d), $\varphi = 0.65$ is fixed while $P$ is varied.
}
\label{random_pack}
\end{figure}

Next, in Figure~\ref{random_pack}, we solve the same problem for varying $P$ and $\phi$, with force and torque entries sampled from a standard normal
%AHB normal? see above.
distribution. We first fix $P = 500$ and vary the packing density $\phi$, reporting the iteration count in panel~\ref{random_iters} and the average number of close neighbors per body in panel~\ref{random_neigh}. As $\phi$ increases, so does the iteration count. We then fix $\phi$ and vary $P$ in Figures~\ref{p_random_iters} and~\ref{p_random_neigh}: they suggest that the iteration count depends primarily on the number of near neighbors per body and remains essentially independent of $P$. All panels display the mean, minimum, and maximum over ten runs for each $(P,\phi)$ combination. In every case, the relative residual evaluated at newly sampled boundary nodes remains below $10^{-5}$.

For these computations, the pair problems were solved using $N_{fp} = 150$, $N_{\text{ell}} = 80$ and $M_p = 200$. Truncation levels for the SVDs were set to $10^{-14}$ ($10^{-11}$ when solving for fine sources) and in the global problem, the GMRES tolerance was $10^{-7}$.

%%%fffffffffffffffffffffffffffffffffffffffffff

\begin{remark}[Clustering]
We observe iteration counts independent of $P$ provided particles do not form tightly clustered configurations with multiple simultaneous near contacts. When close triangles occur, with all pairwise separations small, the iteration count increases by approximately $\mathcal O(1)$ per triangle. This mechanism explains the growth in iteration count observed in Figure~\ref{fig:compare_iters_mob} for the smallest values of $\delta$, as well as in Figure~\ref{random_iters}. Eigenvalue analysis of the preconditioned system matrices shows that each such triangle introduces $\mathcal O(1)$ eigenvalues near zero.
\end{remark}

\section{Conclusions}\label{sec:conc}
%*** Could you help me out here?
We introduce a two-body preconditioned method of fundamental solutions for the 2D Stokes resistance and mobility problems in close-to-touching suspensions of circular particles. Starting from a local basis constructed for each particle in isolation, a new two-body basis is built that incorporates corrections obtained from local high-resolution two-body boundary value problems. These corrections resolve the lubrication-driven fine scales that arise near close contacts while simultaneously regularizing the ill-conditioning caused by shrinking particle gaps.

A central feature of the scheme is the compression of the local fine-grid pair representation into an equivalent coarse representation on a geometrically natural ``peanut'' separation surface. The geometry of the particles therefore directly dictates the compression of the fine grid into a coarse one, allowing the global discretization to remain uniformly coarse while still resolving particle separations down to a thousandth of the radius. As a result, the solve is dominated by all-to-all evaluation of coarse sources at coarse collocation points, accelerated in our implementation by a linearly scaling fast multipole method.

Compared to standard block-diagonal right preconditioning, the two-body preconditioner not only dramatically reduces the GMRES iteration count, but also stabilizes the solve sufficiently to achieve higher accuracy at significantly lower computational cost. For the mobility problem, we observe iteration counts essentially independent of the number of bodies, $P$, depending primarily on the local near-neighbor count per particle. Although not shown explicitly in Section~\ref{sec:num}, the iteration count for resistance grows with $P$, consistent with earlier observations in \cite{USABIAGA2016,Quaife2018,Broms2024c}. An important direction for future work is therefore to couple the present approach with long-range preconditioning strategies \cite{Helsing2011,Helsing2011b} that address ill-conditioning arising from far-field effects.

Our framework is not specific to Stokes flow or to the method of fundamental solutions. We believe that the idea of local pair corrections extends naturally to other elliptic PDEs and to other boundary-based solvers, provided an accurate technique is available for solving the local two-body boundary value problems. Precomputing corrections for all near-contact pairs is currently the dominant cost of the scheme, and an important direction for future work is the development of strategies to reduce or amortize this setup cost. % *** Anna: How do we want to expand on this?

Future directions also include extending the method to more general geometries, such as non-circular particles and polydisperse suspensions, as well as to three-dimensional problems, building on our earlier work for spheres \cite{Broms2025,Broms2024c}. Since the global solve remains coarse while all singular near-contact interactions are resolved locally, the proposed scheme is fast, accurate, and robust, providing a scalable and flexible foundation for large-scale simulations of dense suspensions.

\section*{Acknowledgments}
Broms and Tornberg acknowledge support from the Swedish Research Council: grant no.~2023-04269. 
The Flatiron Institute is a division of the Simons Foundation. Broms gratefully acknowledges support from the Flatiron Institute through two research visits to the Center for Computational Mathematics during this work. All authors participated in the Fall 2025 program \emph{Interfaces and Unfitted Discretization Methods} at Institut Mittag-Leffler. We also benefited greatly from discussions with Daan Huybrechs, Daniel Fortunato, Dhairya Malhotra, and Leslie Greengard.

%\appendix
%\section{The enhanced discretization}
%\alert{Here we show something about the proper choice of fundamental solutions for close contacts.}

%	\bibliographystyle{myIEEEtran} %including doi

\bibliographystyle{abbrv} %including doi
\bibliography{MFS_refs}

\end{document}